\documentclass[10pt,a4paper]{article}
\usepackage[margin=1in]{geometry}
\usepackage{lmodern}
\usepackage{amsmath,mathtools,amssymb,amsthm,booktabs,array}
\usepackage{xcolor}
\definecolor{linkc}{rgb}{0,0.2,0.55}
\usepackage{enumitem}
\usepackage{graphicx}
\usepackage{float}
\usepackage{listings}
\usepackage{tikz}
\usetikzlibrary{arrows.meta,positioning,intersections}

\usepackage[colorlinks=true,linkcolor=linkc,citecolor=linkc,urlcolor=linkc]{hyperref}
\usepackage{aliascnt}
\newtheorem{theorem}{Theorem}[section]
\newaliascnt{proposition}{theorem}
\newtheorem{proposition}[proposition]{Proposition}
\aliascntresetthe{proposition}
\newaliascnt{lemma}{theorem}
\newtheorem{lemma}[lemma]{Lemma}
\aliascntresetthe{lemma}
\newaliascnt{corollary}{theorem}
\newtheorem{corollary}[corollary]{Corollary}
\aliascntresetthe{corollary}
\theoremstyle{remark}
\newaliascnt{remark}{theorem}
\newtheorem{remark}[remark]{Remark}
\aliascntresetthe{remark}
\theoremstyle{definition}
\newaliascnt{definition}{theorem}
\newtheorem{definition}[definition]{Definition}
\aliascntresetthe{definition}

\numberwithin{equation}{section}

\date{}

\newcommand{\cst}[1]{\ensuremath{\csname const@#1\endcsname}}
\newcommand{\registerconstant}[2]{%
  \expandafter\gdef\csname const@#1\endcsname{#2}}
\registerconstant{admissible}{\mathcal O}
\registerconstant{bed}{\mathcal B}
\registerconstant{foldCoefficient}{\Theta}
\registerconstant{interface}{\mathcal I}
\registerconstant{interfaceSpeed}{W_{\!I}}
\registerconstant{kernelThreshold}{W_{\!c}}
\registerconstant{surface}{\mathcal S}
\registerconstant{vorticityOne}{\Omega_1}
\registerconstant{vorticityTwo}{\Omega_2}

\begin{document}
\title{Closed bifurcation loops for periodic gravity waves with two layers of constant vorticity}
\author{Yanlong Fan$^{a}$\qquad Yiqiu Mao$^{b}$\qquad Quan Wang$^{c}$\thanks{Corresponding author.}\\[4pt]
\small $^{a}$School of Mathematics, Southwest Jiaotong University, Chengdu, Sichuan, 611756, China\\
\small $^{b}$School of Mathematics and Information Science, Guangzhou University, Guangzhou, 510000, China\\
\small $^{c}$College of Mathematics, Sichuan University, Chengdu, Sichuan, 610065, China}
\maketitle

\begin{abstract}
We study the bifurcation of periodic gravity waves in two layers of equal density and different constant vorticities. The free surface and the material interface are graphs, and the mean area of each layer is prescribed. First, we consider a collision of two positive fundamental dispersion roots for which the relative surface and interface velocities satisfy $s<0<W_I$. On the two-root side near the collision, we prove that the nearby nonlaminar waves and the two laminar endpoints form a simple closed analytic curve in the even, phase-fixed space. To leading order this curve is an ellipse in the flux parameter and the first cosine coefficient of the interface. We determine the sign of the cubic term by an exact polynomial inequality on the physical fold domain. Second, if both layers have a critical level, we show that each nonlaminar wave on the loop has two centres and two saddles per period. Finally, a nodal argument identifies the endpoint whenever the global graph branch returns to a laminar flow. The initial bending at the larger root can change sign at least twice before the roots coalesce.
\end{abstract}
\section{Introduction}\label{sec_intro}
In this paper we study periodic travelling gravity waves in a fluid of constant density. A material streamline separates two layers of different constant vorticities. The velocity is continuous across the interface, and both the free surface and the interface are unknown. We prescribe the wavelength and the mean area of each layer. For certain values of the shear, two distinct laminar flows admit a bifurcation in the same fundamental mode. We ask whether, with the shear held fixed, a branch issuing from one laminar flow reaches the other. We prove that this connection exists near a nondegenerate collision of fundamental roots with positive vertical profile, provided the relative horizontal velocity is negative at the surface and positive at the interface. In the space of even graph waves, the two laminar states and the waves between them form a simple closed analytic curve.

The existence theory for rotational waves provides a starting point for this question. Constantin and Strauss established a large-amplitude existence theory and subsequently considered bounded discontinuous vorticity \cite{CS,CS11}. By means of spatial dynamics, Groves and Wahl\'en obtained small-amplitude Stokes and solitary gravity waves \cite{GW08}. For local results on capillary--gravity waves with general or piecewise constant vorticity, see \cite{W06,MM}. Of particular relevance here is the work of Martin and Matioc \cite{MM16} on two layers of different constant vorticities. They constructed small-amplitude periodic gravity waves allowing stagnation and described the streamlines near bifurcation. We compare their surface-speed dispersion relation with ours in \autoref{rem_fourier_normalization}. At two distinct simple roots the local theory yields one arc near each root. Whether these arcs belong to the same branch requires information on the waves between the two laminar states.

A related question concerns the critical layers along a closed branch. Wahl\'en constructed waves with a critical layer \cite{Wahlen09}. Ehrnstr\"om, Escher and Wahl\'en treated multiple critical layers, and Ehrnstr\"om, Escher and Villari studied the corresponding interior dynamics \cite{EEW,EEV}. For analytic vorticity, Varholm obtained a global bifurcation result allowing multiple critical layers \cite{Varholm}. In our problem the vorticity jumps across a moving streamline. If the limiting shear has a critical height in each layer, we show that each layer contains one centre and one saddle per period for every nonlaminar wave on the small loop.

The first difficulty is to determine the dimension of the kernel at the collision. Several Fourier modes can participate in bifurcation for rotational waves, as in the trimodal and $N$-modal constructions of \cite{EW15,KL18}. Here both roots belong to the fundamental mode. We first determine the interval on which the normalized vertical profile is positive. On this interval the scalar dispersion equation is strictly convex in the square-root flux coordinate. Comparing the profiles across the interface by means of their Wronskian excludes all higher modes. Thus the kernel remains one-dimensional at the parameter where simple-root transversality fails.

The second difficulty is to determine the sign of the nonlinear term. Let $\lambda$ denote the flux parameter, $\mu$ the physical parameter being varied and $\tau$ the first cosine coefficient of the interface. Translation by half a period changes the sign of $\tau$. Solving the range equation and using this symmetry, we obtain $\tau G(\lambda,\mu,\tau^2)=0$. If $2\varepsilon$ denotes the leading separation of the two roots, the equation for nonlaminar waves takes the form
\[
 (\lambda-\lambda^\ast)^2+\Theta\tau^2=\varepsilon^2+O(\varepsilon^3).
\]
Here the sign of $\Theta$ determines whether the two laminar states lie on a small oval. To calculate this coefficient we must take into account both the mean flow and the second harmonic forced at second order. Each contributes to the fundamental mode at third order. We obtain a Stokes solvability determinant whose sign, on the physical fold manifold, reduces to a polynomial inequality in four variables. An exact integer certificate gives $\Theta>0$ throughout the admissible fold domain. Consequently the squared amplitude is strictly concave in $\lambda$ near the fold. Its two simple zeros are the laminar endpoints. The positive and negative amplitude arcs meet there, and the reconstruction accounts for all nearby nonlaminar waves.

Global bifurcation gives continuations with critical layers \cite{CSV,WW}; the fate of a particular branch depends on the boundary conditions and on which physical quantities remain fixed. Doak, Matthies, Sewell and Wheeler studied a vorticity interface between rigid walls \cite{DMSW}. In the present problem a free upper boundary supplies a Bernoulli condition, and the Bernoulli constant varies along the branch. Kozlov's loop-exclusion criteria illustrate the importance of the parameter choice, in particular when the period itself varies \cite{Ko23}. For two vortical layers, the computations of Sh\=oji and Okamoto exhibit capillary--gravity branches approaching another laminar state or a branch of different mode \cite[Figures~8--13]{SO}. Such diagrams motivate an endpoint question that a turning point in a numerical projection cannot decide.

At fixed period and fixed layer areas we prove that the connection exists near the fold. Away from the fold we obtain a conditional result. Under the continuation hypotheses, if the graph branch first returns to a laminar state, a nodal argument identifies that state. We also derive a signed flux bound independent of the surface slope. Even for a fixed pair of depths the initial bending at the larger root changes sign at least twice as the lower shear varies. Thus its direction alone does not determine whether the branch returns.

\subsection{The free-boundary problem}\label{sub_01_formation_of_the_problem}

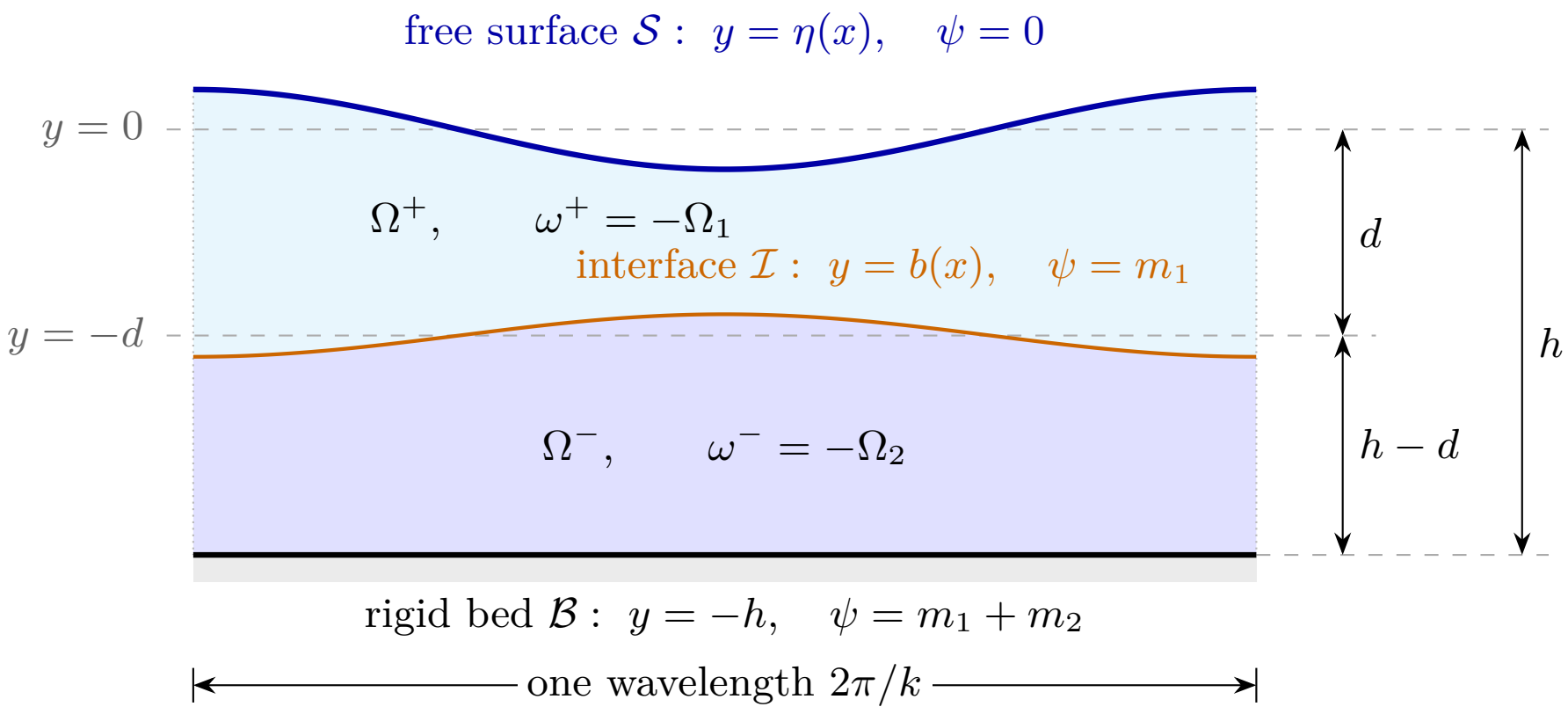
\begin{figure}[!htbp]
\centering
\begin{tikzpicture}[x=1cm,y=1cm,>=Stealth,font=\small]
\fill[cyan!8]
  plot[domain=0:8,samples=100] (\x,{0.30*cos(45*\x)})
  -- plot[domain=8:0,samples=100] (\x,{-1.55-0.16*cos(45*\x)}) -- cycle;
\fill[blue!12]
  plot[domain=0:8,samples=100] (\x,{-1.55-0.16*cos(45*\x)})
  -- (8,-3.2) -- (0,-3.2) -- cycle;
\fill[black!8] (0,-3.2) rectangle (8,-3.4);
\draw[gray!65,dashed] (-0.2,0) -- (10.2,0);
\draw[gray!65,dashed] (-0.2,-1.55) -- (8.9,-1.55);
\draw[gray!65,dashed] (8,-3.2) -- (10.2,-3.2);
\draw[gray!55,densely dotted] (0,0.30) -- (0,-3.2);
\draw[gray!55,densely dotted] (8,0.30) -- (8,-3.2);
\draw[very thick,blue!65!black]
  plot[domain=0:8,samples=100] (\x,{0.30*cos(45*\x)});
\draw[thick,orange!80!black]
  plot[domain=0:8,samples=100] (\x,{-1.55-0.16*cos(45*\x)});
\draw[very thick] (0,-3.2) -- (8,-3.2);
\node[blue!65!black] at (4,0.72)
  {free surface ${\mathcal S}:\ y=\eta(x),\quad\psi=0$};
\node[orange!80!black] at (5.2,-1.04)
  {interface ${\mathcal I}:\ y=b(x),\quad\psi=m_1$};
\node at (4,-3.66)
  {rigid bed ${\mathcal B}:\ y=-h,\quad\psi=m_1+m_2$};
\node at (2.7,-0.67)
  {$\Omega^+,\qquad\omega^+=-{\Omega_1}$};
\node at (4,-2.38)
  {$\Omega^-,\qquad\omega^-=-{\Omega_2}$};
\node[anchor=east,text=gray!75!black] at (-0.25,0) {$y=0$};
\node[anchor=east,text=gray!75!black] at (-0.25,-1.55) {$y=-d$};
\draw[<->] (8.65,0) -- (8.65,-1.55) node[midway,right] {$d$};
\draw[<->] (8.65,-1.55) -- (8.65,-3.2) node[midway,right] {$h-d$};
\draw[<->] (10,0) -- (10,-3.2) node[midway,right] {$h$};
\draw[<->] (0,-4.18) -- (8,-4.18)
  node[midway,fill=white,inner sep=2pt] {one wavelength $2\pi/k$};
\draw (0,-4.05) -- (0,-4.31);
\draw (8,-4.05) -- (8,-4.31);
\end{tikzpicture}
\caption{The equal-density, two-vorticity-layer geometry over one wavelength. The dashed horizontal lines mark the fixed mean levels $y=0$ and $y=-d$; the mean upper and lower depths are $d$ and $h-d$. The profiles are schematic.}
\label{fig_intro_model}
\end{figure}

Let $g>0$, $h>0$ and $k>0$ denote gravity, the mean total depth and the wavenumber, and set \(L=2\pi/k\) and \(\mathbb S=\mathbb R/L\mathbb Z\). For graph profiles, the free surface, material interface and rigid bed are illustrated in \autoref{fig_intro_model} and defined by
\begin{equation}\label{eq_intro_geometry}
\begin{aligned}
&{\mathcal S}=\{(x,\eta(x)):x\in\mathbb S\},\qquad {\mathcal I}=\{(x,b(x)):x\in\mathbb S\},\\
&{\mathcal B}=\{(x,-h):x\in\mathbb S\},\qquad -h<b(x)<\eta(x),\\
&\Omega^+=\{(x,y):x\in\mathbb S,\ b(x)<y<\eta(x)\},\\
&\Omega^-=\{(x,y):x\in\mathbb S,\ -h<y<b(x)\}.
\end{aligned}
\end{equation}

The densities are equal, $\rho^+=\rho^-=\rho_0>0$, and the shear rates satisfy ${\Omega_1}\ne{\Omega_2}$. In a frame moving at the wave speed $c$, the stream functions $\psi^\pm$ determine the velocity and vorticity by
\begin{equation}\label{eq_intro_velocity}
\begin{aligned}
&u^\pm-c=\psi^\pm_y,\qquad v^\pm=-\psi^\pm_x,\\
&\omega^+=-\Delta\psi^+=-{\Omega_1},\qquad
\omega^-=-\Delta\psi^-=-{\Omega_2}.
\end{aligned}
\end{equation}

We write $\psi$ for the piecewise stream function and $[f]=f^+-f^-$ for the jump across ${\mathcal I}$. The free-boundary transmission problem is
\begin{equation}\label{eq_intro_model}
\begin{aligned}
&\Delta\psi^+={\Omega_1}&&\text{in }\Omega^+,\\
&\Delta\psi^-={\Omega_2}&&\text{in }\Omega^-,\\
&\psi=0,\qquad |\nabla\psi|^2+2gy=Q&&\text{on }{\mathcal S},\\
&\psi^+=\psi^-=m_1,\qquad [\nabla\psi]=0&&\text{on }{\mathcal I},\\
&\psi=m_1+m_2&&\text{on }{\mathcal B}.
\end{aligned}
\end{equation}
Here $m_1,m_2$ are the relative layer fluxes and $Q$ is the Bernoulli constant. We seek periodic profiles and stream functions that are even in $x$. The physical parameters $g,k,h,{\Omega_1},{\Omega_2}$ and the mean upper-layer thickness $d\in(0,h)$ are fixed, with
\begin{equation}\label{eq_intro_mean_levels}
\langle\eta\rangle=0,\qquad \langle b\rangle=-d,\qquad
\langle f\rangle:=\frac{k}{2\pi}\int_0^{2\pi/k}f(x)\,dx.
\end{equation}

In the moving frame the velocity is divergence free, so we may introduce a stream function. Since the vorticity is constant in each layer, the stream function satisfies a Poisson equation there. The kinematic conditions prescribe its traces. The two fluids have equal density and there is no vortex sheet, so the velocity is continuous at the interface. At the free surface we use Bernoulli's law with atmospheric pressure. We write $\Omega_\eta=\{(x,y):-h<y<\eta(x)\}$ and $V=-\psi_x$. Since the traces on $y=b(x)$ are constant,
\begin{equation}\label{eq_trace_reduction}
 \psi_x^\pm=-b'\psi_y^\pm,\qquad [\nabla\psi]=0\quad\Longleftrightarrow\quad[\psi_y]=0.
\end{equation}
The surface has the corresponding tangential identity with $b$ replaced by $\eta$. The areas per period are $Lh$ and $L(h-d)$ for the whole fluid and lower layer, respectively. We normalize the horizontal period to $2\pi$ when indexing Fourier modes. The half-period translation $\mathcal T_\pi$ preserves evenness and reverses the first cosine coefficient. The involution
\begin{equation}\label{eq_mirror_symmetry}
(\psi,{{\Omega_1}},{{\Omega_2}},m_1,m_2,Q)
\mapsto
(-\psi,-{{\Omega_1}},-{{\Omega_2}},-m_1,-m_2,Q).
\end{equation}
allows either sign of the laminar surface velocity to be represented.
\subsection{Laminar flows and positive modes}\label{sec_laminar}

For laminar flows we have \(\eta\equiv0\), \(b\equiv-d\). The interface relative velocity is $W_I:=W(-d)$, and the velocity profile is piecewise affine:
\begin{equation}\label{eq_w}
\begin{aligned}
W(y)&=W_I+{{\Omega_1}}(y+d) &&(-d<y<0),\\
W(y)&=W_I+{{\Omega_2}}(y+d) &&(-h<y<-d).
\end{aligned}
\end{equation}
Its values at the free surface and at the bed are
\begin{equation}\label{eq_s_wb}
s:=W(0)=W_I+{{\Omega_1}}d,\qquad
W_b:=W(-h)=W_I-{{\Omega_2}}(h-d).
\end{equation}

Since \(\psi_y=W\) and \(\psi=0\) at \(y=0\), the layer fluxes are the integrals of this profile,
\[
\begin{aligned}
m_1&=-\int_{-d}^{0}W(y)\,dy=-W_Id-\frac12{{\Omega_1}}d^2,\\
m_2&=-\int_{-h}^{-d}W(y)\,dy=-W_b(h-d)-\frac12{{\Omega_2}}(h-d)^2,
\end{aligned}
\]
and we write \(\Psi^{\pm}(y)\) for the laminar stream function:
\begin{equation}\label{eq_psi_laminar}
\begin{aligned}
\Psi^{+}(y) &=m_1 + W_I(y+d)+\frac12{{\Omega_1}}(y+d)^2, && -d<y<0,\\
\Psi^{-}(y) &= m_1 + W_I(y+d)+\frac12{{\Omega_2}}(y+d)^2, && -h<y<-d.
\end{aligned}
\end{equation}

On the free surface, \(y=0\) and \(\psi_y=s\), so \eqref{eq_intro_model} reduces to $
Q=s^2$.
For \(Q>0\), the two branches are
$s=\sqrt{Q}\qquad\text{and}\qquad s=-\sqrt{Q}.$
The branches meet at \(s=0\), where \(Q=0\) and \(\partial_s Q=0\).

On the branch with \(u(0)<c\), we use
\begin{equation}\label{eq_lambda_definition}
\lambda:=(c-u(0))^2=s^2>0,\qquad \sqrt\lambda=c-u(0),\qquad s=-\sqrt\lambda.
\end{equation}
Substituting \(s=-\sqrt\lambda\) into \eqref{eq_s_wb} and the formulas for the fluxes and \(Q\), we obtain
\begin{equation}\label{eq_lampar}
\begin{aligned}
&s(\lambda)=-\sqrt\lambda,\qquad {{W_{\!I}}}(\lambda)=-\sqrt\lambda-{{\Omega_1}}d,\\
& W_b(\lambda)=-\sqrt\lambda-{{\Omega_1}}d-{{\Omega_2}}(h-d),\qquad Q(\lambda)=\lambda,\\
&m_1(\lambda)=d\sqrt\lambda+\frac12{{\Omega_1}}d^2,\qquad m_2(\lambda)=(h-d)\big(\sqrt\lambda+{{\Omega_1}}d\big)+\frac12{{\Omega_2}}(h-d)^2 .
\end{aligned}
\end{equation}
For a nonlaminar wave only $m_2=m_2(\lambda)$ is prescribed in this way. We determine the upper flux $m_1$ and the Bernoulli constant $Q$ along the branch. The laminar quantities are analytic for $\lambda>0$. Differentiation gives
\[
\begin{aligned}
m_2'(\lambda)&=\frac{h-d}{2\sqrt\lambda}>0,\qquad
{{W_{\!I}}}'(\lambda)=s'(\lambda)&=-\frac1{2\sqrt\lambda}.
\end{aligned}
\]
With the physical parameters fixed, either \(\lambda\) or \({{W_{\!I}}}\) determines the laminar flow on this branch. For a quantity \(G\) on this family, write \(\partial_WG\) for its derivative with respect to \({{W_{\!I}}}\). The chain rule gives
\[
\partial_\lambda G={{W_{\!I}}}'(\lambda)\,\partial_WG=-\frac{\partial_WG}{2\sqrt\lambda}.
\]

\begin{lemma}\label{lem_cl}
For the laminar flow \eqref{eq_w}, \(W\) has a simple zero in \((-d,0)\) if and only if \(W_Is<0\). Its location is
$
y=-\frac{s}{{{\Omega_1}}}$.
It has a simple zero in \((-h,-d)\) if and only if \(W_IW_b<0\). Its location is $
y=-d-\Bigl|\frac{W_I}{{{\Omega_2}}}\Bigr|$.
If \(s<0<W_I\) and both layers contain a zero, then \({{\Omega_1}}<0<{{\Omega_2}}\).
\end{lemma}

For a laminar flow with a zero of \(W\) in \((-d,0)\), the mirror symmetry \eqref{eq_mirror_symmetry} allows us to choose \(s<0\). Then \autoref{lem_cl} and \eqref{eq_lampar} give
${{\Omega_1}}<0,\qquad s<0<{{W_{\!I}}},\qquad 0<\lambda<{{\Omega_1}}^2d^2.$
Both \((-d,0)\) and \((-h,-d)\) contain a zero of \(W\) if and only if
${{\Omega_2}}>0,\qquad {{W_{\!I}}}<{{\Omega_2}}(h-d).$

\begin{remark}\label{rem_lambdazero}
We set $H_-=h-d$ and $m_{20}=\Omega_1 dH_-+\Omega_2H_-^2/2$. In the unshifted physical equations the bed datum is $m_1+m_2$ and depends affinely on $m_2$. The signed laminar parameter $\sigma=(m_2-m_{20})/H_-$ is analytic on both sides of zero. Along the chosen laminar branch $\lambda=\sigma^2$ with $\sigma>0$. Thus $m_2=m_{20}$ is a coordinate endpoint of this parametrization. The laminar surface is stagnant there; continuation of a particular nonlaminar branch beyond that level would require the admissibility hypotheses in the enlarged parameter domain.
\end{remark}

Consider a laminar state with \(s\,{W_{\!I}}\ne0\). Linearizing in \autoref{sec_modes}, we obtain the following system for an even perturbation \(\varphi(y)\cos nx\), \(n\ge1\).
\begin{equation}\label{eq_intro_mode}
\begin{aligned}
&\varphi''=n^2\varphi
  &&\text{on }(-h,-d)\text{ and }(-d,0),\\
&\varphi(-h)=0,\\
&[\varphi]=0
  &&\text{at }y=-d,\\
&[\varphi']-\frac{{{\Omega_1}}-{{\Omega_2}}}{{{W_{\!I}}}}\varphi=0
  &&\text{at }y=-d.
\end{aligned}
\end{equation}
The normalization in \eqref{eq_f} fixes the scale by prescribing the lower-layer profile
$\varphi(y)=\sinh n(y+h),\qquad -h<y<-d .$
The free-surface condition for this mode is \(D_n(\lambda)=0\), where
\begin{equation}\label{eq_intro_dispersion}
D_n(\lambda):=s^2\varphi'(0)-(g+{\Omega_1}s)\varphi(0).
\end{equation}

We call a root \(\lambda\) of \(D_n\) \emph{positive} if its normalized profile satisfies \(\varphi>0\) on \((-h,0]\). We say that \(\lambda^\ast\) is a \emph{simple root} of \(D_1\) if
\begin{equation}\label{eq_intro_simple}
D_1(\lambda^\ast)=0,\qquad \partial_\lambda D_1(\lambda^\ast)\ne0.
\end{equation}

For \(D_1\), under \({\Omega_1}<\min(0,{\Omega_2})\), \autoref{prop_pairing} gives a positive-kernel window in the component \(s<0<{W_{\!I}}\) and counts the positive roots there. We write \({\Omega_2}^\ast\) for the threshold defined by the minimum criterion of \autoref{prop_pairing} at the chosen wavenumber, with \(k=1\) for \(D_1\). If the reduced scalar function has an interior minimum, the positive roots are
\[
\begin{cases}
\text{two simple roots}, & {\Omega_1}<{\Omega_2}<{\Omega_2}^\ast,\\
\text{one double root}, & {\Omega_2}={\Omega_2}^\ast,\\
\text{no roots}, & {\Omega_2}>{\Omega_2}^\ast .
\end{cases}
\]
If the reduced scalar function in \autoref{prop_pairing} has no interior minimum, then \({\Omega_2}^\ast={\Omega_1}\) and there are no positive roots in the component under the stated vorticity condition.
When two positive roots exist, write \(\lambda_-<\lambda_+\) for the pair.

At a simple root, the local bifurcation theorems also require laminar admissibility and higher-mode nonresonance. Their respective conditions are
\[
\lambda^\ast>0,\quad {W_{\!I}}(\lambda^\ast)\ne0,
\qquad\qquad
D_n(\lambda^\ast)\ne0\quad(n\ge2).
\]
At an admissible positive root of \(D_1\), higher-mode nonresonance follows from \autoref{lem_sl}.

For global graph continuation we assume
\begin{equation*}\tag{GC1}\label{eq_intro_continuation_signs}
{\Omega_1}<0<{\Omega_2}.
\end{equation*}
Let \(\varphi\) be the \(n=1\) profile at \(\lambda^\ast\), normalized by \eqref{eq_f}. The starting-root conditions are
\begin{equation*}\tag{GC2}\label{eq_intro_continuation_root}
\begin{aligned}
&D_1(\lambda^\ast)=0,\qquad \partial_\lambda D_1(\lambda^\ast)\ne0,\\
&\varphi>0\ \text{on }(-h,0],\qquad s(\lambda^\ast)<0<{W_{\!I}}(\lambda^\ast).
\end{aligned}
\end{equation*}
Under (GC2), \autoref{lem_sl} supplies higher-mode nonresonance. Let \(\mu\) be a varying physical parameter. We call an admissible, nonresonant double root \((\lambda^\ast,\mu_0)\) of \(D_1(\lambda;\mu)\) a \emph{nondegenerate fold} if
\[
\partial_\lambda^2 D_1(\lambda^\ast;\mu_0)\ne0,\qquad
\partial_\mu D_1(\lambda^\ast;\mu_0)\ne0 .
\]
See \autoref{def_fold}. The \emph{two-root side} of \(\mu_0\) is the side with two nearby simple roots \(\lambda_-<\lambda_+\).

On a local graph branch, let \(\tau=\zeta_1\) be the first cosine coefficient of the interface displacement \(b+d\). The normalization \eqref{eq_local_interface_amplitude} relates \(\tau\) to the local kernel coordinate \(a\). At an admissible, nonresonant simple root, let \(\lambda_2\) be the coefficient of \(\zeta_1^2\) in the expansion of \(\lambda\).

The \emph{fold coefficient} \({\Theta}\) in \eqref{eq_theta} is the common limit of \(\lambda_2^{(\pm)}(\lambda_\mp-\lambda_\pm)\) as \(\mu\to\mu_0\) on the two-root side, where \(\lambda_2^{(\pm)}\) is the coefficient on the branch through \(\lambda_\pm\) (\autoref{lem_theta}).

\subsection{Main results}\label{sub_01_main_results}

We first state the result on small closed branches. We then give the first-return result and the bounds on the physical quantities along the continuation.

\begin{theorem}\label{thm_b}
Let $(\lambda^\ast,\mu_0)$ be a nondegenerate fold at which $s<0<{W_{\!I}}$ and the normalized mode-one profile satisfies $\varphi>0$ on $(-h,0]$. For $\mu$ sufficiently close to $\mu_0$ on the two-root side, the nearby nontrivial solutions and the two laminar flows form a simple closed loop of graph waves. Under \eqref{eq_intro_continuation_signs} it is the image of the graph continuation through either root. Let $\varepsilon=\varepsilon(\mu)$ be the leading half-separation of the two roots given by \eqref{eq_eps}. In $(\lambda,\zeta_1)$ coordinates each loop lies near the ellipse
\begin{equation}\label{eq_intro_ellipse}
(\lambda-\lambda^\ast)^2+{{\Theta}}\zeta_1^2=\varepsilon^2,\qquad 2\varepsilon\sim C\,|\mu-\mu_0|^{1/2},
\end{equation}
where the asymptotic relation is taken as $\mu\to\mu_0$ on the two-root side.
Uniformly along the loop, the left-hand side of \eqref{eq_intro_ellipse} divided by $\varepsilon^2$ is $1+O(\varepsilon)$ by \autoref{thm_loops}. At fixed fold data its maximal interface coefficient satisfies
$\max|\zeta_1|=(\varepsilon/\sqrt\Theta)(1+O(\varepsilon^2))$.
Here and below, the term \emph{isola} refers to the nonlaminar waves together with the two laminar endpoints. The laminar line also belongs to the full local zero set.

Suppose also that the laminar flow at the fold has critical heights
$W(y_c^+)=W(y_c^-)=0,\qquad y_c^+\in(-d,0),\qquad y_c^-\in(-h,-d).$

\begin{samepage}
Each nontrivial solution has exactly four nondegenerate stagnation points per period. Their locations and types are listed below on the positively oriented half of the loop.
\begin{equation}
\begin{array}{@{}lcc@{}}
\text{layer}& x=0 &x=\pi \\
\Omega^+ & \text{centre }(0,Y_0^+) & \text{saddle }(\pi,Y_\pi^+)\\
\Omega^- & \text{saddle }(0,Y_0^-) & \text{centre }(\pi,Y_\pi^-).
\end{array}
\end{equation}
\par
\end{samepage}
Here $Y_0^\pm$ and $Y_\pi^\pm$ lie near $y_c^\pm$. Translating by half a period interchanges $x=0$ and $x=\pi$. In this way we obtain the stagnation points on the second half of the loop.
\end{theorem}

After identifying half-period translates, this circle becomes an interval whose endpoints are the two laminar flows.
\begin{theorem}\label{thm_c}
For all $g,k,{{\Omega_1}},h,d$, every fold arising from a double root of the fundamental-mode dispersion function \(D_k\) in \eqref{eq_d} with $s<0<{W_{\!I}}$ and normalized profile $\varphi>0$ on $(-h,0]$ satisfies $
{{\Theta}}>0$.
Hence as ${{\Omega_2}}$ decreases through ${{\Omega_2}}^\ast$, every pair of positive roots with an upper critical layer is born with a small loop joining its two members. If the laminar flow at the fold satisfies \(W(y_c^+)=W(y_c^-)=0\) with \(y_c^+\in(-d,0)\) and \(y_c^-\in(-h,-d)\), each nontrivial wave on the loop has the four stagnation points of \autoref{thm_b}.
\end{theorem}

\begin{theorem}\label{thm_a}
Assume \eqref{eq_intro_continuation_signs} and \eqref{eq_intro_continuation_root}. The graph branch issuing from $\lambda^\ast$ has positive orientation up to its first positive laminar return. In the event of a return, the branch reaches the other positive root. Its image then consists of the connecting arc, the half-period translate of that arc, and the two laminar endpoints. Otherwise, in each direction there is a sequence along which one of the following occurs: surface or interface stagnation is approached, a layer thickness tends to zero, a profile slope becomes unbounded, or $\lambda=0$. In either case $0<\sqrt\lambda<(-\Omega_1)d+\Omega_2(h-d)/2$.
\end{theorem}
\begin{theorem}\label{thm_d}
Let \(\lambda^\ast\) be an admissible simple root of \(D_1\) with \(D_n(\lambda^\ast)\ne0\) for every \(n\ge2\).
The local branch through the laminar flow at \(\lambda^\ast\) is a pitchfork, possibly degenerate, with expansion
$\lambda=\lambda^\ast+\lambda_2\zeta_1^2+O(\zeta_1^4).$
For $g=k=1$, \eqref{eq_lambda2formula} gives $\lambda_2$ at a simple positive root in $s<0<{W_{\!I}}$. For general $g,k>0$, the coefficient in the original variables is $gk$ times the coefficient for the scaled problem.
Let $\lambda_-<\lambda_+$ be the two positive roots of \autoref{prop_pairing} in this component. As the lower-layer shear ${\Omega_2}$ approaches ${\Omega_2}^\ast$ from below,
$\lambda_2^{(\pm)}(\lambda_\mp-\lambda_\pm)\to{{\Theta}}>0,$
hence \(\lambda_2^{(+)}<0<\lambda_2^{(-)}\) close to the fold. In this limit the $\mathbb Z_2$-symmetric isola of \autoref{thm_type} shrinks to the laminar flow at the fold. At ${\Omega_2}={\Omega_2}^\ast$ all nearby solutions are trivial by \autoref{thm_loops}.
At $g=k=1$ and $(\Omega_1,h,d)=(-1,4,3)$ one has $\lambda_2^{(+)}>0$ at $\Omega_2=1$. Along the larger-root branch its sign changes at least twice, once in $(0.2100,0.2102)$ and once in $(1.7115,1.7117)$.
\end{theorem}

The paper is organized as follows. In Section~2 we formulate the problem on fixed strips and study its linearization. Section~3 deals with the dispersion roots and local bifurcation. In Section~4 we reduce the wave equation to a scalar equation near a fold. Section~5 determines the sign of its nonlinear coefficient and gives four examples with interval bounds. In Section~6 we identify the first laminar return, derive flux bounds and present a numerical example.
\section{The graph problem and its physical linearization}\label{sec_operator_form}
\subsection{The equivalent problem on a fixed domain}\label{sec_flat}

With positive layer thicknesses, a vertical affine map identifies each moving layer with the fixed strip \(\mathbb S\times[0,1]\). We pull back the fluid and boundary equations to these strips.

\begin{equation}\label{eq_flat}
\begin{aligned}
\Omega^+&:\quad y=Y^+(x,t)=b(x)+tT^+(x),\qquad T^+=\eta-b,\\
\Omega^-&:\quad y=Y^-(x,t)=-h+tT^-(x),\qquad T^-=b+h,
\end{aligned}\qquad t\in[0,1],
\end{equation}
For \(Y=Y^\pm\), write \(\Pi_Y(x,t)=(x,Y(x,t))\) for the map from strip coordinates to physical coordinates. The stream functions on the fixed strips are
\begin{equation}\label{eq_strip_phi}
\Phi^\pm(x,t)=\psi^\pm\bigl(x,Y^\pm(x,t)\bigr).
\end{equation}

\begin{figure}[htbp]
\centering
\begin{tikzpicture}[>=Stealth,scale=1.0]

\def\xL{-4}
\def\xR{4}
\def\hval{1.3}

\def\bfun#1{0.32*sin(90*#1)+0.10*sin(200*#1)}
\def\etafun#1{\bfun{#1}+0.95+0.20*cos(70*#1)} 

\node at (0,1.95) {$(x,y)$};

\draw[thick,black] (\xL,-\hval) -- (\xR,-\hval)
    node[midway,below] {rigid bed $\mathcal B$, $y=-h$};

\fill[red!10] plot[domain=\xL:\xR,samples=200] (\x,{\etafun{\x}})
    -- plot[domain=\xR:\xL,samples=200] (\x,{\bfun{\x}}) -- cycle;
\fill[blue!10] plot[domain=\xL:\xR,samples=200] (\x,{\bfun{\x}})
    -- (\xR,-\hval) -- (\xL,-\hval) -- cycle;

\draw[thick,blue] plot[domain=\xL:\xR,samples=200] (\x,{\bfun{\x}});
\draw[thick,red] plot[domain=\xL:\xR,samples=200] (\x,{\etafun{\x}});

\node[red] at (-2.2,0.35) {$\Omega^+$};
\node[blue] at (-2.6,-0.8) {$\Omega^-$};
\node[red, anchor=west] at (\xR+0.25,1.05) {free surface $\mathcal S$, $y=\eta(x)$};
\node[blue, anchor=west] at (\xR+0.25,-0.15) {interface $\mathcal I$, $y=b(x)$};

\foreach \xpos in {-2.5,0,2.5}{
  \pgfmathsetmacro{\by}{\bfun{\xpos}}
  \pgfmathsetmacro{\ey}{\etafun{\xpos}}
  \draw[dashed] (\xpos,-\hval) -- (\xpos,\by);
  \draw[dashed] (\xpos,\by) -- (\xpos,\ey);
  \fill (\xpos,-\hval) circle (1.3pt);
  \fill (\xpos,\by) circle (1.3pt);
  \fill (\xpos,\ey) circle (1.3pt);
}

\draw[thick,->] (0,-2.0) -- (0,-3.2)
    node[midway,right] {flattening $\Pi_{Y^\pm}^{-1}$};

\node[above=0.3em] at (0,-4.4) {$(x,t)$};

\draw[thick,red] (\xL,-4.2) rectangle (\xR,-5.2);
\node[red] at (0,-4.7) {$\mathbb S\times[0,1],\;\Phi^+$};
\draw[red] (\xL,-4.2) node[left] {$t=1$} -- (\xR,-4.2) node[right] {$t=1\mapsto y=\eta(x)$};
\draw[red] (\xL,-5.2) node[left] {$t=0$} -- (\xR,-5.2) node[right] {$t=0\mapsto y=b(x)$};

\draw[thick,blue] (\xL,-5.2) rectangle (\xR,-6.4);
\node[blue] at (0,-5.8) {$\mathbb S\times[0,1],\;\Phi^-$};
\draw[blue] (\xL,-5.2) node[left] {$t=1$} -- (\xR,-5.2) node[right] {$t=1\mapsto y=b(x)$};
\draw[blue] (\xL,-6.4) node[left] {$t=0$} -- (\xR,-6.4) node[right] {$t=0\mapsto y=-h$};

\foreach \xpos in {-2.5,0,2.5}{
  \draw[red,dashed] (\xpos,-5.2) -- (\xpos,-4.2);
  \draw[blue,dashed] (\xpos,-6.4) -- (\xpos,-5.2);
  \fill (\xpos,-4.2) circle (1.3pt);
  \fill (\xpos,-5.2) circle (1.3pt);
  \fill (\xpos,-6.4) circle (1.3pt);
}

\node[below=1.2em of current bounding box]
{$
\begin{aligned}
Y^+(x,t)&=b(x)+tT^+(x),\quad T^+=\eta-b,\\
Y^-(x,t)&=-h+tT^-(x),\quad T^-=b+h.
\end{aligned}
$};
\end{tikzpicture}
\caption{Coordinate flattening for a two-layer fluid. The inverse maps \(\Pi_{Y^\pm}^{-1}\) send the closures of the physical layers \(\Omega^\pm\) to the fixed strips \(\mathbb S\times[0,1]\).}
\label{fig_flattening_map}
\end{figure}

By the chain rule, the physical derivatives of \(\psi^\pm\) at \(\Pi_{Y^\pm}(x,t)\) are as follows.
\begin{equation}\label{eq_strip_vel}
T^\pm=Y^\pm_t,\qquad
\psi^\pm_y=\frac{\Phi^\pm_t}{T^\pm},\qquad
\psi^\pm_x=\Phi^\pm_x-\frac{\Phi^\pm_t Y^\pm_x}{T^\pm}.
\end{equation}
A second differentiation yields the following expression for the Laplacian on the fixed strip
\begin{equation}\label{eq_lapflat}
\begin{aligned}
&(T^\pm)^2\,(\Delta\psi^\pm)\circ\Pi_{Y^\pm}
=\mathcal P_{Y^\pm}\Phi^\pm\\
&:=(T^\pm)^2\Phi^\pm_{xx}-2T^\pm Y^\pm_x\Phi^\pm_{xt}
+(1+(Y^\pm_x)^2)\Phi^\pm_{tt}
+(2Y^\pm_x T^\pm_x-T^\pm Y^\pm_{xx})\Phi^\pm_t.
\end{aligned}
\end{equation}

If \(T^\pm>0\), each \(\Pi_{Y^\pm}\) maps \(\mathbb S\times[0,1]\) bijectively onto \(\overline{\Omega^\pm}\). Using \eqref{eq_lapflat} for the Poisson equations and \eqref{eq_trace_reduction} for velocity continuity, we pull back \eqref{eq_intro_model} to the equivalent system \eqref{eq_fixed_1}--\eqref{eq_fixed_bed}.
\begin{align}
&\mathcal P_{Y^+}\Phi^+={{\Omega_1}}(T^+)^2
\quad\text{on }\mathbb S\times(0,1),\qquad
\mathcal P_{Y^-}\Phi^-={{\Omega_2}}(T^-)^2
\quad\text{on }\mathbb S\times(0,1),\label{eq_fixed_1}\\
&\Phi^+(\cdot,1)=0,\qquad
(1+\eta'^2)\Phi^+_t(\cdot,1)^2+(2g\eta-Q)(T^+)^2=0
\quad\text{on }{{\mathcal S}},\label{eq_fixed_surf}\\
&\Phi^+(\cdot,0)=\Phi^-(\cdot,1)=m_1,\qquad
\frac{\Phi^+_t(\cdot,0)}{T^+}-\frac{\Phi^-_t(\cdot,1)}{T^-}=0
\quad\text{on }{{\mathcal I}},\label{eq_fixed_int}\\
&\Phi^-(\cdot,0)=m_1+m_2
\quad\text{on }{{\mathcal B}}.\label{eq_fixed_bed}
\end{align}
The unknowns in \eqref{eq_fixed_1}--\eqref{eq_fixed_bed} are \((\eta,b,\Phi^\pm,m_1,Q)\).

\subsection{The nonlinear operator and admissible set}\label{sec_formulation}

We take the departure from the laminar flow as the unknown. In these variables the graph equation consists of two interior equations and two boundary conditions. The layer thicknesses and boundary speeds specify the open set on which we work.

We fix \(\lambda>0\), prescribe \(m_2=m_2(\lambda)\), and let \(\Psi_\lambda\) be the laminar stream function of \eqref{eq_psi_laminar}. On the two reference strips the laminar profiles are
\[
\Phi_0^+(t)=\Psi_\lambda(-d+td),\qquad
\Phi_0^-(t)=\Psi_\lambda(-h+t(h-d)).
\]
We use the deviations from the laminar state at the same \(\lambda\) as unknowns for \eqref{eq_fixed_1}--\eqref{eq_fixed_bed}. The reconstruction is
\begin{align}\label{un_know}
&\begin{aligned}
&v=(\eta,\zeta,\delta\varphi^+,\delta\varphi^-,a,q),\\
&b=-d+\zeta,\qquad m_1=m_1(\lambda)+a,\qquad Q=\lambda+q,
\end{aligned}\\
&\Phi^+=\Phi_0^++a(1-t)+\delta\varphi^+,\qquad
\Phi^-=\Phi_0^-+a+\delta\varphi^-.\label{eq_reconstruct_phi}
\end{align}
It follows from the endpoint values of \(\Phi_0^\pm\) that \(\delta\varphi^\pm\) vanish at \(t=0,1\) .

We use the complete state $w$ for estimates on the reconstructed wave. Its stream-function components are
\[
\varphi^+=\Phi^+-m_1(1-t),\qquad
\varphi^-=\Phi^--m_1-m_2(\lambda)(1-t).
\]
The complete state and its laminar value at $v=0$ are
\begin{equation}\label{eq_graph_states}
\begin{aligned}
&w=(\eta,\zeta,\varphi^+,\varphi^-,m_1,Q),\\
&w_0(\lambda)=(0,0,\varphi_0^+,\varphi_0^-,m_1(\lambda),\lambda),\\
&\varphi_0^+(t)=\Phi_0^+(t)-m_1(\lambda)(1-t),\\
&\varphi_0^-(t)=\Phi_0^-(t)-m_1(\lambda)-m_2(\lambda)(1-t).
\end{aligned}
\end{equation}
Thus \(\delta\varphi^\pm=\varphi^\pm-\varphi_0^\pm\). The map from perturbations to complete states is
$\mathcal A_{\rm g}(\lambda,v)=(\lambda,w),\qquad w=w_0(\lambda)+v.$
By the laminar formulas, \(w_0(\lambda)\) is polynomial in \(\sqrt\lambda\) and hence real-analytic for \(\lambda>0\).

The domain and target spaces are
\begin{equation}\label{eq_spaces}
\begin{aligned}
X&=C^{3,\alpha}_{e,0}(\mathbb S)\times C^{3,\alpha}_{e,0}(\mathbb S)
\times C^{3,\alpha}_{e,D}(\mathbb S\times[0,1])
\times C^{3,\alpha}_{e,D}(\mathbb S\times[0,1])\times\mathbb R^2,\\
Y&=C^{1,\alpha}_{e}(\mathbb S\times[0,1])
\times C^{1,\alpha}_{e}(\mathbb S\times[0,1])
\times C^{2,\alpha}_{e}(\mathbb S)
\times C^{2,\alpha}_{e}(\mathbb S),
\end{aligned}
\end{equation}
with subscripts: \(e\) for evenness in \(x\), \(0\) for zero mean, and \(D\) for vanishing boundary values at \(t=0,1\).

With the fields reconstructed from \(v\) by \eqref{un_know}--\eqref{eq_reconstruct_phi}, define \(F=(F_1^+,F_1^-,F_2,F_3)\) by
\begin{align}
&F_1^\pm(\lambda,v):=\mathcal P_{Y^\pm}\Phi^\pm
-\Omega_{1,2}(T^\pm)^2,\label{eq_f1}\\
&F_2(\lambda,v):=(1+\eta'^2)\Phi^+_t(\cdot,1)^2
+(2g\eta-Q)(T^+)^2,\label{eq_f2}\\
&F_3(\lambda,v):=\frac{\Phi^+_t(\cdot,0)}{T^+}
-\frac{\Phi^-_t(\cdot,1)}{T^-}.\label{eq_f3}
\end{align}

For each \(\delta>0\), evaluate the thicknesses and boundary speeds on the complete state \(w=w_0(\lambda)+v\). The open admissible sets are
\begin{equation}\label{eq_odelta}
\begin{aligned}
{\mathcal O}_\delta
&=\Bigl\{(\lambda,v)\in(0,\infty)\times X:\quad
  \lambda>\delta,\quad \inf T^+>\delta,\quad \inf T^->\delta,\\
&\hspace{1.5em} \inf_{\mathcal S}|\nabla\psi^+|>\delta,\quad
  \inf_{\mathcal I}|\nabla\psi^+|>\delta\Bigr\},\quad
{\mathcal O}=\bigcup_{\delta>0}{\mathcal O}_\delta.
\end{aligned}
\end{equation}
The gradient magnitudes on the surface and interface are given by the trace formulas
\[
\begin{aligned}
|\nabla\psi^+|^2\big|_{\mathcal S}
&=(1+\eta'^2)\frac{\Phi^+_t(\cdot,1)^2}{(T^+)^2},\qquad
|\nabla\psi^+|^2\big|_{\mathcal I}
=(1+b'^2)\frac{\Phi^+_t(\cdot,0)^2}{(T^+)^2}.
\end{aligned}
\]
These formulas give continuity of the speed margins on the set of positive thicknesses (cf.\ \eqref{eq_trace_reduction}).

\begin{proposition}\label{prop_operator}
For the operator \(F=(F_1^+,F_1^-,F_2,F_3)\), the following statements hold:
\begin{enumerate}[label=\textup{(\roman*)},leftmargin=2.3em]
\item \(F:{\mathcal O}\to Y\) is real-analytic, and \(F(\lambda,0)=0\) whenever \((\lambda,0)\in{\mathcal O}\).
The precise laminar admissibility conditions are
\[
\begin{aligned}
&(\lambda,0)\in{\mathcal O}
\quad\Longleftrightarrow\quad \lambda>0,\quad {{W_{\!I}}}(\lambda)\ne0,\\
&{{\Omega_1}}<0\quad\Longrightarrow\quad
 \bigl({{W_{\!I}}}(\lambda)\ne0\ \Longleftrightarrow\ \lambda\ne{{\Omega_1}}^2d^2\bigr).
\end{aligned}
\]
For every $\lambda>0$, $s(\lambda)=-\sqrt\lambda\ne0$ automatically.
\item For \((\lambda,v)\in{\mathcal O}\), reconstruct \(w=w_0(\lambda)+v\) and \(\psi^\pm=\Phi^\pm\circ\Pi_{Y^\pm}^{-1}\). The residual components satisfy
\begin{equation}\label{eq_identities}
\begin{aligned}
&F_1^\pm=(T^\pm)^2
\bigl[(\Delta\psi^\pm-\Omega_{1,2})\circ\Pi_{Y^\pm}\bigr],\\
&F_2=(T^+)^2
\bigl[|\nabla\psi^+|^2+2gy-Q\bigr]\big|_{y=\eta(x)},\\
&F_3=\bigl(\psi^+_y-\psi^-_y\bigr)\big|_{y=b(x)}.
\end{aligned}
\end{equation}
\item \(F(\lambda,v)=0\) with \((\lambda,v)\in{\mathcal O}\) if and only
if the reconstructed \(\eta,b\in C^{3,\alpha}\) and
\(\psi^\pm\in C^{3,\alpha}(\overline{\Omega^\pm})\) form an even classical graph solution to
\eqref{eq_intro_model} with \(m_2=m_2(\lambda)\), \(\lambda>0\),
the mean conditions \eqref{eq_intro_mean_levels}, \(-h<b<\eta\), and no stagnation on \(\mathcal S\cup\mathcal I\).
The zeros with \(\eta\equiv0\), \(b\equiv-d\) are precisely the laminar line \(v=0\), with complete states \(w_0(\lambda)\).
\end{enumerate}
\end{proposition}
\begin{proof}
For (i), each term in \eqref{eq_lapflat} is a product of strip functions and derivatives of order at most two and hence depends analytically on the unknowns with values in $C^{1,\alpha}$. The first derivative traces belong to $C^{2,\alpha}$, and inversion of a positive thickness is analytic by the geometric series in the H\"older algebra. The laminar profiles are analytic in $\sqrt\lambda$ for $\lambda>0$ and in the physical parameters for $0<d<h$. Substitution gives $F(\lambda,0)=0$ and the stated admissibility conditions. We turn to (ii) and (iii).
Equation \eqref{eq_lapflat} gives the bulk identities in \eqref{eq_identities}. The constant surface trace gives $\psi_x^+=-\eta'\psi_y^+$ there, and \eqref{eq_strip_vel} turns its squared gradient into the expression in $F_2$. Evaluating $\psi_y^\pm=\Phi_t^\pm/T^\pm$ at $t=0$ in the upper layer and $t=1$ in the lower layer, we obtain $F_3$.
The identities \eqref{eq_identities} give the two directions of the correspondence. For each $x$, the positive thickness $T^\pm=\partial_tY^\pm$ makes $Y^\pm(x,\cdot)$ a bijection from $[0,1]$ onto its fluid column, so $\Pi_{Y^\pm}$ is a $C^{3,\alpha}$ diffeomorphism up to the boundary.

From a zero $F(\lambda,v)=0$ in ${\mathcal O}$, reconstruct $\psi^\pm=\Phi^\pm\circ\Pi_{Y^\pm}^{-1}\in C^{3,\alpha}(\overline{\Omega^\pm})$.
The boundary traces and mean conditions are built into $X$.
Identities \eqref{eq_identities} give the Poisson equations, Bernoulli
condition and $[\psi_y]=0$; \eqref{eq_trace_reduction} gives
$[\nabla\psi]=0$. 
Conversely, flatten a classical solution satisfying the conditions in (iii) by setting
\begin{equation}\label{eq_inverse_flat_construction}
\begin{aligned}
&\zeta=b+d,\qquad \Phi^\pm=\psi^\pm\circ\Pi_{Y^\pm},\\
&a=m_1-m_1(\lambda),\qquad q=Q-\lambda,\\
&\delta\varphi^+=\Phi^+-\Phi_0^+-a(1-t),\\
&\delta\varphi^-=\Phi^--\Phi_0^--a,\\
&v=(\eta,\zeta,\delta\varphi^+,\delta\varphi^-,a,q).
\end{aligned}
\end{equation}
By the kinematic conditions we have $\delta\varphi^\pm|_{t=0,1}=0$. Evenness and the two mean conditions then give $v\in X$. Taking $\delta>0$ smaller than the boundary and thickness minima and smaller than $\lambda$, we obtain $(\lambda,v)\in{\mathcal O}_\delta$. It follows from \eqref{eq_identities} that $F(\lambda,v)=0$.

Let $F(\lambda,v)=0$ with $\eta=0$ and $b=-d$. Compare the stream functions $\Phi^\pm$ with $\Phi_0^\pm$ at the same $\lambda$. Since both states have $T^+=d$ and $T^-=h-d$, the differences $\Phi^\pm-\Phi_0^\pm$ solve the homogeneous flat bulk equations. With $a=m_1-m_1(\lambda)$, the upper difference has traces $a,0$ and the lower difference has traces $a,a$ at $t=0,1$, because $m_2(\lambda)$ is fixed. Uniqueness for the periodic Dirichlet problem gives
$\Phi^+-\Phi_0^+=a(1-t),\qquad \Phi^--\Phi_0^-=a.$
Subtracting the interface residuals, we obtain $-a/d=0$. Thus the profiles agree, and the surface residual gives $Q=s^2=\lambda$. Combining these identities, we obtain $v=0$.
\end{proof}

\begin{proposition}\label{rem_components}
$\operatorname{sgn}\psi_y|_{{\mathcal S}}$ and $\operatorname{sgn}\psi_y|_{{\mathcal I}}$ are constant on each connected component of ${\mathcal O}$, and in particular along every connected set of solutions. The laminar flows with $s<0<{{W_{\!I}}}$ lie in the component where $\psi_y<0$ on ${{\mathcal S}}$ and $\psi_y>0$ on ${{\mathcal I}}$. The half-period translation $\mathcal T_\pi$ shifts the profiles by $x\mapsto x+\pi$ and fixes the scalars. It is a linear isometry on $X$ and a bounded linear isomorphism on $Y$. It preserves ${\mathcal O}_\delta$. It satisfies
\[
F(\lambda,\mathcal T_\pi v)=\mathcal T_\pi F(\lambda,v),\qquad
\mathcal T_\pi w_0(\lambda)=w_0(\lambda).
\]
On the Fourier mode $\cos nx$, $\mathcal T_\pi$ acts as multiplication by $(-1)^n$.
\end{proposition}

\subsection{Regularity and compactness}\label{sec_reg}

\label{sec_compact}
We first establish higher regularity under fixed admissibility margins. The resulting estimate gives compactness for bounded families of graph solutions.

\subsubsection{Regularity}

For a coordinate of period $\ell$, write
$|a|_{\rm per}=\min_{m\in\mathbb Z}|a-m\ell|.$
For an arclength coordinate, $\ell$ is the length of the curve over one spatial period.

\begin{lemma}\label{lem_reg}
Let $(\lambda,v)\in{\mathcal O}_\delta$ satisfy $F(\lambda,v)=0$ and set $w=w_0(\lambda)+v$. Then
$\eta,b\in C^\infty(\mathbb S)$ and
$\psi^\pm\in C^\infty(\overline{{\Omega}^\pm})$. We write $X^{(4)}$ for
the space obtained from $X$ by replacing $C^{3,\alpha}$ with $C^{4,\alpha}$.
For every $R>0$ there is $C=C(\delta,R,h,d,{{\Omega_1}},{{\Omega_2}},g)$
such that $|\lambda|+\|w\|_X\le R$ implies
\begin{equation*}
\begin{aligned}
\|\eta\|_{C^{4,\alpha}}+\|b\|_{C^{4,\alpha}}+\|w\|_{X^{(4)}}
+\|\psi^+\|_{C^{4,\alpha}(\overline{{\Omega}^+})}
+\|\psi^-\|_{C^{4,\alpha}(\overline{{\Omega}^-})}\le C.
\end{aligned}
\end{equation*}
\end{lemma}
\begin{proof}
We apply \autoref{lem_reglocal} to the free surface and the interface. To this end we verify the separation and parametrization bounds in that lemma.

Constants depend on $R$ and $\delta$. By \autoref{prop_operator}(iii), $\psi^\pm$ solve \eqref{eq_intro_model}, and since $\psi^\pm=\Phi^\pm\circ\Pi_{Y^\pm}^{-1}$ with $T^\pm\ge\delta$, the margins and \eqref{eq_trace_reduction} give
\[
\|\psi^+\|_{C^{3,\alpha}(\overline{{\Omega}^+})}+\|\psi^-\|_{C^{3,\alpha}(\overline{{\Omega}^-})}\le C,\qquad
|\nabla\psi|\ge\delta\quad\text{on }{{\mathcal S}}\cup{{\mathcal I}} .
\]
For $u\in\{\eta,b\}$ let $P(\sigma)=(x(\sigma),u(x(\sigma)))$ be the arclength parametrization with $x$ increasing. The graph bounds imply $\|P'\|_{C^{2,\alpha}}\le C$, and the length of one period is at most $2\pi(1+R^2)^{1/2}$. In particular,
$x'(\sigma)=\bigl(1+u'(x(\sigma))^2\bigr)^{-1/2}\ge(1+R^2)^{-1/2}.$
Integrating this bound and projecting the chord onto the horizontal coordinate, we obtain
\[
|P(\sigma_1)-P(\sigma_2)|\ge(1+R^2)^{-1/2}\min\{|\sigma_1-\sigma_2|_{\rm per},1\}.
\]
The bed is horizontal, so $\operatorname{dist}({{\mathcal S}},{{\mathcal B}})\ge\operatorname{dist}({{\mathcal I}},{{\mathcal B}})=\inf T^-\ge\delta$. For the surface and the interface, if $|x_1-x_2|<\delta/(2R)$ then
$\eta(x_2)-b(x_1)=T^+(x_2)+b(x_2)-b(x_1)\ge\delta-R|x_1-x_2|>\delta/2,$
so $\operatorname{dist}({{\mathcal S}},{{\mathcal I}})\ge\min\{\delta/2,\delta/(2R)\}$. Hence \autoref{lem_reglocal} applies with $j=3$ and constants depending only on $R,\delta$, and gives
\[
\|P_{{\mathcal S}}'\|_{C^{3,\alpha}}+\|P_{{\mathcal I}}'\|_{C^{3,\alpha}}+\|\psi^+\|_{C^{4,\alpha}(\overline{{\Omega}^+})}+\|\psi^-\|_{C^{4,\alpha}(\overline{{\Omega}^-})}\le C .
\]

The arclength estimate for $P$, together with the lower bound on $x'$, yields a uniform $C^{4,\alpha}$ estimate for the inverse $x\mapsto\sigma(x)$. Hence $u=P_2\circ\sigma$ satisfies the same estimate on each graph. The bounds on $\eta,b$ control $Y^\pm,T^\pm$, and those on $\psi^\pm$ control $\Phi^\pm=\psi^\pm\circ\Pi_{Y^\pm}$ in $C^{4,\alpha}$. We have thus obtained the required $X^{(4)}$ bound for $w$.

Iteration of \autoref{lem_reglocal} with $j=4,5,\dots$, followed at each step by the same change to graph coordinates, gives $\eta,b\in C^\infty(\mathbb S)$ and $\psi^\pm\in C^\infty(\overline{{\Omega}^\pm})$.
\end{proof}

We next apply the bounds of \autoref{lem_reg} to obtain compactness on each member of the following exhaustion.

\begin{proposition}\label{prop_compact}
For $\delta>0$ and each positive integer $j$, let $Q_j(\delta)$ be the set of $(\lambda,v)$ for which the reconstructed state $w=w_0(\lambda)+v$ satisfies
\[
\begin{aligned}
Q_j(\delta)=\Bigl\{&(\lambda,v)\in(0,\infty)\times X:\quad
|\lambda|+\|w_0(\lambda)+v\|_X\le j,\\
&\lambda\ge\delta+\frac1j,\quad
\inf T^+\ge\delta+\frac1j,\quad
\inf T^-\ge\delta+\frac1j,\\
&\inf_{\mathcal S}|\nabla\psi^+|\ge\delta+\frac1j,\quad
\inf_{\mathcal I}|\nabla\psi^+|\ge\delta+\frac1j
\Bigr\}.
\end{aligned}
\]
The sets $Q_j(\delta)$ are closed and bounded in ${\mathbb{R}}\times X$, and satisfy
${\mathcal O}_\delta=\bigcup_{j\ge1}Q_j(\delta).$
For every $j\ge1$, $\{(\lambda,v)\in Q_j(\delta):F(\lambda,v)=0\}
\quad\text{is compact in }{\mathbb{R}}\times X.$
\end{proposition}
\begin{proof}
The margins are continuous wherever the thicknesses and $\lambda$ are positive. Their lower bounds persist under limits, so $Q_j(\delta)$ is closed. Boundedness follows from $\|v\|_X\le j+\max_{[\delta+1/j,j]}\|w_0\|_X$ whenever the set is nonempty. At each admissible point take $j$ large enough that all defining inequalities are strict; the interiors therefore exhaust $\mathcal O_\delta$.

For a sequence of zeros $(\lambda_n,v_n)\in Q_j(\delta)$, set $w_n=w_0(\lambda_n)+v_n$. By \autoref{lem_reg} with $R=j$, the sequence $(\lambda_n,w_n)$ is bounded in ${\mathbb{R}}\times X^{(4)}$. The compact embedding of $C^{4,\alpha}$ into $C^{3,\alpha}$ and the continuity of $w_0$ at the limiting $\lambda\ge\delta+1/j$ give a subsequence for which
\begin{equation*}
\begin{aligned}
(\lambda_n,w_n)&\longrightarrow(\lambda,w) \quad\text{in }{\mathbb{R}}\times X,\\
v_n=w_n-w_0(\lambda_n)&\longrightarrow v=w-w_0(\lambda) \quad\text{in }X.
\end{aligned}
\end{equation*}
Since the set is closed, we have $(\lambda,v)\in Q_j(\delta)$, and by continuity of $F$ we obtain $F(\lambda,v)=0$.
\end{proof}

\subsection{The linearization in physical variables}\label{sec_linearization}

\subsubsection{Physical variables: the conjugation}\label{sec_conjugation}

We fix $\delta>0$ and a zero $(\lambda,v)\in{{\mathcal O}}_\delta$ of $F$. Recover its complete state $(\lambda,w)=\mathcal A_{\rm g}(\lambda,v)$, and assume
\begin{equation}\label{eq_fixedsol}
\eta,b\in C^\infty(\mathbb S),\qquad
\psi^\pm\in C^\infty(\overline{{\Omega}^\pm}).
\end{equation}
To study $\partial_vF(\lambda,v)$ we pass to the physical domains of the fixed solution. We first compute the jump at the interface. We then subtract the coordinate displacement from the variations to obtain the unknowns of the linear transmission problem.

On ${{\mathcal I}}$, the unit normal pointing into ${\Omega}^+$ is $
\nu=\frac{(-b',1)}{\sqrt{1+b'^2}}$.
The vorticity jump determines the vertical second-derivative jump used in the linearized interface condition.
\begin{proposition}\label{rem_kappa}
For the fixed smooth graph solution, the jump in the second vertical derivative of $\psi$ across ${{\mathcal I}}$ satisfies $[\psi_{yy}]=({{\Omega_1}}-{{\Omega_2}})\nu_y^2$.
\end{proposition}

At fixed $\lambda$ take $v_\varepsilon=v+\varepsilon\xi$, with $\xi\in X$, and write hats for derivatives at $\varepsilon=0$. The complete-state derivative is ${\widehat w}=\xi=({\widehat\eta},{\widehat\zeta},{\widehat\varphi}^+,{\widehat\varphi}^-,{\widehat m}_1,{\widehat Q})$. Differentiating \eqref{eq_flat} and \eqref{eq_reconstruct_phi} at fixed $(x,t)$, we obtain
\[
\begin{aligned}
{\widehat T}^+&={\widehat\eta}-{\widehat\zeta},& {\widehat T}^-&={\widehat\zeta},\\
{\widehat Y}^+&={\widehat\zeta}+t{\widehat T}^+,& {\widehat Y}^-&=t{\widehat\zeta},\\
{\widehat\Phi}^+&={\widehat m}_1(1-t)+{\widehat\varphi}^+,&
{\widehat\Phi}^-&={\widehat m}_1+{\widehat\varphi}^-.
\end{aligned}
\]

We fix a point $(x,y)$ in the interior of one layer. For sufficiently small $|\varepsilon|$ it remains in that layer, so we can define ${\widehat\Psi}^\pm(x,y)$ as the derivative of $\psi_\varepsilon^\pm$ at this point. Differentiating $\Phi^\pm=\psi^\pm\circ\Pi_{Y^\pm}$ with the strip coordinates held fixed, we obtain the physical perturbation
\begin{equation}\label{eq_psihat}
{\widehat\Psi}^\pm\big(x,Y^\pm(x,t)\big)
={\widehat\Phi}^\pm(x,t)-\psi^\pm_y\big(x,Y^\pm(x,t)\big)\,{\widehat Y}^\pm(x,t).
\end{equation}
By \eqref{eq_spaces}, the strip variation of \eqref{eq_flat}, and the reconstruction \eqref{eq_reconstruct_phi}, ${\widehat\Phi}^\pm$ and ${\widehat Y}^\pm$ are $C^{3,\alpha}$ on the closed strips. The smoothness of $\psi^\pm$ and the strip diffeomorphisms up to the boundaries gives a $C^{3,\alpha}(\overline{{\Omega}^\pm})$ extension of ${\widehat\Psi}^\pm$ by \eqref{eq_psihat}, with boundary traces given by $t=0,1$.

On ${{\mathcal S}}$ and ${{\mathcal I}}$ the admissibility condition \eqref{eq_odelta} gives a positive lower bound for $|\nabla\psi^+|$. Since $\psi$ is constant on each graph, it follows that $\psi_y\ne0$ there.
Eliminating \({\widehat\eta}\) from the linearized Bernoulli condition and \({\widehat\zeta}\) from the linearized transmission condition, we obtain the following coefficients of the fixed background solution. The formula for \(\kappa\) uses \autoref{rem_kappa}.
\begin{equation}\label{eq_betakappa}
\begin{aligned}
\beta :=-\frac{2\nabla\psi\cdot\nabla\psi_y+2g}{\psi_y}
 \text{ on }{{\mathcal S}},\qquad
\kappa :=\frac{[\psi_{yy}]}{\psi_y}
=\frac{{{\Omega_1}}-{{\Omega_2}}}{(1+b'^2)\,\psi_y}
 \text{ on }{{\mathcal I}}.
\end{aligned}
\end{equation}

\begin{lemma}\label{lem_conj}
For a zero $(\lambda,v)\in{{\mathcal O}}_\delta$ of $F$ satisfying \eqref{eq_fixedsol}, use the physical variations in \eqref{eq_psihat} to define
\[
\mathcal C{\widehat w}=({\widehat\Psi}^+,{\widehat\Psi}^-,{\widehat m}_1,{\widehat Q}).
\]
Then $\mathcal C$ is a linear isomorphism of $X$ onto the space
\begin{equation}\label{eq_xtilde}
\begin{aligned}
\widetilde X:=\Bigg\{&({\widehat\Psi}^+,{\widehat\Psi}^-,{\widehat m}_1,{\widehat Q})\in \hat{X}:\ [{\widehat\Psi}]=0\text{ on }{{\mathcal I}},\quad {\widehat\Psi}^-={\widehat m}_1\text{ on }{{\mathcal B}},\\
&\Big\langle\frac{{\widehat\Psi}^+}{\psi_y}\Big\rangle_{{{\mathcal S}}}=0,\quad \Big\langle\frac{{\widehat\Psi}-{\widehat m}_1}{\psi_y}\Big\rangle_{{{\mathcal I}}}=0\Bigg\},\\
& \hat{X}=C^{3,\alpha}(\overline{{\Omega}^+})\times C^{3,\alpha}(\overline{{\Omega}^-})\times{\mathbb{R}}^2,
\end{aligned}
\end{equation}
where all functions in $\widetilde X$ are even in $x$ and ${\widehat\Psi}$ on \({{\mathcal I}}\) denotes the common trace of ${\widehat\Psi}^+$ and ${\widehat\Psi}^-$. The means $\langle\cdot\rangle_{{\mathcal S}}$ and $\langle\cdot\rangle_{{\mathcal I}}$ are taken in $x$ at $(x,\eta(x))$ and $(x,b(x))$, respectively. The inverse gives the shape perturbations
\[
{\widehat\eta}=-{\widehat\Psi}^+/\psi_y \quad\text{on }{{\mathcal S}},\qquad
{\widehat\zeta}=({\widehat m}_1-{\widehat\Psi})/\psi_y \quad\text{on }{{\mathcal I}}.
\]
For the residuals, identify functions on the two graphs with functions of $x$. The physical target space, with $e$ denoting evenness in $x$, is
\[
\widetilde Y:=C_{e}^{1,\alpha}(\overline{{\Omega}^+})\times C_{e}^{1,\alpha}(\overline{{\Omega}^-})\times C_{e}^{2,\alpha}({{\mathcal S}})\times C_{e}^{2,\alpha}({{\mathcal I}}).
\]
The linear isomorphism $\widehat{\mathcal C}:Y\to\widetilde Y$ is given by
\[
\begin{aligned}
\widehat{\mathcal C}(\widehat f_1^+,\widehat f_1^-,\widehat f_2,\widehat f_3)
=\Bigg(&\frac{\widehat f_1^+\circ\Pi_{Y^+}^{-1}}{(T^+)^2},
        \frac{\widehat f_1^-\circ\Pi_{Y^-}^{-1}}{(T^-)^2},
        \frac{\widehat f_2}{(T^+)^2},\ \widehat f_3\Bigg).
\end{aligned}
\]
With these identifications and the coefficients in \eqref{eq_betakappa}, the conjugated linearization is
\begin{equation}\label{eq_ltilde}
\begin{aligned}
\widetilde L&:=\widehat{\mathcal C}\circ\partial_vF(\lambda,v)\circ\mathcal C^{-1},\\
\widetilde L({\widehat\Psi}^+,{\widehat\Psi}^-,{\widehat m}_1,{\widehat Q})&=
\begin{pmatrix}
{\Delta}{\widehat\Psi}^+\\
{\Delta}{\widehat\Psi}^-\\
2\nabla\psi\cdot\nabla{\widehat\Psi}^++\beta{\widehat\Psi}^+-{\widehat Q}\ \text{on }{{\mathcal S}}\\
[{\widehat\Psi}_y]+\kappa({\widehat m}_1-{\widehat\Psi})\ \text{on }{{\mathcal I}}
\end{pmatrix}^T.
\end{aligned}
\end{equation}
\end{lemma}
\begin{proof}
The formula \eqref{eq_psihat}, with its smooth background coefficients and strip diffeomorphisms, gives \({\widehat\Psi}^\pm\in C^{3,\alpha}(\overline{{\Omega}^\pm})\) with $
\|{\widehat\Psi}^\pm\|_{C^{3,\alpha}}\le C\|{\widehat w}\|_X $.
Evaluating \eqref{eq_psihat} at the strip endpoints, we obtain the traces
\begin{equation*}
\begin{aligned}
&{\widehat\Psi}^+(x,\eta(x))=-\psi_y(x,\eta(x)){\widehat\eta},\\
&{\widehat\Psi}^\pm(x,b(x))={\widehat m}_1-\psi_y^\pm(x,b(x)){\widehat\zeta}(x),\\
&{\widehat\Psi}^-(x,-h)={\widehat m}_1.
\end{aligned}
\end{equation*}
The continuity of $\psi_y$ gives $[{\widehat\Psi}]=0$, and the two mean conditions in \eqref{eq_xtilde} express 
$\langle{\widehat\eta}\rangle=\langle{\widehat\zeta}\rangle=0.$

For the inverse, use the shape displacements in the statement. Since $\psi_y$ is smooth and bounded away from zero by \eqref{eq_odelta}, these even profiles belong to $C^{3,\alpha}$ and satisfy the zero means. With ${\widehat Y}^\pm$ given by the strip variation preceding \eqref{eq_psihat}, the remaining components follow from \eqref{eq_psihat} and the differentiated reconstruction:
\begin{equation*}
\begin{aligned}
&{\widehat\Phi}^\pm={\widehat\Psi}^\pm\circ\Pi_{Y^\pm}
+\psi^\pm_y\circ\Pi_{Y^\pm}\,{\widehat Y}^\pm ,\\
&{\widehat\varphi}^+={\widehat\Phi}^+-{\widehat m}_1(1-t),\qquad
{\widehat\varphi}^-={\widehat\Phi}^--{\widehat m}_1 .
\end{aligned}
\end{equation*}
The boundary conditions on ${\widehat\Psi}^\pm$ imply ${\widehat\varphi}^\pm=0$ at $t=0,1$. Hence we obtain ${\widehat w}=({\widehat\eta},{\widehat\zeta},{\widehat\varphi}^+,{\widehat\varphi}^-,{\widehat m}_1,{\widehat Q})\in X$, and the two constructions are mutually inverse. The inverse is bounded by the reconstruction formulas and the positive lower bound for \(|\psi_y|\) on the surface and interface.
Since $\Pi_{Y^\pm}$ are smooth diffeomorphisms and $T^\pm\ge\delta$, the linear map $\widehat{\mathcal C}$ has a bounded inverse.

Differentiate the bulk identities in \eqref{eq_identities} at fixed $\lambda$. The variation of the factor $(T^\pm)^2$ multiplies the zero residual, and the displacement term vanishes because $\partial_y{\Delta}\psi^\pm=0$ in each layer. The resulting bulk components of \eqref{eq_ltilde} follow from
\[
\partial_vF_1^\pm\,{\widehat w}=(T^\pm)^2({\Delta}{\widehat\Psi}^\pm)\circ\Pi_{Y^\pm}.
\]
For the surface component, we use the Bernoulli residual. By \eqref{eq_identities} it is
\[
F_2=(T^+)^2\,\mathcal B(\lambda,v),\qquad
\mathcal B(\lambda,v)(x)=|\nabla\psi^+|^2(x,\eta(x))+2g\eta(x)-Q .
\]
At the solution $\mathcal B=0$, so $\partial_vF_2\,{\widehat w}=(T^+)^2\partial_v\mathcal B\,{\widehat w}$. Differentiating the moving evaluation at $(x,\eta(x))$ and using the surface displacement from Step 1, we obtain
\begin{equation*}
\begin{aligned}
\partial_v\mathcal B\,{\widehat w}
&=2\nabla\psi\cdot\nabla{\widehat\Psi}^+
  +\bigl(\partial_y|\nabla\psi|^2+2g\bigr){\widehat\eta}-{\widehat Q}\\
&=2\nabla\psi\cdot\nabla{\widehat\Psi}^+
  -\frac{2\nabla\psi\cdot\nabla\psi_y+2g}{\psi_y}{\widehat\Psi}^+-{\widehat Q}\\
&=2\nabla\psi\cdot\nabla{\widehat\Psi}^++\beta{\widehat\Psi}^+-{\widehat Q}.
\end{aligned}
\end{equation*}
After division by $(T^+)^2$ in $\widehat{\mathcal C}$, this gives the third component of \eqref{eq_ltilde}.

At \((x,b(x))\), differentiate the interface residual \(F_3=(\psi^+_y-\psi^-_y)(x,b(x))\) and use the displacement from Step 1 and the coefficient \(\kappa\) in \eqref{eq_betakappa}. This gives the fourth component of \eqref{eq_ltilde},
\begin{equation*}
\begin{aligned}
\partial_vF_3\,{\widehat w}
&=[{\widehat\Psi}_y]+[\psi_{yy}]{\widehat\zeta}
=[{\widehat\Psi}_y]+\frac{[\psi_{yy}]}{\psi_y}({\widehat m}_1-{\widehat\Psi})
=[{\widehat\Psi}_y]+\kappa({\widehat m}_1-{\widehat\Psi}).
\end{aligned}
\end{equation*}
\end{proof}

\subsection{The Fredholm property of the linearization}\label{sec_fred}

We now prove that the linearized operator in \eqref{eq_ltilde} has index zero. We begin with the transmission problem for the Laplacian.

\begin{lemma}\label{lem_transmission}
Let ${\Omega}$ be a connected open set with ${\Omega}+(L,0)={\Omega}$, bounded below by the bed ${{\mathcal B}}=\{y=-h\}$ and above by an $L$-periodic $C^\infty$ curve ${{\mathcal S}}\subset\{y>-h\}$ without self-intersections, and let ${{\mathcal I}}\subset{\Omega}$ be an $L$-periodic $C^\infty$ curve without self-intersections that divides ${\Omega}$ into two connected open sets, \[
{\Omega}^+\ \text{with}\ \partial{\Omega}^+={{\mathcal S}}\cup{{\mathcal I}}\ \text{and}\ {\Omega}^-\ \text{with}\ \partial{\Omega}^-={{\mathcal I}}\cup{{\mathcal B}}.
\] Let $n$ denote the unit normal of ${{\mathcal S}}$ pointing out of ${\Omega}$ and the unit normal of ${{\mathcal I}}$ pointing into ${\Omega}^+$, and let $[\cdot]$ denote the jump across ${{\mathcal I}}$. All functions below are $L$-periodic, and norms are taken over one period. Then
\begin{enumerate}[label=\textup{(\alph*)},leftmargin=2.3em]
\item For all $f^\pm\in C^{1,\alpha}(\overline{{\Omega}^\pm})$, $g_{{\mathcal S}}\in C^{2,\alpha}({{\mathcal S}})$ and $g_{{\mathcal I}}\in C^{2,\alpha}({{\mathcal I}})$ there is exactly one function $\Psi$ on $\overline{{\Omega}}$ with $\Psi^\pm:=\Psi|_{\overline{{\Omega}^\pm}}\in C^{3,\alpha}(\overline{{\Omega}^\pm})$ and
\begin{equation}\label{eq_transm}
\begin{aligned}
&{\Delta}\Psi^\pm=f^\pm \quad\text{in }{\Omega}^\pm,\qquad
\partial_n\Psi^+=g_{{\mathcal S}} \quad\text{on }{{\mathcal S}},\\
&[\Psi]=0,\quad [\partial_n\Psi]=g_{{\mathcal I}} \quad\text{on }{{\mathcal I}},\qquad
\Psi^-=0 \quad\text{on }{{\mathcal B}},
\end{aligned}
\end{equation}
and
\[
\|\Psi^+\|_{C^{3,\alpha}}+\|\Psi^-\|_{C^{3,\alpha}}
\le C\bigl(\|f^+\|_{C^{1,\alpha}}+\|f^-\|_{C^{1,\alpha}}
+\|g_{{\mathcal S}}\|_{C^{2,\alpha}}+\|g_{{\mathcal I}}\|_{C^{2,\alpha}}\bigr).
\]
\item If ${\Omega}$ and ${{\mathcal I}}$ are symmetric about $\{x=0\}$ and $f^\pm$, $g_{{\mathcal S}}$, $g_{{\mathcal I}}$ are even in $x$, then $\Psi$ is even in $x$.
\end{enumerate}
\end{lemma}
\begin{proof}
On the quotient by the translation \((L,0)\), the zero Dirichlet condition on \({{\mathcal B}}\) makes the Dirichlet form coercive. The Lax--Milgram theorem therefore gives a unique weak solution. For the interface estimate, freeze the Laplacian in orthonormal tangential and normal coordinates. At frequency $\xi\ne0$ the decaying amplitudes satisfy $A_+=A_-$ and $-|\xi|(A_++A_-)=0$ for the homogeneous transmission data, so both vanish. The surface and bed carry the usual Neumann and Dirichlet conditions, respectively. The boundary and transmission estimates of \cite{GT,ADN2} give the regularity and bound in (a). See also \cite{Schechter} and \cite[Appendix~9, Theorem~9.1]{LV}.
For (b), reflect the solution across $x=0$. The reflected function satisfies \eqref{eq_transm} with the same data, so uniqueness gives $\Psi(-x,y)=\Psi(x,y)$.
\end{proof}

\begin{lemma}\label{lem_model}
For a zero $(\lambda,v)\in{{\mathcal O}}_\delta$ of $F$ satisfying \eqref{eq_fixedsol}, let $\widetilde X$ and $\widetilde Y$ be the spaces of \autoref{lem_conj}. Define $\widetilde L_0:\widetilde X\to\widetilde Y$ by
\[
\widetilde L_0({\widehat\Psi}^+,{\widehat\Psi}^-,{\widehat m}_1,{\widehat Q})
=\bigl({\Delta}{\widehat\Psi}^+,\ {\Delta}{\widehat\Psi}^-,\ 
(2\nabla\psi\cdot\nabla{\widehat\Psi}^+-{\widehat Q})|_{{\mathcal S}},\ 
[{\widehat\Psi}_y]|_{{\mathcal I}}\bigr).
\]
Then $\widetilde L_0$ is Fredholm of index zero.
\end{lemma}
\begin{proof}
We compute the index in two stages. We first adjoin the two scalar unknowns to the normal transmission problem of \autoref{lem_transmission}. Restriction by the two mean conditions then gives $\widetilde L_0$.
On ${{\mathcal S}}$, let $\nu_{{\mathcal S}}$ be the outward unit normal and write $\sigma_{{\mathcal S}}=\operatorname{sgn}\psi_y|_{{\mathcal S}}$. Since $\psi=0$ and $\nabla\psi\ne0$ there,
\[
\nabla\psi=\sigma_{{\mathcal S}}|\nabla\psi|\nu_{{\mathcal S}}
\qquad\text{on }{{\mathcal S}}.
\]
The trace condition $[{\widehat\Psi}]=0$ on ${{\mathcal I}}$ gives $[\nabla{\widehat\Psi}]=[\partial_\nu{\widehat\Psi}]\nu$ there, where $\nu=(-b',1)/\sqrt{1+b'^2}$ is the upward unit normal.
The boundary operators therefore satisfy
\[
2\nabla\psi\cdot\nabla{\widehat\Psi}^+=2\sigma_{{\mathcal S}}|\nabla\psi|\,\partial_{\nu_{{\mathcal S}}}{\widehat\Psi}^+\quad\text{on }{{\mathcal S}},\qquad
[{\widehat\Psi}_y]=\frac{[\partial_\nu{\widehat\Psi}]}{\sqrt{1+b'^2}}\quad\text{on }{{\mathcal I}}.
\]
Let $U$ be the space of even pairs $u=(u^+,u^-)\in C^{3,\alpha}(\overline{{\Omega}^+})\times C^{3,\alpha}(\overline{{\Omega}^-})$ with $[u]=0$ on ${{\mathcal I}}$ and $u^-=0$ on ${{\mathcal B}}$, and put
\[
Au=\bigl({\Delta} u^+,\ {\Delta} u^-,\ 2\nabla\psi\cdot\nabla u^+|_{{\mathcal S}},\ [u_y]|_{{\mathcal I}}\bigr)\in\widetilde Y .
\]
Since $|\nabla\psi|\ge\delta$ on ${{\mathcal S}}$, these identities turn $Au=(f^+,f^-,g_{{\mathcal S}},g_{{\mathcal I}})$ into \eqref{eq_transm} with data
\[
\Bigl(f^+,\ f^-,\ \frac{g_{{\mathcal S}}}{2\sigma_{{\mathcal S}}|\nabla\psi|},\ \sqrt{1+b'^2}\,g_{{\mathcal I}}\Bigr).
\]
Thus $A:U\to\widetilde Y$ is an isomorphism by \autoref{lem_transmission}(a),(b).

Let $\widehat X_0$ be the space defined by \eqref{eq_xtilde} before imposing the two mean conditions.
The map $(u,{\widehat m}_1,{\widehat Q})\mapsto(u+{\widehat m}_1,{\widehat m}_1,{\widehat Q})$ is an isomorphism from $U\times{\mathbb{R}}^2$ onto $\widehat X_0$. In these coordinates, the equation for $f\in\widetilde Y$ is equivalent to
\[
\widetilde L_0(u+{\widehat m}_1,{\widehat m}_1,{\widehat Q})=f
\quad\Longleftrightarrow\quad
u=A^{-1}\bigl(f+(0,0,{\widehat Q},0)\bigr).
\]
For each $f$, both scalars are free. Hence the operator is onto and its kernel has dimension two.

Let $M$ be the pair of bounded mean functionals in \eqref{eq_xtilde}. We choose an even smooth $\chi$ equal to $1$ near ${\mathcal I}$ and to $0$ near ${\mathcal S}\cup{\mathcal B}$. Evaluating $M$ on two states $({\widehat\Psi}^+,{\widehat\Psi}^-,{\widehat m}_1,{\widehat Q})\in\widehat X_0$, we obtain
\[
\begin{aligned}
&M(1,1,1,0)=\bigl(\langle 1/\psi_y\rangle_{{\mathcal S}},0\bigr),\\
&M(\chi|_{\Omega^+},\chi|_{\Omega^-},0,0)
=\bigl(0,\langle 1/\psi_y\rangle_{{\mathcal I}}\bigr).
\end{aligned}
\]
Both means are nonzero since $\psi_y$ has a constant nonzero sign on each graph by \eqref{eq_trace_reduction} and the admissibility margins. Therefore $\widetilde X=\ker M$ has codimension two, and the restriction of $\widetilde L_0$ to $\widetilde X$ has Fredholm index zero.
\end{proof}

\begin{lemma}\label{thm_fred}
For a zero $(\lambda,v)\in{{\mathcal O}}_\delta$ of $F$ satisfying \eqref{eq_fixedsol}, the linearization $\partial_vF(\lambda,v):X\to Y$ is Fredholm of index zero.
\end{lemma}
\begin{proof}
By \autoref{lem_conj}, we have
$\partial_vF(\lambda,v)=\widehat{\mathcal C}^{-1}\circ\widetilde L\circ\mathcal C$
with isomorphisms $\mathcal C:X\to\widetilde X$ and $\widehat{\mathcal C}:Y\to\widetilde Y$.
These isomorphisms preserve the closed range property and the dimensions of the kernel and cokernel. It therefore suffices to show that $\widetilde L=\widetilde L_0+K$ is Fredholm of index zero, where
\[
K({\widehat\Psi}^+,{\widehat\Psi}^-,{\widehat m}_1,{\widehat Q})
=(0,0,\beta{\widehat\Psi}^+|_{{\mathcal S}},\kappa({\widehat m}_1-{\widehat\Psi})|_{{\mathcal I}}).
\]
Since $\beta$ and $\kappa$ are smooth, the two boundary terms are bounded into $C^{3,\alpha}({{\mathcal S}})\times C^{3,\alpha}({{\mathcal I}})$. Arzelà--Ascoli gives the compact embedding
\[
C^{3,\alpha}({{\mathcal S}})\times C^{3,\alpha}({{\mathcal I}})
\hookrightarrow C^{2,\alpha}({{\mathcal S}})\times C^{2,\alpha}({{\mathcal I}}).
\]
Thus $K:\widetilde X\to\widetilde Y$ is compact, so $\widetilde L$ is Fredholm of index zero by \autoref{lem_model} and the compact perturbation theorem \cite[Theorem~IV.5.26]{Kato}.
\end{proof}

\begin{corollary}\label{cor_fredall}
Let $(\lambda,v)\in{\mathcal O}$ satisfy $F(\lambda,v)=0$. Then $\partial_vF(\lambda,v):X\to Y$ is Fredholm of index zero.
\end{corollary}

For the interface Schauder estimate, reflect the lower equation across the flattened interface. At tangential frequency $\xi\ne0$ the two decaying amplitudes agree by continuity; the conormal jump of the reflected pair is $-2|\xi|\sqrt{\det a}$ times this amplitude when the frozen principal matrix $a$ has a common trace. Thus the complementing condition holds. Lower-order terms and tangential derivatives of the common trace leave the principal condition unchanged \cite{ADN2}.
\begin{lemma}\label{lem_reglocal}
Let ${\Omega},{{\mathcal S}},{{\mathcal I}},{\Omega}^\pm$ be as in \textup{\autoref{lem_transmission}}, except that ${{\mathcal S}}$ and ${{\mathcal I}}$ need only be of class $C^{3,\alpha}$. Let $\psi^\pm\in C^{3,\alpha}(\overline{{\Omega}^\pm})$ solve \eqref{eq_intro_model}, and fix $j\ge3$, $R\ge1$ and $\delta\in(0,1)$.

For each $\Gamma\in\{{{\mathcal S}},{{\mathcal I}}\}$, suppose there is an arclength parametrization $P_\Gamma$ of length at most $R$ per period. For all arclength parameters $\sigma_1,\sigma_2$, assume
\[
\begin{aligned}
\|P_\Gamma'\|_{C^{j-1,\alpha}}&\le R,\\
|P_\Gamma(\sigma_1)-P_\Gamma(\sigma_2)|
&\ge\delta\min\{|\sigma_1-\sigma_2|_{\rm per},1\}.
\end{aligned}
\]
Assume the non-stagnation and separation bounds
\[
\begin{aligned}
|\nabla\psi|&\ge\delta\quad\text{on }{{\mathcal S}}\cup{{\mathcal I}},\qquad
\operatorname{dist}({{\mathcal S}},{{\mathcal I}}),\quad
\operatorname{dist}({{\mathcal I}},{{\mathcal B}}),\quad
\operatorname{dist}({{\mathcal S}},{{\mathcal B}})\ge\delta.
\end{aligned}
\]
If $\|\psi^\pm\|_{C^{j,\alpha}(\overline{{\Omega}^\pm})}\le R$, then
\[
P_{{\mathcal S}}',\ P_{{\mathcal I}}'\in C^{j,\alpha},
\qquad
\psi^\pm\in C^{j+1,\alpha}(\overline{{\Omega}^\pm}).
\]
Under the same hypotheses with $j=3$, the gain in regularity satisfies the uniform estimate
\[
\|P_{{\mathcal S}}'\|_{C^{3,\alpha}},\quad
\|P_{{\mathcal I}}'\|_{C^{3,\alpha}},\quad
\|\psi^\pm\|_{C^{4,\alpha}(\overline{{\Omega}^\pm})}\le C,
\]
where $C$ depends only on $R,\delta$ and $g,h,{{\Omega_1}},{{\Omega_2}},\alpha$.
\end{lemma}
\begin{proof}
Boundary nonstagnation gives local stream-function coordinates at the free surface and interface, and the geometric bounds give uniformly sized charts. Schauder estimates for the tangential derivative and the height equation for the remaining normal derivatives give local gains, which we combine with bed reflection and interior estimates. We give the estimates for $j=3$. Raising the indices by $j-3$ gives the regularity gain for $j>3$.

\smallskip\noindent\textit{Step 1. Uniformly separated charts.}\enspace
Let $P=P_{{\mathcal S}}(\sigma_0)$. Since $|P_{{\mathcal S}}''|\le R$, the unit tangent $P_{{\mathcal S}}'(\sigma)$ makes an angle at most $R|\sigma-\sigma_0|$ with $P_{{\mathcal S}}'(\sigma_0)$. We set $\sigma_1=\min\{1/(4R),1\}$. On $|\sigma-\sigma_0|<\sigma_1$ the angle is at most $1/4$, so the projection of $P_{{\mathcal S}}(\sigma)-P$ onto the tangent line at $P$ has derivative at least $\cos\frac14>\frac12$. The arc $\{P_{{\mathcal S}}(\sigma):|\sigma-\sigma_0|<\sigma_1\}$ is therefore the graph of a $C^{3,\alpha}$ function $\mathsf g$ over an interval containing $(-\sigma_1/2,\sigma_1/2)$ of the tangent line. The function has $\mathsf g(0)=\mathsf g'(0)=0$ and $\|\mathsf g\|_{C^{3,\alpha}}\le C$. The slope satisfies
\begin{equation*}
\begin{aligned}
|\mathsf g'|&\le\tan\frac14
<1.
\end{aligned}
\end{equation*}
By the separation hypothesis, points $P_{{\mathcal S}}(\sigma)$ with $|\sigma-\sigma_0|_{\rm per}\ge\sigma_1$ satisfy $|P_{{\mathcal S}}(\sigma)-P|\ge\delta\sigma_1$. Hence, for $r_1:=\delta\sigma_1/2$, the set ${{\mathcal S}}\cap B_{r_1}(P)$ is contained in this graph, and $B_{r_1}(P)$ does not meet ${{\mathcal I}}\cup{{\mathcal B}}$. Thus $\overline{\Omega}\cap B_{r_1}(P)$ is the part of $B_{r_1}(P)$ on one side of the graph. The corresponding graph description holds at every $P\in{{\mathcal I}}$, with one side of the graph in ${\Omega}^+$ and the other in ${\Omega}^-$.

\smallskip\noindent\textit{Step 2. The free-surface regularity gain.}\enspace
Let $P\in{{\mathcal S}}$, and choose coordinates $z'=(x',y')=\mathrm{Rot}(z-P)$. The inverse rotation $\mathrm{Rot}^{-1}$ has angle $\theta$ such that $e_{y'}=\nabla\psi(P)/|\nabla\psi(P)|$, which is $\pm n(P)$; so the tangent line at $P$ is the $x'$-axis and 
${{\mathcal S}}\cap B_{r_1}(P)=\{y'=\mathsf g(x')\}$.
 Put 
$r_0:=\min\{r_1/2,\delta/(4R)\}$.

Every $z\in\overline{{\Omega}^+}\cap B_{r_0}(P)$, $z'=(x',y')$, is joined to $P$ inside $\overline{{\Omega}^+}\cap B_{r_1}(P)$ by the vertical segment from $z'$ to $(x',\mathsf g(x'))$ followed by the arc of the graph from $(x',\mathsf g(x'))$ to $0$. The slope bound in Step~1 gives $|\mathsf g(x')|\le\tan\frac14\,|x'|$ and bounds the path length by
$|y'|+1.3|x'|\le2|z-P|.$
The Hessian bound transfers the derivative at $P$ along this path: since $\partial_{y'}\psi(P)=|\nabla\psi(P)|\ge\delta$ and $|D^2\psi^+|\le R$, it follows that
\begin{equation*}
\begin{aligned}
\partial_{y'}\psi&\ge\delta-2Rr_0
\ge\delta/2\qquad\text{on }\overline{{\Omega}^+}\cap B_{r_0}(P).
\end{aligned}
\end{equation*}
Let $\epsilon\in\{\pm1\}$ be such that ${\Omega}$ lies on the side $\epsilon(y'-\mathsf g(x'))<0$ near $P$, and let
$G=\{|x'|<r_0/4,\ 0\le\epsilon(\mathsf g(x')-y')<r_0/2\}$;
 then $G\subset\overline{{\Omega}^+}\cap B_{r_0}(P)$ because $|\mathsf g(x')|\le|x'|$. On each vertical segment of $G$, $\psi$ is strictly increasing in $y'$, vanishes at the graph and changes by at least $\delta r_0/4$. So $(x',y')\mapsto(q,p)=(x',\psi)$ is injective on $G$, its image contains the box 
$\mathbb B=(-r_0/4,r_0/4)\times\{0\le-\epsilon p<\delta r_0/4\}.$
The inverse $y'=H(q,p)$ is $C^{3,\alpha}$ on $\mathbb B$ with bounds, by the formulas 
$H_p=1/\psi_{y'}$ and $H_q=-\psi_{x'}/\psi_{y'}$. In particular, $H_p\in[1/R,2/\delta]$.

Since ${\Delta}\psi={{\Omega_1}}$ is invariant under rotations, $H$ satisfies in $\mathbb B$
\begin{equation}\label{eq_height}
\mathcal Q(H):=(1+H_q^2)H_{pp}-2H_qH_pH_{qp}+H_p^2H_{qq}+{{\Omega_1}}H_p^3=0 .
\end{equation}
Its principal coefficient matrix has determinant $H_p^2\ge R^{-2}$ and trace bounded in terms of $R,\delta$, so \eqref{eq_height} is uniformly elliptic. With $v=H_q$, tangential differentiation gives in $\mathbb B$
\begin{equation}\label{eq_vlin}
\begin{aligned}
0={}&(1+H_q^2)v_{pp}-2H_qH_pv_{qp}+H_p^2v_{qq}
+\big(2H_qH_{pp}-2H_pH_{qp}\big)v_q\\
&+\big(2H_pH_{qq}-2H_qH_{qp}+3{{\Omega_1}}H_p^2\big)v_p,
\end{aligned}
\end{equation}
while on $p=0$,
\[
\frac{1+H_q^2}{H_p^2}+2g(q\sin\theta+H\cos\theta)=Q-2gy(P).
\]
Tangential differentiation gives the boundary equation for $v$:
\[
-\frac{2(1+H_q^2)}{H_p^3}v_p+\frac{2H_q}{H_p^2}v_q+2g\cos\theta\,v=-2g\sin\theta.
\]
The equation coefficients are uniformly $C^{1,\alpha}$, the boundary coefficients are uniformly $C^{2,\alpha}$, and $H_p\in[1/R,2/\delta]$ gives uniform obliqueness.

The interior equation \eqref{eq_vlin} for $v=H_q$ and its oblique boundary condition give $v\in C^{3,\alpha}$ on the half-size box by the oblique Schauder estimate \cite[Theorem~7.3]{ADN1} with $m=1$. The equation \eqref{eq_height} then recovers
\[
H_{pp}=\frac{2H_qH_pv_p-H_p^2v_q-{{\Omega_1}}H_p^3}{1+H_q^2}\in C^{2,\alpha}.
\]
Together the bounds for $H_q$ and $H_{pp}$ give $H\in C^{4,\alpha}$.

Thus $q\mapsto P+\mathrm{Rot}^{-1}(q,H(q,0))$ has speed at least $1$ and is uniformly $C^{4,\alpha}$, as is the inverse $\psi^+$ because $H_p\ge1/R$. For $z\in\overline{\Omega}\cap B_\rho(P)$ with $\rho\le r_0/8$, the rotated point lies in $G$. The path length bound and $\psi(P)=0$ give
$|q|<\rho,\qquad |p|=|\psi(z)|<3R\rho.$
Thus the inverse coordinate map sends the half-size box onto a set containing $\overline{\Omega}\cap B_{\rho_1}(P)$ for $\rho_1=\min\{r_0/8,\delta r_0/(24R)\}$.

\smallskip\noindent\textit{Step 3. The interface regularity gain.}\enspace
In the rotated stream-function coordinates above at $P\in{{\mathcal I}}$, continuity of $\psi,\nabla\psi$ and the one-sided Hessian bounds give $\partial_{y'}\psi\ge\delta/2$ on $B_{r_0}(P)$. Hence $(q,p)=(x',\psi-m_1)$ has a piecewise $C^{3,\alpha}$ inverse $H$ on
$\mathbb B_{{\mathcal I}}=(-r_0/4,r_0/4)\times(-\delta r_0/4,\delta r_0/4).$
For \(v=H_q\), the equations \eqref{eq_height}, \eqref{eq_vlin} hold with \({{\Omega_1}}\) on the \({\Omega}^+\) side and \({{\Omega_2}}\) on the \({\Omega}^-\) side. Continuity of \(\psi,\nabla\psi\) gives matching traces of \(H,H_q,H_p\), so differentiation along \(p=0\) gives matching traces of \(v_q=H_{qq}\) and \(v_p=\partial_qH_p\). These traces and the resulting flux condition on \(p=0\) are
\begin{equation*}
\begin{aligned}
[H]=[H_q]=[H_p]&=0,\qquad 
[v]=[v_q]=[v_p]=0,\\
\big[(1+H_q^2)v_p-H_qH_pv_q\big]&=0.
\end{aligned}
\end{equation*}
The principal coefficients in \eqref{eq_vlin} depend only on $\nabla H$ and have a common trace. The transmission estimate \cite[Theorem~9.3]{ADN2} with $m=1$ and recovery of $H_{pp}$ from the height equation give one-sided $C^{4,\alpha}$ bounds for $H,\psi^\pm$ and a $C^{4,\alpha}$ bound for the interface parametrization. The image of the inverse coordinate chart on each side contains $\overline{{\Omega}^\pm}\cap B_{\rho_1}(P)$ after a uniform decrease of $\rho_1$.

\smallskip\noindent\textit{Step 4. Global bounds by covering and reflection.}\enspace
Since the lengths of ${{\mathcal S}}$ and ${{\mathcal I}}$ per period are at most $R$, finitely many discs $B_{\rho_1}(P)$, $P\in{{\mathcal S}}\cup{{\mathcal I}}$, at most $2R/\rho_1+2$ per period, cover 
$\{z\in\overline{\Omega}:\operatorname{dist}(z,{{\mathcal S}}\cup{{\mathcal I}})<\rho_1/2\}$.

The strip $\{-h<y<-h+\delta\}$ is connected, is adjacent to ${{\mathcal B}}$ and disjoint from ${{\mathcal S}}\cup{{\mathcal I}}$, since both curves are at distance at least $\delta$ from ${{\mathcal B}}$. The strip therefore lies in ${\Omega}^-$; there the odd reflection of 
$\psi^--\psi^-|_{{\mathcal B}}-\frac12{{\Omega_2}}(y+h)^2$
 across ${{\mathcal B}}$ is harmonic in $\{|y+h|<\delta\}$, and interior estimates bound $\psi^-$ in $C^{4,\alpha}$ on $\{-h\le y\le-h+\delta/2\}$.

At the remaining points, at distance at least $\min\{\rho_1,\delta\}/2$ from ${{\mathcal S}}\cup{{\mathcal I}}\cup{{\mathcal B}}$, interior estimates for the harmonic functions 
$\psi^\pm-\frac12{\Omega}_{1,2}y^2$
 give bounds. Together with the surface and interface gains above this bounds $\psi^\pm$ in $C^{4,\alpha}(\overline{{\Omega}^\pm})$ and the arclength parametrizations of ${{\mathcal S}}$ and ${{\mathcal I}}$ with their derivatives bounded in $C^{3,\alpha}$.

\end{proof}

The hodograph at a non-stagnant free surface with constant vorticity
uses $p=\psi$, as in \cite[\S4]{CV}. At the interface, applying the local
hodograph method of \cite[\S2, Theorem~3.1$'$]{KNS} with $p=\psi-m_1$
gives qualitative interface regularity.

\section{Dispersion relations and local bifurcation}\label{sec_linear}

In this section we determine the laminar flows at which bifurcation can occur. The Fourier equations give the dispersion relation for each mode. We then characterize the roots with positive profiles and establish local bifurcation at a simple root.

\subsection{Fourier modes and the zero-eigenvalue problem}\label{sec_modes}

\subsubsection{The physical kernel and the dispersion function}

Let us first determine the kernel of \(\partial_vF(\lambda,0)\) along the admissible laminar family.

At \(w_0(\lambda)\), the flattening maps and horizontal relative velocity in \autoref{lem_conj} are
\[
Y_0^+=-d+td,\qquad
Y_0^-=-h+t(h-d),\qquad
\psi_y=W(y)\quad\text{from \eqref{eq_w}}.
\]
Since $\psi_x=0$ and $\psi_y=W(y)$,
$\nabla\psi\cdot\nabla\psi_y=W(y)W'(y).$
Consequently, on the free surface \({{\mathcal S}}=\{y=0\}\),
\[
2\nabla\psi\cdot\nabla{\widehat\Psi}=2s{\widehat\Psi}_y,
\qquad
\beta=-\frac{2s{{\Omega_1}}+2g}{s}.
\]
On the flat interface ${{\mathcal I}}=\{y=-d\}$, the corresponding jump coefficient is
\[
\kappa=\frac{{{\Omega_1}}-{{\Omega_2}}}{{{W_{\!I}}}}.
\]

Hence at an admissible laminar point \((\lambda,0)\in{\mathcal O}\), a vector
\(({\widehat\Psi}^+,{\widehat\Psi}^-,{\widehat m}_1,{\widehat Q})\in\widetilde X\) belongs to the kernel of
\(\widetilde L\) if and only if
\begin{equation}\label{eq_kernelprob}
\begin{aligned}
&{\Delta}{\widehat\Psi}^\pm=0 &&\text{in the flat layers},\\
&{\widehat\Psi}^-={\widehat m}_1 &&\text{on }y=-h,\\
&[{\widehat\Psi}]=0 &&\text{on }y=-d,\\
&s^2{\widehat\Psi}^+_y-(g+{{\Omega_1}}s){\widehat\Psi}^+=\frac{s}{2}{\widehat Q}
   &&\text{on }y=0,\\
&[{\widehat\Psi}_y]+\frac{{{\Omega_1}}-{{\Omega_2}}}{{{W_{\!I}}}}({\widehat m}_1-{\widehat\Psi})=0
   &&\text{on }y=-d,\\
&\langle{\widehat\Psi}^+(\cdot,0)\rangle=0,\qquad
 \langle{\widehat\Psi}(\cdot,-d)\rangle={\widehat m}_1 .
\end{aligned}
\end{equation}
At a laminar point we have \(s=-\sqrt\lambda\neq 0\). We obtain the condition on \(y=0\) by multiplying the surface component \(2s{\widehat\Psi}_y+\beta{\widehat\Psi}^+-{\widehat Q}=0\) of \eqref{eq_ltilde} by \(s/2\).
The mean conditions in \eqref{eq_xtilde} give the last line of \eqref{eq_kernelprob} because \(\psi_y\) equals the nonzero constants \(s\) and \({{W_{\!I}}}\) on \({{\mathcal S}}\) and \({{\mathcal I}}\), respectively.

The \(x\)-independent laminar coefficients make both \(\widetilde L\) and \(\partial_vF(\lambda,0)\) preserve Fourier modes, in the kernel and in the range equation. Let \({\widehat\Psi}^\pm_n(y)\) denote the coefficient of \(\cos nx\) in \({\widehat\Psi}^\pm(x,y)\). The scalar variations \({\widehat m}_1,{\widehat Q}\) belong to the zero mode.

We compute the dispersion function from \eqref{eq_kernelprob} at an arbitrary wavenumber \(k>0\). On the \(2\pi\)-periodic graph space, the \(n\)th mode corresponds to \(k=n\). For the single-mode calculation set
${\widehat\Psi}^\pm=\varphi(y)\cos kx,\qquad {\widehat m}_1={\widehat Q}=0.$

And, writing
\[
{\widehat\eta}=\eta_1\cos kx,\qquad
{\widehat\zeta}=\zeta_1\cos kx
\]
for the corresponding shape perturbations, the inverse of \(\mathcal C\) in
\autoref{lem_conj} gives
\begin{equation}\label{eq_etazeta}
\eta_1=-\frac{\varphi(0)}{s},\qquad
\zeta_1=-\frac{\varphi(-d)}{{{W_{\!I}}}}.
\end{equation}
The problem \eqref{eq_kernelprob} becomes
\begin{align}
&\varphi''=k^2\varphi
\qquad \text{on }(-h,-d)\text{ and on }(-d,0),\label{eq_mode_ode}\\
&\varphi(-h)=0,\label{eq_mode_bed}\\
&[\varphi]=0
\qquad \text{at }y=-d,\label{eq_mode_cont}\\
&[\varphi']+({{\Omega_1}}-{{\Omega_2}})\zeta_1=0,\label{eq_mode_jump}\\
&s\varphi'(0)+(g+{{\Omega_1}}s)\eta_1=0.\label{eq_mode_surf}
\end{align}
The interface and surface conditions follow from \eqref{eq_kernelprob} by substituting \eqref{eq_etazeta} and \({\widehat Q}=0\), and dividing the surface equation by \(s\). In \eqref{eq_mode_jump}, the term \(({{\Omega_1}}-{{\Omega_2}})\zeta_1=[W']\zeta_1\) is the contribution of the interface displacement to the jump of \(\psi_y\), with \([W']\) from \eqref{eq_w}.

With the notations
\begin{equation}\label{eq_sc}
S_2=\sinh k(h-d),\qquad C_2=\cosh k(h-d),
\end{equation}
solving \eqref{eq_mode_ode}--\eqref{eq_mode_jump} gives, up to a multiplicative constant,
\begin{equation}\label{eq_f}
\begin{aligned}
&\varphi(y)=
\begin{cases}
\sinh k(y+h), & -h<y<-d,\\
S_2\cosh k(y+d)+\widehat B\sinh k(y+d), & -d<y<0,
\end{cases}\\
&\widehat B=C_2+\frac{({{\Omega_1}}-{{\Omega_2}})S_2}{k{{W_{\!I}}}}.
\end{aligned}
\end{equation}
The surface condition \eqref{eq_mode_surf} then reads \(D_k(\lambda)=0\), where
\begin{equation}\label{eq_d}
\begin{aligned}
D_k(\lambda)
&:=s^2\varphi'(0)-(g+{{\Omega_1}} s)\varphi(0)\\
&= s^2k\left[S_2\sinh kd+\widehat B\cosh kd\right] 
-(g+{{\Omega_1}}s)\left[S_2\cosh kd+\widehat B\sinh kd\right].
\end{aligned}
\end{equation}
We call \(D_k\) the \emph{dispersion function} of the mode \(\cos kx\)
and \(D_k(\lambda)=0\) the \emph{dispersion relation}. Here \(s\) and \({{W_{\!I}}}\) depend on \(\lambda\) through \eqref{eq_lampar}. 

Under the normalization \eqref{eq_f}, the shape amplitudes \eqref{eq_etazeta} are
\begin{equation}\label{eq_amps}
\eta_1=-\frac{\varphi(0)}{s},\qquad
\zeta_1=-\frac{S_2}{{{W_{\!I}}}}.
\end{equation}

\begin{proposition}\label{prop_modes}
Let \((\lambda,0)\in{\mathcal O}\). For \(n\ge1\), let \(\varphi_n^+\) and
\(\varphi_n^-\) be the restrictions of the profile \eqref{eq_f} at \(k=n\) to
\((-d,0)\) and \((-h,-d)\), respectively. With the isomorphism \(\mathcal C\) of \autoref{lem_conj},
\[
\begin{aligned}
\ker\partial_vF(\lambda,0)
=\operatorname{span}\Bigl\{&
\mathcal C^{-1}\bigl(\varphi_n^+(y)\cos nx,\,
\varphi_n^-(y)\cos nx,\,0,\,0\bigr):
n\ge1,\quad D_n(\lambda)=0\Bigr\}.
\end{aligned}
\]
\end{proposition}
\begin{proof}
By \autoref{lem_conj}, it suffices to determine the Fourier modes of the physical kernel \eqref{eq_kernelprob}.
The mean constraints and the bed and interface traces in \eqref{eq_kernelprob} give
\[
{\widehat\Psi}_0^+(0)=0,\qquad
{\widehat\Psi}_0^+(-d)={\widehat\Psi}_0^-(-d)={\widehat\Psi}_0^-(-h)={\widehat m}_1.
\]
From the bulk equations we obtain \(({\widehat\Psi}_0^\pm)''=0\). The lower profile is affine with equal endpoint values. Since the common interface trace equals \({\widehat m}_1\), the jump condition gives
$({\widehat\Psi}_0^+)'(-d)=({\widehat\Psi}_0^-)'(-d)=0.$
The zero-mode upper profile \({\widehat\Psi}_0^+\) is constant and vanishes at \(y=0\). Its trace at \(y=-d\) gives \({\widehat m}_1=0\), so \({\widehat\Psi}_0^-\) also vanishes. The surface equation gives \({\widehat Q}=0\).

For \(n\ge1\), the mean conditions hold automatically and the scalar variations vanish in this mode. The bed condition fixes the lower profile up to a factor, and continuity and the derivative jump determine both upper Cauchy data at \(y=-d\). Thus the mode is \(c_n\varphi_n\) by \eqref{eq_f}, with the surface condition \(c_nD_n(\lambda)=0\) by \eqref{eq_d}.

The modes with \(D_n(\lambda)=0\) give nonzero kernel vectors under \(\mathcal C^{-1}\), linearly independent by their distinct Fourier frequencies. Since \(F(\lambda,0)=0\) by \autoref{prop_operator}, \autoref{cor_fredall} gives a finite-dimensional kernel, so only finitely many such modes occur. Uniqueness of Fourier coefficients gives the claimed span.
\end{proof}

\begin{remark}\label{rem_fourier_normalization}
Our mode is normalized at the bed as in \eqref{eq_f}. In the notation of Martin and Matioc \cite[(2.8), (2.18)]{MM16}, set $d_2=d$, $d_1=h-d$, $\gamma_2=\Omega_1$, $\gamma_1=\Omega_2$, and $R_n=k$ for the physical wavenumber. Their signed surface speed is $\Lambda=s=-\sqrt\lambda$, and $\Lambda-\gamma_2d_2=W_I$. Multiplying \eqref{eq_dsplit} by $W_I/(k\cosh(kh))$, we obtain their monic cubic in $\Lambda$. The two dispersion equations are equivalent for $W_I\ne0$. For the capillary--gravity formulation see also \cite{MM}.
\end{remark}

The formulas \eqref{eq_f}--\eqref{eq_amps} hold for arbitrary $k>0$. The analysis characterizes the roots of $D_k(\lambda)=0$ with positive normalized profiles \eqref{eq_f}, determines when positive roots occur in pairs in $\{s<0<{{W_{\!I}}}\}$, and excludes higher-mode resonance at positive first-mode roots.

Applying the hyperbolic addition formulas to \eqref{eq_d}, with \(S_2\) and \(C_2\) from \eqref{eq_sc}, we obtain
\begin{equation}\label{eq_dsplit}
\begin{aligned}
D_k(\lambda)
&=s^2k\cosh kh-(g+{{\Omega_1}}s)\sinh kh\\
&\quad+\frac{({{\Omega_1}}-{{\Omega_2}})S_2}{{{W_{\!I}}}}
\left[s^2\cosh kd-\frac{g+{{\Omega_1}}s}{k}\sinh kd\right].
\end{aligned}
\end{equation}
The first line of \eqref{eq_dsplit} is the single-layer dispersion function of \cite{CSV} at depth $h$ and vorticity ${{\Omega_1}}$. The term proportional to ${{\Omega_1}}-{{\Omega_2}}$ is the two-layer contribution.

The only possible singularity is ${{W_{\!I}}}=0$, the critical level at the interface. For ${{\Omega_1}}<0$ this occurs at $\lambda={{\Omega_1}}^2d^2$. With $s={{W_{\!I}}}+{{\Omega_1}}d$, the product ${{W_{\!I}}}D_k$ is a cubic in ${{W_{\!I}}}$ with leading coefficient $k\cosh kh$, so $D_k$ has at most three roots.

\subsection{Positive roots}\label{sub_04_positive_roots}

For the normalized profile $\varphi$ defined in \eqref{eq_f} we first find the condition for positivity on $(-h,0]$. We then ask for which parameters the corresponding roots of $D_k(\lambda)=0$ occur in pairs.

\begin{lemma}\label{lem_pos}
Let $\varphi$ be the normalized profile \eqref{eq_f}.
\begin{enumerate}[label=\textup{(\roman*)},leftmargin=2.3em]
\item
\[
\varphi>0\ \text{on }(-h,0]
\quad\Longleftrightarrow\quad \varphi(0)>0.
\]
\item If ${{W_{\!I}}}>0$ and ${{\Omega_2}}>{{\Omega_1}}$, then $\varphi>0$ on $(-h,0]$ if and only if
\begin{equation}\label{eq_wc}
{{W_{\!I}}}>{{W_{\!c}}}:=\frac{{{\Omega_2}}-{{\Omega_1}}}{k\big(\coth kd+\coth k(h-d)\big)}.
\end{equation}
\end{enumerate}
\end{lemma}

\begin{proof}
The lower profile in \eqref{eq_f} is strictly positive above the bed. In the upper layer the solution of $u''=k^2u$ interpolates its values at $-d,0$ with the two positive weights $\sinh(-ky)/\sinh(kd)$ and $\sinh(k(y+d))/\sinh(kd)$. This proves (i). In (ii) evaluation at zero gives $\varphi(0)=\sinh(k(h-d))\sinh(kd)[\coth(kd)+\coth(k(h-d))+(\Omega_1-\Omega_2)/(kW_I)]$. Its sign gives \eqref{eq_wc}.
\end{proof}

The next lemma rewrites the surface condition \eqref{eq_d} on the component $\{s<0<{{W_{\!I}}}\}$. We fix ${{\Omega_1}}<0$ and set $t_k=\tanh kd$, $K=\coth k(h-d)$ and $\sigma=\sqrt\lambda=-s$. For $\lambda>0$, define
\begin{equation}\label{eq_dispersion_q}
q(\lambda)=\frac{k\lambda}{g-{{\Omega_1}}\sqrt\lambda}-t_k
=\frac{ks^2}{g+{{\Omega_1}}s}-t_k,
\end{equation}
so that $ks^2=(g+{{\Omega_1}}s)(q+t_k)$. Differentiating in $\sigma$, we obtain
\begin{equation}\label{eq_qsigma}
\begin{aligned}
q_\sigma
&=\frac{2k\sigma(g-{{\Omega_1}}\sigma)+k{{\Omega_1}}\sigma^2}{(g-{{\Omega_1}}\sigma)^2}
=\frac{k\sigma(2g-{{\Omega_1}}\sigma)}{(g-{{\Omega_1}}\sigma)^2}>0.
\end{aligned}
\end{equation}
Thus $q$ is strictly increasing for $\sigma>0$. Its definition gives
\[
\lim_{\lambda\downarrow0}q(\lambda)=-t_k<0,\qquad
\lim_{\lambda\to+\infty}q(\lambda)=+\infty.
\]
Hence $q$ has a unique positive zero $\lambda_q$. The interval
\begin{equation}\label{eq_positive_root_interval}
\mathcal J=\{\lambda:\ \lambda_q<\lambda<{{\Omega_1}}^2d^2\}
\end{equation}
is the set of $\lambda>0$ in this component on which $q>0$: ${{W_{\!I}}}=-\sqrt\lambda-{{\Omega_1}}d>0$ for $\lambda<{{\Omega_1}}^2d^2$, and \eqref{eq_qsigma} gives $q>0$ for $\lambda>\lambda_q$. For $\lambda>0$ with $q(\lambda)\ne0$, set
\begin{equation}\label{eq_dispersion_rho}
\rho(\lambda)=\frac{1-t_k^2}{q(\lambda)}-t_k,\qquad f(\lambda)=k{{W_{\!I}}}\rho(\lambda).
\end{equation}
On the component $\{s<0<{{W_{\!I}}}\}$, write $\varkappa_k:=({{\Omega_1}}-{{\Omega_2}})/(k{{W_{\!I}}})$. For $q(\lambda)\ne0$, define the reduced dispersion function by
\begin{equation}\label{eq_dispersion_reduction}
\begin{aligned}
\mathcal E(\lambda)
&=-k{{W_{\!I}}}(K+\varkappa_k-\rho)
={{\Omega_2}}-{{\Omega_1}}-kK{{W_{\!I}}}+f(\lambda).
\end{aligned}
\end{equation}

\begin{lemma}\label{lem_convex}
Let ${{\Omega_1}}<0$ and use $q$ from \eqref{eq_dispersion_q}, the interval $\mathcal J$ from \eqref{eq_positive_root_interval}, the functions $\rho,f$ from \eqref{eq_dispersion_rho}, and the reduced dispersion function $\mathcal E$ in \eqref{eq_dispersion_reduction}.
\begin{enumerate}[label=\textup{(\roman*)},leftmargin=2.3em]
\item On $\mathcal J$, the dispersion function factors as
\[
{{W_{\!I}}} D_k(\lambda)=-k^{-1}\sinh k(h-d)\cosh kd\,(g+{{\Omega_1}}s)\,q(\lambda)\,\mathcal E(\lambda).
\]
Thus $D_k=-c\,\mathcal E$ on $\mathcal J$ with $c>0$ real-analytic.
\item For a root $\lambda$ of $D_k$ with $s<0<{{W_{\!I}}}$,
\[
\varphi>0\ \text{on }(-h,0]
\quad\Longleftrightarrow\quad q(\lambda)>0
\quad\Longleftrightarrow\quad \lambda\in\mathcal J.
\]
The positive roots in this component are therefore exactly the zeros of
$\mathcal E$ in $\mathcal J$.
\item Both $f$ and $\mathcal E$ are strictly convex in
$\sigma=\sqrt\lambda$ on $\mathcal J$; equivalently, $\mathcal E$ is strictly
convex in ${{W_{\!I}}}$.
At a double zero,
\[
\partial_\sigma^2\mathcal E>0,\qquad
\partial_\lambda\mathcal E=0<\partial_\lambda^2\mathcal E.
\]
\item If $\mathcal J\ne\emptyset$, the endpoint limits are
\[
\begin{aligned}
\mathcal E(\lambda)&\to+\infty
&&\text{as }\lambda\downarrow\lambda_q,\\
\mathcal E(\lambda)&\to{{\Omega_2}}-{{\Omega_1}}
&&\text{as }\lambda\uparrow{{\Omega_1}}^2d^2
\quad({{W_{\!I}}}\downarrow0).
\end{aligned}
\]
\end{enumerate}
\end{lemma}
\begin{proof}
On \(s<0<{{W_{\!I}}}\), we have \(g+{{\Omega_1}}s>0\). Substituting \eqref{eq_sc}, \eqref{eq_f} and \(ks^2=(g+{{\Omega_1}}s)(q+t_k)\) in \eqref{eq_d}, we obtain
\[
\begin{aligned}
\frac{D_k}{S_2\cosh kd\,(g+{{\Omega_1}}s)}
&=(q+t_k)(t_k+K+\varkappa_k)-1-(K+\varkappa_k)t_k\\
&=q(K+\varkappa_k+t_k)-(1-t_k^2)
=q(K+\varkappa_k-\rho) \qquad(q\ne0).
\end{aligned}
\]
Since \(\mathcal E=-k{{W_{\!I}}}(K+\varkappa_k-\rho)\), this gives (i), with
\[
c=\frac{S_2\cosh kd\,(g+{{\Omega_1}}s)\,q}{k{{W_{\!I}}}}>0
\qquad\text{on }\mathcal J.
\]
All factors are analytic there.

At a root with \(s<0<{{W_{\!I}}}\), the second line of the preceding calculation excludes \(q=0\) because \(1-t_k^2>0\). The third line then gives \(K+\varkappa_k=\rho\). Substitution in \eqref{eq_f} yields
\[
\begin{aligned}
\varphi(0)&=S_2\cosh kd+\widehat B\sinh kd
=S_2\cosh kd\bigl(1+(K+\varkappa_k)t_k\bigr)
=S_2\cosh kd\,(1-t_k^2)\frac{q+t_k}{q}.
\end{aligned}
\]
Since \(q+t_k=k\lambda/(g+{{\Omega_1}}s)>0\), \(\varphi(0)\) has the sign of \(q\). By \autoref{lem_pos}, positivity is equivalent to \(q>0\), which in this component means \(\lambda\in\mathcal J\) by \eqref{eq_positive_root_interval}. Part~(i) identifies these roots with the zeros of \(\mathcal E\) on \(\mathcal J\).

Write primes for differentiation in $\sigma=\sqrt\lambda$. By \eqref{eq_dispersion_rho}, only the term $k(1-t_k^2){W_{\!I}}/q$ contributes to $f''$, since ${W_{\!I}}=-\sigma-{\Omega_1}d$. Differentiating \eqref{eq_qsigma}, we obtain
\[
q''=\frac{2kg^2}{(g-{\Omega_1}\sigma)^3}>0.
\]
On $\mathcal J$, the bound $0<q<k\sigma^2/(g-{\Omega_1}\sigma)$ yields
\[
\begin{aligned}
qq''&<\frac{2k^2g^2\sigma^2}{(g-{\Omega_1}\sigma)^4}
<\frac{2k^2\sigma^2(2g-{\Omega_1}\sigma)^2}{(g-{\Omega_1}\sigma)^4}
=2q'^2,
\end{aligned}
\]
where $0<g<2g-{\Omega_1}\sigma$. Using ${W_{\!I}}'=-1$, we obtain
\[
f''=k(1-t_k^2)
\left({W_{\!I}}\frac{2q'^2-qq''}{q^3}+\frac{2q'}{q^2}\right)>0
\qquad\text{on }\mathcal J.
\]
The other terms in $\mathcal E$ in \eqref{eq_dispersion_reduction} are affine in $\sigma$, so $\mathcal E$ is strictly convex in $\sigma$ and hence in ${W_{\!I}}$.

At a double zero $\sigma_0$, $\mathcal E''(\sigma_0)=f''(\sigma_0)>0$. Since $\partial_\lambda=(2\sigma)^{-1}\partial_\sigma$, at a double zero $\partial_\lambda\mathcal E=0$ and $\partial^2_\lambda\mathcal E=\mathcal E''/(4\lambda)>0$; the multiplicity of a zero is the same in $\sigma$, in ${{W_{\!I}}}$ and in $\lambda$.

If \(\mathcal J\ne\emptyset\), then \(q\downarrow0\) and \({{W_{\!I}}}\to-\sqrt{\lambda_q}-{{\Omega_1}}d>0\) at the left endpoint. Thus the term \(k(1-t_k^2){{W_{\!I}}}/q\) in \(\mathcal E\) diverges to \(+\infty\) while the other terms stay bounded. At the right endpoint, \({{W_{\!I}}}\downarrow0\) and \(q\to q({{\Omega_1}}^2d^2)>0\), so \(f=k{{W_{\!I}}}\rho\to0\) and \(\mathcal E\to{{\Omega_2}}-{{\Omega_1}}\).
\end{proof}

\begin{proposition}\label{prop_pairing}
Let ${{\Omega_1}}<\min(0,{{\Omega_2}})$. In the component
$\{s<0<{{W_{\!I}}}\}$, equivalently $0<\lambda<{{\Omega_1}}^2d^2$,
set
$\lambda_c:=\big((-{{\Omega_1}}d-{{W_{\!c}}})_+\big)^2,$
where ${{W_{\!c}}}$ is given by \eqref{eq_wc}.
Then a root of $D_k(\lambda)=0$ satisfies $\varphi>0$ on $(-h,0]$ if and only if $0<\lambda<\lambda_c$, equivalently
$0<\sqrt\lambda<-{{\Omega_1}}d-{{W_{\!c}}}.$
The positive roots in this component are the zeros of $\mathcal E$ on $\mathcal J$.

Write the dependence on ${{\Omega_2}}$ as
\[
\mathcal E=\mathcal E_0+{{\Omega_2}},\qquad
\mathcal E_0=-{{\Omega_1}}-kK{{W_{\!I}}}+f,
\]
where $\mathcal E_0$ is independent of ${{\Omega_2}}$. Then the minimum of $\mathcal E_0$ is attained
in $\mathcal J$ exactly when
\[
\mathcal J\neq\emptyset,\qquad
\coth k(h-d)>\rho({{\Omega_1}}^2d^2).
\]
We define the threshold by
\[
{{\Omega_2}}^\ast=
\begin{cases}
-\min_{\mathcal J}\mathcal E_0,&\text{minimum in }\mathcal J,\\
{{\Omega_1}},&\text{otherwise}.
\end{cases}
\]
Then the number of distinct positive roots is
\[
\#\{\lambda\in\mathcal J:\mathcal E(\lambda)=0\}=
\begin{cases}
2,&{{\Omega_2}}<{{\Omega_2}}^\ast,\\
1,&{{\Omega_2}}={{\Omega_2}}^\ast,\\
0,&{{\Omega_2}}>{{\Omega_2}}^\ast.
\end{cases}
\]
The two roots in the first case are simple. At equality, the unique root is double.
\end{proposition}
\begin{proof}
By \autoref{lem_convex}(ii), the positive roots in \(\{s<0<{{W_{\!I}}}\}\) are the zeros of \(\mathcal E\) in \(\mathcal J\).
Since \({{\Omega_2}}>{{\Omega_1}}\), \autoref{lem_pos} gives \({{W_{\!I}}}>{{W_{\!c}}}\) as the criterion for a root in this component to have a positive profile. Substituting \({{W_{\!I}}}=-\sqrt\lambda-{{\Omega_1}}d\), we obtain the stated interval, which is empty when \({{W_{\!c}}}\ge-{{\Omega_1}}d\). If $\mathcal J=\emptyset$, there are no positive roots and the definition gives ${{\Omega_2}}^\ast={{\Omega_1}}<{{\Omega_2}}$, as required. Assume henceforth that $\mathcal J\ne\emptyset$. The multiplicity of a zero of $D_k$ in $\mathcal J$ equals its multiplicity as a zero of $\mathcal E$, since $D_k=-c\,\mathcal E$ with $c>0$ analytic on $\mathcal J$ (\autoref{lem_convex}(i)). In the rest of the proof we regard functions on $\mathcal J$ as functions of
$\sigma=\sqrt\lambda\in(\sqrt{\lambda_q},-{{\Omega_1}}d)$. By \autoref{lem_convex}(iii),(iv), $\mathcal E_0=\mathcal E-{{\Omega_2}}$ is strictly convex on this interval, tends to $+\infty$ at its left end and to $-{{\Omega_1}}$ at its right end $\sigma=-{{\Omega_1}}d$, where ${{W_{\!I}}}=0$.

The derivative of $\mathcal E_0$ is strictly increasing and is negative near the left endpoint because $\mathcal E_0$ diverges there. At the right endpoint, $q({{\Omega_1}}^2d^2)>0$, so $\rho$ is smooth there. Using $\partial_\sigma{{W_{\!I}}}=-1$,
\begin{equation*}
\begin{aligned}
\partial_\sigma\mathcal E_0
&=k(K-\rho)+k{W_{\!I}}\,\partial_\sigma\rho
\longrightarrow k\bigl(K-\rho({\Omega_1}^2d^2)\bigr)
\qquad\text{as }\sigma\uparrow-{\Omega_1}d .
\end{aligned}
\end{equation*}
Thus the derivative vanishes in $\mathcal J$ exactly when this limit is positive, which is the criterion for an attained minimum. If there is no interior minimum, $\mathcal E_0$ decreases to $-{{\Omega_1}}$ at the right endpoint, so $\mathcal E_0>-{{\Omega_1}}>-{{\Omega_2}}$. There is no root, and the threshold is ${{\Omega_2}}^\ast={{\Omega_1}}$. If the minimum $-{{\Omega_2}}^\ast$ is attained at $\sigma_m$, it lies below $-{{\Omega_1}}$. The function strictly decreases from $+\infty$ to this minimum and strictly increases from it to $-{{\Omega_1}}$. For ${{\Omega_1}}<{{\Omega_2}}<{{\Omega_2}}^\ast$, the level $-{{\Omega_2}}$ occurs once on each side of $\sigma_m$, and both roots are simple since $\partial_\sigma\mathcal E_0\ne0$ there. At ${{\Omega_2}}={{\Omega_2}}^\ast$, the sole root is double by \autoref{lem_convex}(iii). For ${{\Omega_2}}>{{\Omega_2}}^\ast$ there are no roots.
\end{proof}

\begin{lemma}\label{lem_depth}
Suppose that $D_k(\lambda)=0$, its normalized profile satisfies $\varphi>0$ on $(-h,0]$, and $s<0<{{W_{\!I}}}$. Then ${{\Omega_1}}<0$.
The upper critical layer lies at depth
$\delta_c=s/{{\Omega_1}}{}=\sqrt\lambda/|{{\Omega_1}}|$ below the mean surface, with
\begin{equation}\label{eq_depth}
\delta_c>\frac{\tanh kd}{k}\Big(1+\frac{g}{{{\Omega_1}}s}\Big)>\frac{\tanh kd}{k}.
\end{equation}
\end{lemma}
\begin{proof}
Since $s<0<W_I$ and $s=W_I+\Omega_1d$, one has $\Omega_1<0$. The root has $q>0$ by \autoref{lem_convex}(ii). Thus $ks^2>(g+\Omega_1s)\tanh(kd)$; dividing by the positive number $k\Omega_1s$ gives \eqref{eq_depth}.
\end{proof}

\begin{corollary}\label{cor_fold_depth}
We set $a=-\Omega_1>0$ and $H_-=h-d$. For $\Omega_2>\Omega_1$ the threshold in \autoref{prop_pairing} has an interior minimum exactly when
\[
 a^2d\,[kd-\tanh(kh)]>g\tanh(kh).
\]
Here, at fixed $a,d,g,k$ the critical lower-layer shear satisfies
$\partial_{H_-}\Omega_2^\ast=-k^2W_{I,\ast}\operatorname{csch}^2(kH_-)<0$.
\end{corollary}
\begin{proof}
We write $t=\tanh(kd)$ and $K=\coth(kH_-)$. At the right endpoint $\sigma=ad$, the minimum criterion requires $q>0$ and $K>(1-t^2)/q-t$. Together these are equivalent to $k(ad)^2/(g+a^2d)>t+(1-t^2)/(K+t)=\tanh(kh)$. Rearrangement gives the stated inequality. At the unique nondegenerate maximum of $\Omega_1+kW_I(K-\rho)$ its derivative in $\sigma$ vanishes. Differentiation in $H_-$ therefore gives $kW_{I,\ast}\partial_{H_-}K$, as asserted.
\end{proof}

At $s=0$, the left end of a nonempty window $0<\lambda<\lambda_c$, \eqref{eq_lampar} gives ${{W_{\!I}}}=-{{\Omega_1}}d>{{W_{\!c}}}$. Substituting $s=0$ into \eqref{eq_d}, we obtain $D_k(0)=-g\varphi(0)$. By \autoref{lem_pos}(ii), $\varphi(0)>0$, so
\begin{equation}\label{eq_dzero}
D_k(0)=-g\varphi(0)<0.
\end{equation}

\begin{lemma}\label{lem_sl}
Suppose that $D_k(\lambda)=0$, the normalized profile $\varphi$ of \eqref{eq_f} is positive on $(-h,0]$, and ${{W_{\!I}}}(\lambda)\neq0$.
Then all higher modes satisfy
$D_{k'}(\lambda)>0\qquad(k'>k).$
\end{lemma}
\begin{proof}
We fix \(\lambda\) and compare \(\varphi=\varphi_k\) with \(\varphi_{k'}\) by means of their Wronskian. Related ground-state arguments can be found in \cite[Lemma~3.7]{CS} and \cite{GWYZ}. By \eqref{eq_f} both profiles have the same transmission coefficient in \([\varphi_j']=(({{\Omega_1}}-{{\Omega_2}})/{{W_{\!I}}})\varphi_j(-d)\) for \(j\in\{k,k'\}\). It follows that their Wronskian satisfies
\[
\begin{aligned}
M&=\varphi\varphi_{k'}'-\varphi'\varphi_{k'},&
M(-h)&=0,\qquad [M]=0,\\
M'&=(k'^2-k^2)\varphi\varphi_{k'}
&&\text{on }(-h,-d)\cup(-d,0).
\end{aligned}
\]
The profile \(\varphi_{k'}\) is positive on \((-h,-d]\). At a first zero \(y_0\in(-d,0]\), integration gives \(M(y_0)>0\), whereas \(M(y_0)=\varphi(y_0)\varphi_{k'}'(y_0)\le0\). Thus \(\varphi_{k'}>0\) on \((-h,0]\), and integration gives \(M(0)>0\). Since \(D_k(\lambda)=0\), \eqref{eq_d} gives \(s^2\varphi'(0)=(g+{{\Omega_1}}s)\varphi(0)\), with \(s\ne0\) because \(g\varphi(0)>0\). Hence
\[
D_{k'}(\lambda)=s^2\varphi_{k'}'(0)-(g+{{\Omega_1}}s)\varphi_{k'}(0)
=\frac{s^2}{\varphi(0)}M(0)>0.
\]
\end{proof}

Each parameter set in \autoref{tab_roots} has two positive roots of \(D_1(\lambda)=0\) in \(\{s<0<{{W_{\!I}}}\}\). In this component, \eqref{eq_amps} gives
\begin{equation}\label{eq_mode_amplitude_product}
\eta_1\zeta_1=\frac{S_2\varphi(0)}{s{{W_{\!I}}}}.
\end{equation}
At both roots, $\varphi(0)>0$ and $s{{W_{\!I}}}<0$, so the surface and interface amplitudes have opposite signs. We call the root with larger \(|\eta_1/\zeta_1|\) the \emph{surface-mode root} and the one with smaller \(|\eta_1/\zeta_1|\) the \emph{interface-mode root}.

\begin{table}[htbp]
\centering\small
\caption{Roots of the dispersion expression \eqref{eq_d} as a function of $s$, with $g=k=1$. The columns CL upper and CL lower give the depth below the surface and the height above the bed of the respective laminar critical layers. The first root of each case has $s>0$, so it is not of the form \eqref{eq_lampar}; its entry in the $\lambda$ column is $s^2$, in parentheses.}\label{tab_roots}
\setlength{\tabcolsep}{5pt}
\begin{tabular}{@{}l c c c c c c c c@{}}
\toprule
case $({{\Omega_1}},{{\Omega_2}},h,d)$ & $\lambda$ & $s$ & ${{W_{\!I}}}$ & $m_1$ & $m_2$ & CL upper & CL lower & $\eta_1/\zeta_1$\\
\midrule
A $(-1,1,5,3)$ & $(0.382)$ & $0.618$ & $3.6183$ & $-6.355$ & $-5.237$ & -- & -- & $+87.3$\\
 & $2.6651$ & $-1.6325$ & $1.3675$ & $0.398$ & $-0.735$ & $1.63$ & $0.63$ & $-4.87$\\
 & $4.0067$ & $-2.0017$ & $0.9983$ & $1.505$ & $0.003$ & $2.00$ & $1.00$ & $-0.195$\\
\midrule
B $(-2,2,5,2.5)$ & $(0.172)$ & $0.415$ & $5.4147$ & $-7.287$ & $-7.287$ & -- & -- & $+101.8$\\
 & $6.4037$ & $-2.5306$ & $2.4694$ & $0.076$ & $0.076$ & $1.27$ & $1.27$ & $-2.41$\\
 & $8.3182$ & $-2.8841$ & $2.1159$ & $0.960$ & $0.960$ & $1.44$ & $1.44$ & $-0.607$\\
\midrule
C $(-1,1,4,3)$ & $(0.382)$ & $0.618$ & $3.6181$ & $-6.354$ & $-3.118$ & -- & -- & $+103.5$\\
 & $2.6432$ & $-1.6258$ & $1.3742$ & $0.377$ & $-0.874$ & $1.63$ & (none) & $-7.30$\\
 & $4.5164$ & $-2.1252$ & $0.8748$ & $1.876$ & $-0.375$ & $2.13$ & $0.13$ & $-0.131$\\
\bottomrule
\end{tabular}
\end{table}

A scan of \(5\times6\times6\times7\) parameter sets \((h,d/h,{{\Omega_1}},{{\Omega_2}})\) with \({{\Omega_1}}<0<{{\Omega_2}}\) found four pairs of counts in \(\{s<0<{{W_{\!I}}}\}\). The first entry counts roots whose profile \(\varphi\) in \eqref{eq_f} is positive on \((-h,0]\), and the second counts roots whose profile changes sign,
$(0,0),\qquad(0,1),\qquad(0,2),\qquad(2,0).$
The first root listed in \autoref{tab_roots} for each example also has \(\varphi>0\) on \((-h,0]\). It has \(s>0\), or \(u(0)>c\), and no critical layer, so it lies outside the family \eqref{eq_lampar} and the pairing result.

\subsection{Local bifurcation}\label{sec_local}

We next construct the small-amplitude branch at a simple root and determine the action of half-period translation on its two halves.

At an admissible laminar point $(\lambda^\ast,0)$ we denote the linearization of the graph residual $F$ by
\begin{equation}\label{eq_graph_laminar_linearization}
L_\ast:=\partial_vF(\lambda^\ast,0).
\end{equation}

We recall the analytic form of the simple-eigenvalue bifurcation theorem of Crandall and Rabinowitz \cite{CR}.
\begin{theorem}\label{thm_cr}
For an analytic equation $F(\lambda,0)=0$, assume at a zero that $L=F_v$ is Fredholm of index zero, $\ker L=\mathbb Rw_1$, and $F_{\lambda v}w_1\notin\operatorname{Ran}L$. For any fixed closed complement of $w_1$, all nearby nontrivial solutions lie on the unique local analytic arc with kernel coordinate $a$, tangent in the state variable to $w_1$.
\end{theorem}

The profile \(\varphi(y)\cos x\) from \eqref{eq_f} reduces transversality to the condition \(\partial_\lambda D_1(\lambda^\ast)\ne0\).

For a laminar parameter $\lambda^\ast$, the root and nonresonance conditions, and the simple-root condition, are
\begin{align}
\label{eq_h}
D_1(\lambda^\ast)=0,\qquad D_n(\lambda^\ast)\ne0\quad(n\ge2),\\
\label{eq_hsimple}
\partial_\lambda D_1(\lambda^\ast)\ne0.
\end{align}
\begin{theorem}\label{thm_trans}
Let $(\lambda^\ast,0)\in{\mathcal O}$. Suppose that \eqref{eq_h} holds.
Use the mode-one profiles $\varphi_1^\pm$ of \autoref{prop_modes} at $\lambda^\ast$
and the isomorphism $\mathcal C$ of \autoref{lem_conj}.
Then $L_\ast$ is Fredholm of index zero and
\[
\ker L_\ast=\operatorname{span}\{w_1\},\qquad
\mathcal C(\lambda^\ast)w_1
=\bigl(\varphi_1^+(y)\cos x,\ \varphi_1^-(y)\cos x,\ 0,\ 0\bigr).
\]
The surface and interface displacements of $w_1$ are
\[
{\widehat\eta}=-\frac{\varphi_1^+(0)}{s}\cos x,\qquad
{\widehat\zeta}=-\frac{\varphi_1^-(-d)}{{W_{\!I}}}\cos x.
\]
If \eqref{eq_hsimple} holds as well, then
$\partial_\lambda\partial_vF(\lambda^\ast,0)w_1\notin\operatorname{Ran} L_\ast.$
\end{theorem}

\begin{proof}
Fourier modes identify $\ker L_\ast$. The action \eqref{eq_trans_family} of $L(\lambda)$ on the analytic mode-one family then gives transversality at a simple root.
By \autoref{prop_operator}(i), $F(\lambda^\ast,0)=0$, so \autoref{cor_fredall}
gives Fredholm index zero. By \eqref{eq_h}, only $n=1$ contributes to the kernel
in \autoref{prop_modes}, giving the stated generator with its normalization.
The inverse displacement formulas in \autoref{lem_conj} give ${\widehat\eta}$ and ${\widehat\zeta}$
from the surface trace of $\varphi_1^+$ and the common interface trace
$\varphi_1^-(-d)=\varphi_1^+(-d)=S_2$ from \eqref{eq_f}, since ${\widehat m}_1=0$.

We construct an analytic family through $w_1$ whose only residual is the surface term in \eqref{eq_trans_family}. Along the laminar family, the flat layers, $Y_0^\pm$ and $T_0^\pm$ are fixed, so the target isomorphism $\widehat{\mathcal C}$ is independent of $\lambda$. The means in $\widetilde X$ reduce to the last line of \eqref{eq_kernelprob}, so $\widetilde X$ is fixed as well. Since $W(y)$ is analytic in $\lambda$ by \eqref{eq_lampar}, \eqref{eq_psihat} gives an analytic family of isomorphisms $\mathcal C(\lambda):X\to\widetilde X$. With $\widetilde L(\lambda)$ as in \eqref{eq_ltilde}, we have
\[
L(\lambda):=\partial_vF(\lambda,0)
=\widehat{\mathcal C}^{-1}\widetilde L(\lambda)\mathcal C(\lambda),
\qquad L(\lambda^\ast)=L_\ast.
\]
For $\lambda$ near $\lambda^\ast$, let $\varphi_\lambda$ be the continuous profile \eqref{eq_f} at $k=1$, with layer restrictions $\varphi_\lambda^\pm$. The coefficient $\widehat B$ is analytic where ${W_{\!I}}\ne0$, so the following vectors depend analytically on $\lambda$:
\[
u_\lambda:=\bigl(\varphi_\lambda^+(y)\cos x,\ \varphi_\lambda^-(y)\cos x,\ 0,\ 0\bigr)
\in\widetilde X,\qquad
v_\lambda:=\mathcal C(\lambda)^{-1}u_\lambda.
\]
At $\lambda=\lambda^\ast$ we recover $v_{\lambda^\ast}=w_1$. For all nearby $\lambda$, the profile solves the bulk, bed and interface equations \eqref{eq_mode_ode}--\eqref{eq_mode_jump}. By \eqref{eq_d}, undoing the factor $s/2$ used in the surface equation of \eqref{eq_kernelprob} gives the residual identities
\begin{equation}\label{eq_trans_family}
\widetilde L(\lambda)u_\lambda=\frac{2D_1(\lambda)}{s}\,\mathbf e,\qquad
L(\lambda)v_\lambda=\frac{2D_1(\lambda)}{s}\,\widehat{\mathcal C}^{-1}\mathbf e,\qquad
\mathbf e:=(0,\ 0,\ \cos x,\ 0).
\end{equation}

The fixed vector in \eqref{eq_trans_family} is outside the range:
\begin{equation}\label{eq_trans_offrange}
\mathbf e\notin\operatorname{Ran}\widetilde L(\lambda^\ast).
\end{equation}
Suppose that $\widetilde L(\lambda^\ast)\boldsymbol\chi=\mathbf e$ and project onto $\cos x$ with $P_1$. The operator commutes with $P_1$ because its coefficients are independent of $x$. The right-hand side is already in mode one, and the projection eliminates the scalars ${\widehat m}_1,{\widehat Q}$, giving
$P_1\boldsymbol\chi=\bigl(\chi^+(y)\cos x,\ \chi^-(y)\cos x,\ 0,\ 0\bigr),$
where $\chi^\pm$ are the restrictions of a continuous profile $\chi$. The homogeneous bulk, bed and interface equations determine $\chi=c\varphi_{\lambda^\ast}$ by \eqref{eq_f}. By \eqref{eq_trans_family}, its surface residual is $2cD_1(\lambda^\ast)\cos x/s=0$, contradicting the prescribed $\cos x$.

Differentiating \eqref{eq_trans_family} at $\lambda^\ast$, where $D_1(\lambda^\ast)=0$, we obtain
\[
\partial_\lambda\partial_vF(\lambda^\ast,0)w_1
=L'(\lambda^\ast)w_1
=\frac{2\,\partial_\lambda D_1(\lambda^\ast)}{s}\,\widehat{\mathcal C}^{-1}\mathbf e
-L_\ast\,\partial_\lambda v_\lambda\big|_{\lambda=\lambda^\ast}.
\]
The last term lies in $\operatorname{Ran} L_\ast=\widehat{\mathcal C}^{-1}\operatorname{Ran}\widetilde L(\lambda^\ast)$, and $\widehat{\mathcal C}^{-1}\mathbf e\notin\operatorname{Ran} L_\ast$ by \eqref{eq_trans_offrange}. \eqref{eq_hsimple} makes the coefficient of this vector nonzero and gives
$\partial_\lambda\partial_vF(\lambda^\ast,0)w_1\notin\operatorname{Ran} L_\ast.$
\end{proof}

The first interface cosine mode fixes the kernel coefficient on the local graph branch.

\begin{corollary}\label{cor_local}
Let $(\lambda^\ast,0)\in{\mathcal O}$ satisfy \eqref{eq_h} and \eqref{eq_hsimple}, and let $w_1$ be the kernel vector normalized as in \autoref{thm_trans}. We write $\zeta_1$ for the first cosine coefficient of the interface displacement in a graph state or perturbation and take $Z=\ker\zeta_1\subset X$.
Then $F=0$ has a real-analytic branch through $(\lambda^\ast,0)$ whose recovered graph states are
\[
\begin{aligned}
a\mapsto(\lambda(a),w(a))
&=\bigl(\lambda(a),w_0(\lambda(a))+a w_1+a z(a)\bigr),
&&|a|<a_0,\\
z(a)&\in Z,\qquad z(0)=0.
\end{aligned}
\]
Here $a$ is the coefficient of $w_1$ in this decomposition. The interface amplitude along the branch is
\begin{equation}\label{eq_local_interface_amplitude}
\tau:=\zeta_1(w(a))=-\frac{S_2}{{{W_{\!I}}}(\lambda^\ast)}\,a
\qquad (|a|<a_0).
\end{equation}
This curve and the laminar family exhaust all graph solutions near
$(\lambda^\ast,w_0(\lambda^\ast))$.

Moreover, the half-period translation $\mathcal T_\pi$ of
\autoref{rem_components} gives
\[
\lambda(-a)=\lambda(a),\qquad
\mathcal T_\pi w(a)=w(-a).
\]
\end{corollary}
\begin{proof}
The functional $\zeta_1$ is bounded on $X$ and by \eqref{eq_amps} takes the nonzero value $-S_2/{{W_{\!I}}}(\lambda^\ast)$ on $w_1$. Hence
$X=\operatorname{span}\{w_1\}\oplus Z.$
The admissible laminar curve of \autoref{prop_operator}(i) and the kernel and transversality conclusions of \autoref{thm_trans} verify the hypotheses of \autoref{thm_cr} for $F$ on a neighbourhood in ${\mathcal O}$. The theorem gives the local branch $v=a w_1+a z(a)$ and its uniqueness. Since $z(a)\in Z$ and the laminar interface is flat, only $a w_1$ contributes to $\zeta_1(w(a))$, which gives \eqref{eq_local_interface_amplitude}.

By \autoref{rem_components}, the translation $\mathcal T_\pi$ maps solutions to solutions, fixes the laminar state, and reverses the first cosine coefficient. Thus $\mathcal T_\pi w_1=-w_1$ and $\mathcal T_\pi Z=Z$, giving the translated perturbation
\[
\mathcal T_\pi\bigl(a w_1+a z(a)\bigr)
=(-a)w_1+(-a)\bigl(-\mathcal T_\pi z(a)\bigr).
\]
Local uniqueness at the parameter $-a$ gives
\[
\lambda(-a)=\lambda(a),\qquad
z(-a)=-\mathcal T_\pi z(a),\qquad
\mathcal T_\pi w(a)=w(-a)
\]
for small nonzero $a$. Analyticity extends these identities throughout $|a|<a_0$.
\end{proof}

\begin{remark}
For related local branches, see \cite[Theorems~3.4--3.6]{MM16} for short waves with a different normalization, \cite[Theorem~3.8]{Wh19} for distinct densities and fixed mean depths, and \cite[Corollary~4.3]{DMSW} for a rigid lid.
The simple-root transversality argument parallels \cite[Lemma~3.3]{MM16}, \cite[Lemma~3.5]{Wh19} and, for a rigid lid, \cite[Lemma~4.2]{DMSW}. Here \eqref{eq_trans_family} maps the mode-one profiles to scalar multiples of the fixed vector \(\widehat{\mathcal C}^{-1}\mathbf e\), which lies outside \(\operatorname{Ran}L_\ast\). The scalar factor \(2D_1(\lambda)/s\) gives transversality when \(\partial_\lambda D_1(\lambda^\ast)\ne0\) in the parameter \eqref{eq_lambda_definition}.
\end{remark}

\section{Folds of the dispersion relation and small loops}\label{sec_folds}

We now consider the two roots $\lambda_-(\mu)<\lambda_+(\mu)$ near a double root of the mode-one dispersion relation $D_1(\lambda;\mu)=0$. To decide whether the branches through these roots form a small loop, we determine the sign of the coefficient ${{\Theta}}$ defined in \eqref{eq_theta}.

Throughout this section $\mu$ denotes one of the physical parameters ${{\Omega_1}},{{\Omega_2}},h,d,g$, with the others fixed. The examples vary $\mu={{\Omega_2}}$. We write $F(\lambda,\mu,v)$ for the graph residual with components \eqref{eq_f1}--\eqref{eq_f3} at parameter $\mu$, evaluated on the complete state $w=w_0(\lambda,\mu)+v$. The residual and the laminar quantities $w_0(\lambda,\mu)$ and ${{W_{\!I}}}(\lambda,\mu)$ are analytic in $\mu$ on the admissible set. The spaces $X,Y$ are independent of $\mu$. The mode-one dispersion function of the laminar family is
$\mathcal D(\lambda,\mu):=D_1(\lambda;\mu).$

\begin{definition}\label{def_fold}
A point $(\lambda^\ast,\mu_0)$ with ${{W_{\!I}}}^\ast={{W_{\!I}}}(\lambda^\ast,\mu_0)$ is a \emph{nondegenerate fold} if $(\lambda^\ast,0)\in\mathcal O$ at $\mu=\mu_0$, $D_n(\lambda^\ast;\mu_0)\ne0$ for $n\ge2$, and
\begin{equation}\label{eq_fold}
\mathcal D=\partial_\lambda\mathcal D=0,\qquad \partial_\lambda^2\mathcal D\ne0,\qquad \partial_\mu\mathcal D\ne0\qquad\text{at }(\lambda^\ast,\mu_0).
\end{equation}
Here $(\lambda^\ast,0)\in\mathcal O$ gives ${{W_{\!I}}}^\ast\ne0\ne s^\ast$.
\end{definition}

Since $\partial_\lambda{{W_{\!I}}}=-1/(2\sqrt\lambda)\ne0$, the interface velocity gives a local $W$ coordinate. The chain rule gives
\begin{equation*}
\begin{aligned}
\partial_\lambda\mathcal D&=\partial_\lambda{{W_{\!I}}}\,\partial_WD_1,\\
\partial_\lambda^2\mathcal D&=(\partial_\lambda{{W_{\!I}}})^2\partial_W^2D_1
\qquad\text{at the double root}.
\end{aligned}
\end{equation*}
Thus \eqref{eq_fold} at $(\lambda^\ast,\mu_0)$ becomes
\[
D_1=\partial_WD_1=0,\qquad \partial_W^2D_1\ne0,\qquad \partial_\mu\mathcal D\ne0.
\]
Here $\partial_\mu\mathcal D$ is taken at fixed $\lambda$.
The \emph{two-root side} of $\mu_0$ is specified by
$\partial_\mu\mathcal D\,(\mu-\mu_0)/\partial^2_\lambda\mathcal D<0.$
We define the leading-order half-separation of the roots by
\begin{equation}\label{eq_eps}
\varepsilon=\varepsilon(\mu):=\Big|\frac{2\,\partial_\mu\mathcal D\,(\mu-\mu_0)}{\partial^2_\lambda\mathcal D}\Big|^{1/2}.
\end{equation}
On the two-root side, rescale $\lambda-\lambda^\ast=\varepsilon\xi_1$. The implicit function theorem gives two simple zeros $\lambda_-(\mu)<\lambda_+(\mu)$ of $\mathcal D(\cdot,\mu)$ near $\lambda^\ast$, with
$\lambda_\pm=\lambda^\ast\pm\varepsilon+O(\varepsilon^2).$
Since $\partial^2_\lambda\mathcal D$ stays nonzero nearby, Rolle's theorem rules out further zeros there. On the other side of $\mu_0$, there are no zeros near $\lambda^\ast$.

We next exclude resonance with a higher mode at the two nearby roots. Choose a compact neighbourhood of the fold where the parameters remain bounded and $h,d,h-d,|s|$ and $|{{W_{\!I}}}|$ are bounded away from zero. Set $\kappa=({{\Omega_1}}-{{\Omega_2}})/{{W_{\!I}}}$. In \eqref{eq_dsplit} we take $k=n$ and bound $\sinh(n(h-d))\cosh(nd)$ and $\sinh(n(h-d))\sinh(nd)$ by $\cosh(nh)$. We obtain
\begin{equation*}
\begin{aligned}
\left|\frac{D_n(\lambda;\mu)}{n\cosh(nh)}-s^2\right|
&\le \frac{|g+{{\Omega_1}}s|+|\kappa|s^2}{n}
     +\frac{|\kappa|\,|g+{{\Omega_1}}s|}{n^2}
     =O(n^{-1}).
\end{aligned}
\end{equation*}
The estimate is uniform on this neighbourhood, so $D_n\to+\infty$ uniformly as $n\to\infty$. The finitely many remaining modes $n\ge2$ stay nonzero near the fold by continuity because $D_n(\lambda^\ast;\mu_0)\ne0$. The nearby roots $\lambda_\pm(\mu)$ are simple and admissible, so they satisfy \eqref{eq_h} and \eqref{eq_hsimple}. If the fold lies in the window of \autoref{prop_pairing}, both roots are positive.

\begin{proposition}\label{rem_foldwindow}
Every double root $\lambda^\ast$ of $D_1$ in the component $\{s<0<{{W_{\!I}}}\}$ whose kernel is positive is a nondegenerate fold for $\mu={{\Omega_2}}$.
\end{proposition}
\begin{proof}
The signs $s<0<{{W_{\!I}}}$ give laminar admissibility, and \autoref{lem_sl} gives
$D_n(\lambda^\ast)>0\qquad(n\ge2).$
From \autoref{lem_convex}(i) we obtain $D_1=-c\,\mathcal E$, with $c>0$ analytic near $\lambda^\ast$ and independent of ${{\Omega_2}}$. At the double root, using the strict convexity of \autoref{lem_convex}(iii), we obtain
\[
\begin{aligned}
\partial_\lambda^2D_1&=-c\,\partial_\lambda^2\mathcal E\ne0,\qquad
\partial_{{{\Omega_2}}}D_1&=-c\,\partial_{{{\Omega_2}}}\mathcal E=-c\ne0.
\end{aligned}
\]
\end{proof}

\subsection{The reduced equation}\label{sub_11_the_reduced_equation}

The amplitude used below is the first interface cosine coefficient $\tau=\zeta_1(v)$. Translation by half a period changes its sign. We use this symmetry to obtain the scalar equation.

Along the laminar line define
$L(\lambda,\mu):=\partial_vF(\lambda,\mu,0).$
At the fold $L_\ast=L(\lambda^\ast,\mu_0)$ has a kernel spanned by a mode-one element $w_1$ by \autoref{prop_modes}. With the physical profile \eqref{eq_f}, the kernel vector has interface coefficient $-S_2/{{W_{\!I}}}^\ast\ne0$ by \eqref{eq_amps}. For the fold reduction we rescale the vector $w_1$ so that $\zeta_1(w_1)=1$, using the interface coefficient defined in \autoref{cor_local}.

\begin{lemma}\label{lem_ls}
Let $(\lambda^\ast,\mu_0)$ be a nondegenerate fold, with the mode-one kernel vector $w_1$ normalized by $\zeta_1(w_1)=1$. There are open neighbourhoods of $(\lambda^\ast,\mu_0,0)$, denoted by $\mathcal N$ in ${\mathbb{R}}\times{\mathbb{R}}\times X$ and by $\mathcal N'$ and $\mathcal N''$ in ${\mathbb{R}}^3$. The neighbourhood $\mathcal N''$ contains $\{(\lambda,\mu,\tau^2):(\lambda,\mu,\tau)\in\mathcal N'\}$. With $Z=\{v\in X:\zeta_1(v)=0\}$, there are real-analytic maps $\varrho:\mathcal N'\to Z$ and $G:\mathcal N''\to{\mathbb{R}}$ with the following properties.
\begin{enumerate}[label=\textup{(\roman*)},leftmargin=2.3em]
\item For $(\lambda,\mu,v)\in\mathcal N$, set $\tau=\zeta_1(v)$. Then $F(\lambda,\mu,v)=0$ if and only if $(\lambda,\mu,\tau)\in\mathcal N'$ and
$v=\tau w_1+\varrho(\lambda,\mu,\tau),\qquad \tau\,G(\lambda,\mu,\tau^2)=0.$
\item The range correction $\varrho$ satisfies
\begin{equation*}
\begin{aligned}
\varrho(\lambda,\mu,0)&=0,\qquad
\varrho(\lambda,\mu,-\tau)=\mathcal T_\pi\varrho(\lambda,\mu,\tau);
\end{aligned}
\end{equation*}
In particular, the perturbation $v$ in (i) vanishes exactly when $\tau=0$.
\item At zero amplitude the scalar factor $G$ satisfies
$G(\lambda,\mu,0)=\varpi(\lambda,\mu)\,\mathcal D(\lambda,\mu)$
 with $\varpi$ real-analytic and $\varpi(\lambda^\ast,\mu_0)\ne0$.
\end{enumerate}
\end{lemma}
\begin{proof}
\smallskip\noindent\textit{Step 1. Equivariant splittings of the kernel and range.}\enspace
By the normalization of the amplitude functional we have the direct sum
\[
\ker L_\ast=\operatorname{span}\{w_1\},\qquad X=\operatorname{span}\{w_1\}\oplus Z.
\]
The Fredholm index is zero by \autoref{thm_fred}. Put
$y_1:=(0,\ 0,\ \cos x,\ 0)\in Y.$
For the vector \(\mathbf e\) of \eqref{eq_trans_family}, the argument for \eqref{eq_trans_offrange} uses only \(\mathcal D(\lambda^\ast,\mu_0)=0\). With \(L=\widehat{\mathcal C}^{-1}\widetilde L\,\mathcal C\) from \autoref{lem_conj}, this gives
\begin{equation*}
\begin{aligned}
\mathbf e&\notin\operatorname{Ran}\widetilde L(\lambda^\ast,\mu_0),\qquad
\widehat{\mathcal C}^{-1}\mathbf e=(T_0^+)^2y_1=d^2y_1
\notin\operatorname{Ran} L_\ast.
\end{aligned}
\end{equation*}
Thus
$Y=\operatorname{Ran} L_\ast\oplus\operatorname{span}\{y_1\}.$
Let \(\mathsf Q\) be the projection onto \(\operatorname{Ran} L_\ast\) along \(y_1\), and let \(\ell\in Y^\ast\) satisfy \(\ell(y_1)=1\), \(\ell(\operatorname{Ran} L_\ast)=0\). The half-period translation commutes with \(L_\ast\) by equivariance and acts on mode \(n\) by \((-1)^n\). Hence
\begin{equation*}
\begin{aligned}
\mathcal T_\pi y_1&=-y_1, &\qquad \mathcal T_\pi w_1&=-w_1,\\
\mathcal T_\pi\operatorname{Ran} L_\ast&=\operatorname{Ran} L_\ast, &\qquad \mathcal T_\pi Z&=Z.
\end{aligned}
\end{equation*}
The invariance of \(\operatorname{span}\{y_1\}\) and \(\operatorname{Ran} L_\ast\) gives
\begin{equation*}
\begin{aligned}
\mathsf Q\mathcal T_\pi&=\mathcal T_\pi\mathsf Q,\qquad
\ell\circ\mathcal T_\pi&=-\ell.
\end{aligned}
\end{equation*}

\smallskip\noindent\textit{Step 2. The analytic range correction.}\enspace
Since $\mathsf QL_\ast|_Z:Z\to\operatorname{Ran} L_\ast$ is an isomorphism, the analytic implicit function theorem gives the unique correction $\varrho\in Z$ solving
$\mathsf QF(\lambda,\mu,\tau w_1+\varrho(\lambda,\mu,\tau))=0.$
The laminar line and equivariance under $\mathcal T_\pi$ (\autoref{rem_components}) give the identities in (ii) by uniqueness. See \cite[Chapter~VII, Propositions~1.5 and~3.3]{GS}.
The full equation reduces to $g_{\rm red}=0$, where
\[
g_{\rm red}(\lambda,\mu,\tau)
:=\ell\big(F(\lambda,\mu,\tau w_1+\varrho(\lambda,\mu,\tau))\big).
\]
Equivariance of $F$ and $\varrho$ and the identity $\ell\circ\mathcal T_\pi=-\ell$ make $g_{\rm red}$ odd in $\tau$, so its analytic expansion factors as
\begin{equation}\label{eq_reduced_scalar_factorization}
g_{\rm red}(\lambda,\mu,\tau)=\tau\,G(\lambda,\mu,\tau^2).
\end{equation}
Indeed the convergent joint Taylor series contains only odd powers of $\tau$. Replacing each $\tau^{2j+1}$ by $\sigma^j$ gives a series converging locally for $|\sigma|<r^2$ whenever $r$ is smaller than the original amplitude radius, after shrinking the parameter neighbourhood. Thus $G$ is analytic on a full neighbourhood of $(\lambda^\ast,\mu_0,0)$ \cite[Chapter~VI, Corollary~2.2 and the analytic remark after Lemma~2.1]{GS}.

\smallskip\noindent\textit{Step 3. The dispersion factor at zero amplitude.}\enspace
Differentiating \(g_{\rm red}\) in \(\tau\) at zero amplitude, we obtain
\[
G(\lambda,\mu,0)=\partial_\tau g_{\rm red}(\lambda,\mu,0)=\ell\big(L(\lambda,\mu)(w_1+z')\big),
\]
where \(z'=\partial_\tau\varrho(\lambda,\mu,0)\in Z\). Since \(\mathsf QL(\lambda,\mu)|_Z\) is an isomorphism onto \(\operatorname{Ran} L_\ast\) near \((\lambda^\ast,\mu_0)\), the differentiated range equation makes \(w_1+z'\) the unique vector satisfying
$\zeta_1(w_1+z')=1,\qquad \mathsf QL(\lambda,\mu)(w_1+z')=0.$
Let \(v_{\lambda,\mu}=\mathcal C(\lambda,\mu)^{-1}u_{\lambda,\mu}\) be the mode-one family of \eqref{eq_trans_family} for the physical parameter \(\mu\). By \eqref{eq_amps} and \eqref{eq_trans_family}, using \(\widehat{\mathcal C}^{-1}\mathbf e=d^2y_1\),
\begin{equation*}
\begin{aligned}
\zeta_1(v_{\lambda,\mu})&=-\frac{\sinh(h-d)}{{{W_{\!I}}}},\qquad
L(\lambda,\mu)v_{\lambda,\mu}&=\frac{2d^2\mathcal D(\lambda,\mu)}{s}\,y_1.
\end{aligned}
\end{equation*}
The normalized vector \(v_{\lambda,\mu}/\zeta_1(v_{\lambda,\mu})\) satisfies the same conditions as \(w_1+z'\), so
\begin{equation*}
\begin{aligned}
G(\lambda,\mu,0)
&=\frac{\ell\big(L(\lambda,\mu)v_{\lambda,\mu}\big)}{\zeta_1(v_{\lambda,\mu})}
=\varpi(\lambda,\mu)\mathcal D(\lambda,\mu),
\qquad \varpi=-\frac{2d^2\,{{W_{\!I}}}}{s\,\sinh(h-d)}.
\end{aligned}
\end{equation*}
At the fold \({{W_{\!I}}}\ne0\ne s\), so \(\varpi(\lambda^\ast,\mu_0)\ne0\).
\end{proof}

\subsection{The fold coefficient}\label{sub_11_the_fold_coefficient}

To determine the initial bending at the two roots, we retain the terms of order $(\lambda-\lambda^\ast)^2$ and $\tau^2$ in $G(\lambda,\mu,\tau^2)$.

Let
\begin{equation}\label{eq_abc}
\mathsf a:=\frac12\partial_\lambda^2G,\qquad \mathsf b:=\partial_\mu G,\qquad \mathsf c:=\partial_\sigma G\qquad\text{at }(\lambda^\ast,\mu_0,0),
\end{equation}
where $\sigma$ is the third argument of $G$. With $\varpi_\ast=\varpi(\lambda^\ast,\mu_0)$, \autoref{lem_ls}(iii) and \eqref{eq_fold} give
\[
G(\lambda^\ast,\mu_0,0)=\partial_\lambda G(\lambda^\ast,\mu_0,0)=0,\quad
\mathsf a=\frac12\varpi_\ast\partial_\lambda^2\mathcal D\ne0,\quad
\mathsf b=\varpi_\ast\partial_\mu\mathcal D\ne0.
\]
We define the \emph{fold coefficient}
\begin{equation}\label{eq_theta}
{{\Theta}}:=\frac{\mathsf c}{\mathsf a} .
\end{equation}

\begin{lemma}\label{lem_theta}
Let $(\lambda^\ast,\mu_0)$ be a nondegenerate fold with coefficient ${{\Theta}}$.
\begin{enumerate}[label=\textup{(\roman*)},leftmargin=2.3em]
\item For $\mu$ on the two-root side near $\mu_0$, let $\lambda_\pm=\lambda_\pm(\mu)$ be the two simple zeros of $\mathcal D(\cdot,\mu)$. On each Crandall--Rabinowitz branch through $(\lambda_\pm,0)$, with parameter $\tau=\zeta_1$,
$\lambda(\tau)=\lambda_\pm+\lambda_2^{(\pm)}\tau^2+O(\tau^4).$
\item On the two-root side, the \emph{bending coefficients} for the two simple roots in (i) converge to the fold coefficient,
\[
\chi_\pm(\mu):=\lambda_2^{(\pm)}(\lambda_\mp-\lambda_\pm)
\longrightarrow {{\Theta}}\qquad (\mu\to\mu_0).
\]
If ${{\Theta}}>0$, the local branch through $\lambda_\pm$ bends towards $\lambda_\mp$ near the fold.
\item Fix $g=k=1$, with the varying physical parameter $\mu$ chosen from the other parameters. The coefficient $m_{2,2}$ of $\tau^2$ in the lower flux in the third-order Stokes expansion at $\lambda_\pm$ (\autoref{lem_b2} and~\autoref{lem_b3}) gives
\[
\lambda_2^{(\pm)}=\frac{2\sqrt{\lambda_\pm}}{h-d}\,m_{2,2}.
\]
Let $\widetilde J_5^{(4)}$ consist of the first four columns of the mode-one Stokes matrix $\widetilde J_5$ in \eqref{eq_j5t}. For the third-order mode-one residual $\mathbf t^\flat$, evaluate \eqref{eq_order3} with the second-order correction to $m_1$ set to zero as in \autoref{lem_b2}(a).
We set $\widetilde\Delta=\det\widetilde J_5$ and $\mathfrak g=\det[\widetilde J_5^{(4)}|\mathbf t^\flat]$. Then
\begin{equation}\label{eq_thetaformula}
{{\Theta}}=\frac{2\,\mathfrak g}{\partial_\lambda^2\widetilde\Delta}=\frac{8\lambda^\ast\,\mathfrak g}{\partial_W^2\widetilde\Delta}\qquad\text{at the fold},
\end{equation}
with both determinants evaluated at $(\lambda^\ast,\mu_0)$ and $\partial_W$ denoting differentiation with respect to ${{W_{\!I}}}$ along the laminar family.
\end{enumerate}
\end{lemma}
\begin{proof}
\smallskip\noindent\textit{Step 1. The quadratic branch coefficient.}\enspace
Near $(\lambda_\pm,\mu,0)$,
\[
G(\lambda_\pm,\mu,0)=0,\qquad
\partial_\lambda G(\lambda_\pm,\mu,0)=\varpi\,\partial_\lambda\mathcal D\ne0.
\]
The implicit function theorem gives $\lambda=\Lambda(\mu,\tau^2)$ with $\Lambda$ analytic and $\Lambda(\mu,0)=\lambda_\pm$. By \autoref{lem_ls}(i), $\tau=\zeta_1(w)$. Differentiating $G(\Lambda(\mu,\sigma),\mu,\sigma)=0$ at $\sigma=0$, we obtain
\begin{equation}\label{eq_lambda2_derivative}
\lambda_2^{(\pm)}=-\frac{\partial_\sigma G}{\partial_\lambda G}(\lambda_\pm,\mu,0).
\end{equation}

\smallskip\noindent\textit{Step 2. The paired-root bending limit.}\enspace
For (ii), Taylor expansion of the derivatives in \eqref{eq_lambda2_derivative} at $(\lambda^\ast,\mu_0,0)$, with $\lambda_\pm-\lambda^\ast=\pm\varepsilon+O(\varepsilon^2)$ and $|\mu-\mu_0|=O(\varepsilon^2)$, gives
\[
\begin{aligned}
\partial_\lambda G(\lambda_\pm,\mu,0)&=\pm2\mathsf a\varepsilon+O(\varepsilon^2),\qquad
\partial_\sigma G(\lambda_\pm,\mu,0)&=\mathsf c+O(\varepsilon).
\end{aligned}
\]
Using \eqref{eq_lambda2_derivative} and the root separation gives
\begin{equation*}
\begin{aligned}
\lambda_2^{(\pm)}&=\mp\frac{\mathsf c}{2\mathsf a\varepsilon}+O(1),\qquad
\lambda_\mp-\lambda_\pm=\mp2\varepsilon(1+O(\varepsilon)),\\
\chi_\pm=\lambda_2^{(\pm)}(\lambda_\mp-\lambda_\pm)
&=\frac{\mathsf c}{\mathsf a}+O(\varepsilon).
\end{aligned}
\end{equation*}
These estimates also hold when $\mathsf c=0$.

\smallskip\noindent\textit{Step 3. The Stokes formula for the fold coefficient.}\enspace
For \eqref{eq_thetaformula}, use the square-root coordinate of \autoref{lem_convex}. Let $W_\pm={{W_{\!I}}}(\lambda_\pm)$, ${{W_{\!I}}}^\ast={{W_{\!I}}}(\lambda^\ast)$ and $\sigma_\pm=\sqrt{\lambda_\pm}$. By the paragraph after \autoref{def_fold}, \eqref{eq_h} and \eqref{eq_hsimple} hold at $\lambda_\pm$, and by \autoref{lem_b2}(b) and~\autoref{lem_b3},
\[
\lambda_2^{(\pm)}=\frac{2\sigma_\pm}{h-d}\Big[L_2^\flat-(h-d)C_2^\flat+(h-d)\frac{\mathfrak g}{\partial_W\widetilde\Delta}\Big]_{(W_\pm,\mu)} .
\]
The quantities $L_2^\flat$, $C_2^\flat$ and $\mathfrak g$ are given by explicit formulas (\autoref{lem_b2}) whose only denominators are $d$, $\Sigma_1$, $s$ and 
$\det\widetilde J_5(2)=4{{W_{\!I}}} D_2/(\sinh2d\,\sinh2(h-d))$;
 so $L_2^\flat$, $C_2^\flat$ and $\mathfrak g$ are continuous functions of $({{W_{\!I}}},\mu)$ near $({{W_{\!I}}}^\ast,\mu_0)$. For the determinant quotient, $\widetilde\Delta(\cdot;\mu)$ is analytic and has the simple zeros $W_+\ne W_-$. By \eqref{eq_det_jn},
\[
\widetilde\Delta=\frac{2{{W_{\!I}}}D_1}{\sinh d\sinh(h-d)}.
\]
Since ${{W_{\!I}}}^\ast\ne0$ and $D_1=\partial_WD_1=0\ne\partial^2_WD_1$ at the fold, $\partial_W^2\widetilde\Delta({{W_{\!I}}}^\ast;\mu_0)\ne0$. By the Hermite--Genocchi formula, 
$\widetilde\Delta(W)=(W-W_+)(W-W_-)\,\mathfrak d(W,\mu).$
The analytic factor is the second divided difference
\begin{equation*}
\begin{aligned}
\mathfrak d(W,\mu)&=\widetilde\Delta[W_+,W_-,W]\\
&=\int_0^1\!\!\int_0^{1-\theta_1}\partial_W^2\widetilde\Delta\big(W_++\theta_1(W_--W_+)+\theta_2(W-W_+);\mu\big)\,d\theta_2\,d\theta_1 ,
\end{aligned}
\end{equation*}
which is continuous and tends to $\frac12\partial^2_W\widetilde\Delta({{W_{\!I}}}^\ast;\mu_0)$ as $(W,\mu)\to({{W_{\!I}}}^\ast,\mu_0)$. Hence 
$\partial_W\widetilde\Delta(W_\pm)=(W_\pm-W_\mp)\mathfrak d(W_\pm,\mu).$
Since $\sigma=-{{W_{\!I}}}-{{\Omega_1}}d$ at fixed $\mu$, we convert the $W$-separation to the $\lambda$-separation by
\[
\lambda_\mp-\lambda_\pm=(\sigma_\mp-\sigma_\pm)(\sigma_\mp+\sigma_\pm)=(W_\pm-W_\mp)(\sigma_++\sigma_-).
\]
Substituting this separation into $\chi_\pm$, we obtain
\begin{equation*}
\begin{aligned}
\chi_\pm
&=\lambda_2^{(\pm)}(\lambda_\mp-\lambda_\pm)
=\frac{2\sigma_\pm}{h-d}
\big[L_2^\flat-(h-d)C_2^\flat\big]_{(W_\pm,\mu)}(\lambda_\mp-\lambda_\pm) \\
& \quad +\frac{2\sigma_\pm(\sigma_++\sigma_-)\,\mathfrak g(W_\pm,\mu)}{\mathfrak d(W_\pm,\mu)}
\longrightarrow\frac{8\lambda^\ast\,\mathfrak g}{\partial^2_W\widetilde\Delta}\Big|_{({{W_{\!I}}}^\ast,\mu_0)}
\end{aligned}
\end{equation*}
as $\mu\to\mu_0$ on the two-root side, since $L_2^\flat-(h-d)C_2^\flat$ stays bounded and $\lambda_\mp-\lambda_\pm\to0$. The limit $\chi_\pm\to{{\Theta}}$ from (ii) gives the formula with $\partial_W^2\widetilde\Delta$ in \eqref{eq_thetaformula}. The formula with $\partial_\lambda^2\widetilde\Delta$ follows from the chain rule at the fold, where $\partial_\lambda\widetilde\Delta=\partial_W\widetilde\Delta=0$ and 
\[
\partial_\lambda^2\widetilde\Delta=(\partial_\lambda{{W_{\!I}}})^2\partial_W^2\widetilde\Delta=\partial_W^2\widetilde\Delta/(4\lambda^\ast).
\]
\end{proof}

With the interface amplitude $\tau=\zeta_1(w)$ fixed, \autoref{lem_theta}(ii) expresses ${{\Theta}}$ as the limit of bending coefficients defined by the Crandall--Rabinowitz branches alone. Thus ${{\Theta}}$ is independent of the reduction choices in \autoref{lem_ls}. At the same fold, ${{\Theta}}$ is also independent of which physical parameter varies. The coefficients $\mathsf a$ and $\mathsf c$ are derivatives of $G(\cdot,\mu_0,\cdot)$, and at $\mu=\mu_0$ the reduction can be carried out with the same choices of $w_1$, $Z$ and $y_1$ for every parameter. Consequently \autoref{lem_theta}(iii) holds at every nondegenerate fold with $g=k=1$, in particular in the window, where $\mu={{\Omega_2}}$ is admissible (\autoref{rem_foldwindow}).

\subsection{Small loops}\label{sub_11_small_loops}

\begin{theorem}\label{thm_loops}
At a positive-kernel fold with $s<0<W_I$, sufficiently near the fold on the two-root side, the nontrivial graph solutions together with the two bifurcating laminar states form one simple closed analytic curve $\mathcal K_\mu$. Its maximal $|\tau|$ is $(\varepsilon/\sqrt\Theta)(1+O(\varepsilon^2))$. At the fold and on the other side all sufficiently nearby solutions are laminar.
\end{theorem}
\begin{proof}
We use the complement and the physical amplitude introduced in \autoref{lem_ls}. At the fold $G=G_\lambda=0$, $G_{\lambda\lambda}=2\mathsf a$, $G_\mu=\mathsf b$, and $G_\sigma=\mathsf c=\mathsf a\Theta\ne0$ by \autoref{thm_thetapos}. Thus $G=0$ is precisely $\tau^2=\Sigma(\lambda,\mu)$ near the fold. By implicit differentiation we obtain $\Sigma_{\lambda\lambda}=-2/\Theta<0$ and $\Sigma_\mu=-\mathsf b/\mathsf c$. The unique local maximum is attained at $\lambda_c=\lambda^\ast+O(|\mu-\mu_0|)$. Taylor expansion at $(\lambda^\ast,\mu_0)$ gives $\Sigma(\lambda_c,\mu)=-(\mathsf b/\mathsf c)(\mu-\mu_0)+O(|\mu-\mu_0|^2)$ since $\Sigma_\lambda=0$ at the fold. On the two-root side this is $\varepsilon^2/\Theta+O(\varepsilon^4)$ by \eqref{eq_eps}; taking the positive square root gives the asserted relative $O(\varepsilon^2)$ error since $\Theta>0$ is fixed. On the opposite side the maximum is strictly negative. At the fold the maximum is zero.
On the two-root side, $\Sigma$ vanishes precisely at the simple dispersion roots $\lambda_\pm$. It follows from strict concavity that $\Sigma$ is positive between these roots. The two graphs $\tau=\pm\sqrt\Sigma$ meet smoothly at their endpoints since $\Sigma_\lambda(\lambda_\pm,\mu)\ne0$. They are an analytic circle; reconstruction in \autoref{lem_ls} is injective with inverse given by $\lambda$ and the interface moment. Taylor expansion about $\lambda_c$ yields
$(\lambda-\lambda^\ast)^2+\Theta\tau^2=\varepsilon^2+O(\varepsilon^3)$
uniformly on the circle. Half-period translation interchanges the two signs of $\tau$. Since the reduction accounts for every nontrivial solution in this neighbourhood, the circle is isolated among the bifurcating waves.
\end{proof}
\begin{lemma}\label{lem_catseye}
Let $w_0=w_0(\lambda^\ast,\mu_0)$ be a laminar flow with ${{\Omega_1}}<0<{{\Omega_2}}$. Assume the laminar profile $W$ has critical heights as in \autoref{lem_cl}, with
\[
W(y_c^+)=W(y_c^-)=0,\qquad
y_c^+\in(-d,0),\quad y_c^-\in(-h,-d).
\]
There is a neighbourhood $\mathcal U'$ of $(\lambda^\ast,\mu_0,w_0)$ in ${\mathbb{R}}\times{\mathbb{R}}\times X$ such that each positively oriented solution $(\lambda,w)$ at parameter $\mu$ with $(\lambda,\mu,w)\in\mathcal U'$ has exactly four nondegenerate stagnation points per period. The four stagnation points have the following locations and types.
\[
\begin{array}{@{}lcc@{}}
 & x=0 & x=\pi\\
\text{upper layer, near }y_c^+ & \text{centre} & \text{saddle}\\
\text{lower layer, near }y_c^- & \text{saddle} & \text{centre}.
\end{array}
\]
In each layer, closed streamlines around the centre fill the region bounded by the saddle level set, forming a Kelvin cat's eye.
\end{lemma}
\begin{proof}
As \((\lambda,\mu,w)\to(\lambda^\ast,\mu_0,w_0)\) in \({\mathbb{R}}\times{\mathbb{R}}\times X\), the graph pullbacks satisfy
\[
\psi^\pm\circ\Pi_{Y^\pm}\longrightarrow\Psi\circ Y_0^\pm
\qquad\text{in }C^{3,\alpha}([0,2\pi]\times[0,1]).
\]
Since \(Y^\pm_t=T^\pm\) converges to the positive laminar thickness, the first two physical vertical derivatives converge uniformly after pullback. Choose disjoint neighbourhoods $N^\pm$ of $y_c^\pm$ inside the respective layers. Outside $N^\pm$ in the respective layer, $|\psi_y^\pm|$ is bounded away from zero, so there are no stagnation points. In $N^\pm$,
\[
\psi_{yy}^+={{\Omega_1}}+o(1)<0,\qquad
\psi_{yy}^-={{\Omega_2}}+o(1)>0.
\]
In each neighbourhood, $\psi_y^\pm$ has opposite signs at the two ends of each vertical segment by convergence to the laminar profile. Strict monotonicity gives exactly one zero $y=Y_c^\pm(x)$ on each segment, and uniqueness makes the graph periodic.

For $m^\pm(x):=\psi^\pm(x,Y_c^\pm(x))$, (N1) and evenness give
$(m^\pm)'(x)=-V^\pm(x,Y_c^\pm(x))<0\quad(0<x<\pi),$
with the opposite sign on $(\pi,2\pi)$. The only stagnation points are $(0,Y_c^\pm(0))$ and $(\pi,Y_c^\pm(\pi))$. At the stagnation points $\psi_{xy}^\pm=0$ and $\psi_{xx}^\pm=-V_x^\pm$, so (N2) and the signs of $\psi_{yy}^\pm$ give a non-degenerate upper centre at $0$ and saddle at $\pi$, and a non-degenerate lower saddle at $0$ and centre at $\pi$.

With vertical convexity or concavity and strict monotonicity of $m^\pm$ verified, the level-set argument in \cite[Theorem~4.1 and Lemma~4.3]{Wahlen09} gives closed curves around each centre, bounded by the level set of $\psi^\pm$ through the saddle.
\end{proof}

\begin{corollary}\label{cor_loops}
Let ${{\Omega_1}}<0$ and let $(\lambda^\ast,\mu_0)$ be a nondegenerate fold in the window of \autoref{prop_pairing}, with ${{\Theta}}>0$ and laminar critical heights $W(y_c^\pm)=0$ in both layers. For $\mu$ close to $\mu_0$ on the two-root side, every nontrivial solution on the closed curve $\mathcal K_\mu$ of \autoref{thm_loops} has exactly one Kelvin cat's eye in each layer. The solutions with $\tau=\zeta_1(w)<0$ are positively oriented.
\end{corollary}
\begin{proof}
Let \({\widehat V}_1\) be the vertical velocity of the kernel vector \(w_1\) normalized by \(\zeta_1(w_1)=1\). The graph velocity map \((\lambda,\mu,w)\mapsto(\widetilde V^+,\widetilde V^-)\) of \eqref{eq_vflatg} is analytic with values in \(C^2\) and vanishes at laminar flows. Condition (N) is equivalent to the fixed-strip inequalities \eqref{eq_nflat}, which are open by \autoref{lem_open}. On \(\mathcal K_\mu\),
$w-w_0(\lambda,\mu)=\tau\big(w_1+\widetilde\varrho(\lambda,\mu,\tau)\big).$
Here \(\widetilde\varrho=\varrho/\tau\) extends analytically to \(\tau=0\) and vanishes at \((\lambda^\ast,\mu_0,0)\) because \(\partial_\tau\varrho(\lambda^\ast,\mu_0,0)=0\). As \(\mathcal K_\mu\) shrinks to \((\lambda^\ast,w_0(\lambda^\ast,\mu_0))\), this gives
$\widetilde V/\tau\longrightarrow{\widehat V}_1\circ\Pi_{Y_0}\quad\text{in }C^2$
uniformly on the nontrivial part of \(\mathcal K_\mu\) as \(\mu\to\mu_0\). By \eqref{eq_amps} and $\zeta_1(w_1)=1$, the kernel velocity is
\[
{\widehat V}_1=-\frac{{{W_{\!I}}}^\ast}{\sinh(h-d)}\varphi_1(y)\sin x.
\]
In the fold window, $\varphi_1>0$ on $(-h,0]$ and ${{W_{\!I}}}^\ast>0$. Thus $-{\widehat V}_1$ is a positive multiple of $\varphi_1(y)\sin x$, and $-{\widehat V}_1\circ\Pi_{Y_0}$ satisfies \eqref{eq_nflat} by the kernel calculation in the proof of \autoref{prop_local_ng}.

From \autoref{lem_open} we obtain positive orientation for $\tau<0$, so \autoref{lem_catseye} gives one cat's eye in each layer for $\tau<0$. By \autoref{thm_loops}, the half-period translation $\mathcal T_\pi$ preserves the loop. Since $\mathcal T_\pi$ reverses $\tau$, the solutions with $\tau>0$ also have one cat's eye in each layer.
\end{proof}

\begin{proof}[Proof of \autoref{thm_b}]
Combine \autoref{thm_loops} with \autoref{cor_loops}. The identification with the graph continuation follows from isolation: analytic continuation stays on the circle at its regular points, and at either laminar point the unique bifurcating arc joins its two halves. The circle is compact in the admissible set, so the global alternative of \autoref{thm_bt} gives periodicity.
\end{proof}
\section{The nonlinear coefficient and its sign}\label{sec_bending_new}

We now calculate the coefficient which determines the shape of the small branch. We expand \eqref{eq_intro_model} to third order in the first cosine coefficient $\tau=\zeta_1(w)$ of the interface. In the Stokes expansion this coefficient is denoted by $z_1$. The calculation gives both the quadratic coefficient $\lambda_2$ in \eqref{eq_pitchfork} and the determinant formula \eqref{eq_thetaformula} for $\Theta$.

Take $g=k=1$ and fix a root $\lambda$ of $D_1(\lambda)=0$. We write $W={{W_{\!I}}}(\lambda)$ and $s=W+{{\Omega_1}}d$, with $W\ne0\ne s$. Let $m_2^0=m_2(\lambda)$ be the lower flux \eqref{eq_lampar} of the laminar flow at this root. The expansion uses $m_2$ as parameter and recovers $\lambda$ from $m_2=m_2(\lambda)$ (\autoref{lem_b3}). The rational-function identities hold for either sign of $s$ and $W$. We write
$G_n=\coth nd,\qquad \Sigma_n=\operatorname{csch}nd,\qquad K_n=\coth n(h-d),$
so that $G_n^2-\Sigma_n^2=1$ and $K=K_1=\coth(h-d)$.

\subsection{Boundary expansion in the interface amplitude}\label{sub_15_ansatz_and_relations}

Expanding the two surface and three interface conditions by Fourier mode, we obtain the systems \eqref{eq_j5t}--\eqref{eq_order3} through cubic order.

For the harmonic modes in the two layers we use the following basis functions
\begin{equation*}
\begin{aligned}
\mathsf S_n(y)&=\sinh n(y+d)/\sinh nd=\cosh ny+G_n\sinh ny\\
\mathsf T_n(y)&=-\sinh ny/\sinh nd=-\Sigma_n\sinh ny\\
\mathsf R_n(y)&=\sinh n(y+h)/\sinh n(h-d)=\cosh n(y+d)+K_n\sinh n(y+d).
\end{aligned}
\end{equation*}
Let $m_1^0=-Wd-\frac12{{\Omega_1}}d^2$ be the laminar value of $m_1$ at the root.
For the expansion through order $\tau^3$ in the first interface cosine coefficient $\tau=z_1$, the stream functions consist of quadratic polynomials in $y$ and harmonic corrections. The sums collect the modes retained at each order.
\begin{equation}\label{eq_ansatz}
\begin{aligned}
\psi^+&=\frac12{{\Omega_1}}y^2+Ay+B+\sum_{n\ge1}\big(a_n\mathsf S_n(y)+b_n\mathsf T_n(y)\big)\cos nx,\\
\psi^-&=m_1^0+L+C(y+d)+\frac12{{\Omega_2}}(y+d)^2+\sum_{n\ge1}r_n\mathsf R_n(y)\cos nx.
\end{aligned}
\end{equation}
The surface and interface are $y=\eta(x)$ and $y=b(x)$, with
\[
\eta(x)=\sum_{n\ge1}e_n\cos nx,\qquad
b(x)=-d+\zeta(x),\qquad
\zeta(x)=\sum_{n\ge1}z_n\cos nx.
\]
We set $m_1=m_1^0+M$. The coefficients $a_n$ and $b_n$ are the values of the upper-layer mode at the laminar reference levels $y=0$ and $y=-d$, respectively, since $\mathsf S_n(0)=1=\mathsf T_n(-d)$ and $\mathsf S_n(-d)=0=\mathsf T_n(0)$.

The ansatz solves the two interior Poisson equations term by term and its lower modes vanish at the bed. Its bed trace is
\[
\psi^-(x,-h)=m_1^0+L-C(h-d)+\frac12{{\Omega_2}}(h-d)^2.
\]
Hence \eqref{eq_intro_model} defines the lower flux by
\begin{equation}\label{eq_lambda_s}
m_2:=\psi^-(\cdot,-h)-m_1=m_2^0+L-M-(h-d)(C-W),\qquad m_2^0=-W(h-d)+\frac12{{\Omega_2}}(h-d)^2 .
\end{equation}
The means of $\eta$ and $\zeta$ vanish, so \eqref{eq_intro_mean_levels} holds. The remaining conditions are the \hypertarget{stokes_boundary_relations}{five boundary relations}
\[
\begin{aligned}
&\mathrm{S1}:\ \psi^+(x,\eta)=0,\qquad \mathrm{S2}:\ |\nabla\psi^+|^2(x,\eta)+2\eta-Q=0,\\
&\mathrm{I1}:\ \psi^+(x,b)-m_1=0,\qquad \mathrm{I2}:\ \psi^-(x,b)-m_1=0,\qquad \mathrm{I3}:\ \psi^+_y(x,b)-\psi^-_y(x,b)=0,
\end{aligned}
\]
where I3 is equivalent to $[\nabla\psi]=0$ given I1 and I2 (\eqref{eq_trace_reduction}). At $\tau=0$ the laminar flow at the root is $A=s$, $B=0$, $C=W$, $L=M=0$, $Q=s^2$, all modes zero; it satisfies S1--I3.

Half-period translation reverses $\tau$ and multiplies mode $n$ by $(-1)^n$. We take the first-order vector $(e_1,a_1,b_1,r_1,1)$ in $\ker\widetilde J_5$. Quadratic products of mode one force the mean and mode $\cos2x$. At third order the cubic terms and the products of first- and second-order terms force modes one and three. \autoref{lem_b3} identifies the quadratic coefficient of $\lambda$ from this expansion with $\lambda_2$ in the exact local branch \eqref{eq_pitchfork}.

We denote the coefficients in the expansion as follows.
\begin{equation*}
\begin{aligned}
&(A-s,B,C-W,L,M,Q-s^2)=O(\tau^2)\quad\text{(even)},\\
&\quad\text{with respective }\tau^2\text{-coefficients }(A_2,B_2,C_2,L_2,M_2,Q_2);\\
&e_1=\tau e_1^{(1)}+\tau^3e_1^{(3)}+O(\tau^5)\quad\text{(odd)},\quad\text{and likewise for }a_1,b_1,r_1;\\
&(e_2,z_2,a_2,b_2,r_2)=O(\tau^2)\quad\text{(even)},\\
&\quad\text{with respective }\tau^2\text{-coefficients }(e_2^{(2)},z_2^{(2)},a_2^{(2)},b_2^{(2)},r_2^{(2)}).
\end{aligned}
\end{equation*}
We drop the superscripts $(1)$ and $(2)$ and write $e_1,a_1,b_1,r_1$ for the $\tau$-coefficients of mode one and $e_2,z_2,a_2,b_2,r_2$ for the $\tau^2$-coefficients of mode two; the primes in $e_1',a_1',b_1',r_1'$ label the $\tau^3$-coefficients of mode one.

Substituting into S1--I3 and expanding $\cosh(n\eta)$, $\sinh(n\eta)$, $\cosh(n\zeta)$, $\sinh(n\zeta)$ in powers of $\tau$ makes each relation a power series in $\tau$ whose coefficients are trigonometric polynomials. We record the coefficients that are needed. We set $\mathsf u_1:=G_1a_1-\Sigma_1b_1+{{\Omega_1}}e_1$, the $\tau$-coefficient of the mode-one part of $\psi^+_y$ at the surface.

\emph{Order $\tau$, mode $1$.} For mode $\cos nx$ at any order, the linear part of the five boundary relations has matrix $\widetilde J_5(n)$ in the unknowns $(e_n,a_n,b_n,r_n,z_n)$.
\begin{equation}\label{eq_j5t}
\widetilde J_5(n)=\begin{pmatrix}
s&1&0&0&0\\
2(1+{{\Omega_1}}s)&2nG_ns&-2n\Sigma_ns&0&0\\
0&0&1&0&W\\
0&0&0&1&W\\
0&n\Sigma_n&-nG_n&-nK_n&{{\Omega_1}}-{{\Omega_2}}
\end{pmatrix},\qquad \widetilde J_5:=\widetilde J_5(1).
\end{equation}
Since $z_1=\tau$, the first-order interface coefficient is $1$. Taking $n=1$, we obtain
$\widetilde J_5\,(e_1,a_1,b_1,r_1,1)^{\!\top}=0.$

\emph{Order $\tau^2$, mode $0$.} Define
\[
\begin{gathered}
\mathsf g_1=\frac12e_1(G_1a_1-\Sigma_1b_1)+\frac14{{\Omega_1}}e_1^2,\qquad
\mathsf g_3=\frac12(\Sigma_1a_1-G_1b_1)+\frac14{{\Omega_1}},\\
\mathsf g_4=\frac12K_1r_1+\frac14{{\Omega_2}},\qquad
\mathsf g_5=\frac12(b_1-r_1).
\end{gathered}
\]
The coefficients of $\tau^2$ in the means of S1--I3 give the five relations
\begin{equation}\label{eq_order2mode0}
\begin{aligned}
B_2+\mathsf g_1&=0,\\
2sA_2-Q_2+\frac12\mathsf u_1^2+sa_1e_1+\frac12a_1^2&=0,\\
-dA_2+B_2-M_2+\mathsf g_3&=0,\\
L_2-M_2+\mathsf g_4&=0,\\
A_2-C_2+\mathsf g_5&=0.
\end{aligned}
\end{equation}

\emph{Order $\tau^2$, mode $2$.} The $\tau^2\cos 2x$ coefficients of S1--I3 give
\begin{equation}\label{eq_order2mode2}
\widetilde J_5(2)\,(e_2,a_2,b_2,r_2,z_2)^{\!\top}+\big(\mathsf g_1,\ \frac12\mathsf u_1^2+sa_1e_1-\frac12a_1^2,\ \mathsf g_3,\ \mathsf g_4,\ \mathsf g_5\big)^{\!\top}=0 .
\end{equation}

We write $F_j$ for the five entries of the forcing column in \eqref{eq_order2mode2}. Three rows give $a_2=-F_1-se_2$, $b_2=-F_3-Wz_2$, $r_2=-F_4-Wz_2$. Consequently the only inversion needed for the second harmonic is
\[
 \begin{pmatrix}2(1+\Omega_1s)-4G_2s^2&4\Sigma_2sW\\
 -2\Sigma_2s&\Omega_1-\Omega_2+2W(G_2+K_2)\end{pmatrix}
 \binom{e_2}{z_2}
 =\binom{-F_2+4s(G_2F_1-\Sigma_2 F_3)}
 {-F_5+2(\Sigma_2 F_1-G_2 F_3-K_2 F_4)}.
\]
The same elimination applies at the fundamental fold. It requires invertibility only in the second harmonic, which persists when transversality fails in the fundamental mode.

\emph{Order $\tau^3$, mode $1$.} Since $z_1=\tau$ has no cubic term, the linear part uses the first four columns $\widetilde J_5^{(4)}$ of $\widetilde J_5$. The equations are
\[
\widetilde J_5^{(4)}(e_1',a_1',b_1',r_1')^{\!\top}+\mathbf t=0,
\qquad \mathbf t=(\mathsf t_1,\dots,\mathsf t_5)^{\!\top},
\]
with residual components
\begin{equation}\label{eq_order3}
\begin{aligned}
\mathsf t_1&=A_2e_1+G_2a_2e_1-\Sigma_2b_2e_1+\frac38a_1e_1^2+\frac12e_2\mathsf u_1,\\
\mathsf t_2&=2A_2\mathsf u_1+\frac34se_1^2(G_1a_1-\Sigma_1b_1)+2a_1e_1(G_1a_1-\Sigma_1b_1)+\frac32{{\Omega_1}}a_1e_1^2\\
&\quad+a_2(2G_2\mathsf u_1+4se_1+2a_1)-2\Sigma_2\mathsf u_1b_2+e_2({{\Omega_1}}\mathsf u_1+sa_1),\\
\mathsf t_3&=A_2+\Sigma_2a_2-G_2b_2+\frac38b_1+\frac12z_2(\Sigma_1a_1-G_1b_1+{{\Omega_1}}),\\
\mathsf t_4&=C_2+K_2r_2+\frac38r_1+\frac12z_2(K_1r_1+{{\Omega_2}}),\\
\mathsf t_5&=\frac38(\Sigma_1a_1-G_1b_1-K_1r_1)+2b_2-2r_2+\frac12z_2(b_1-r_1).
\end{aligned}
\end{equation}
Only $A_2$ and $C_2$ enter from order $\tau^2$, mode $0$.

\emph{Order $\tau^3$, mode $3$.} Let $\mathbf x_3$ be the vector of $\tau^3$-coefficients of $(e_3,a_3,b_3,r_3,z_3)$ in that order. We set $\mathbf x_3=0$ in the left-hand sides of S1, S2, I1, I2 and I3. The $\tau^3\cos3x$ coefficients of the five residuals, in the same order, form the source vector $\mathbf y_3$, a polynomial in the mode-one coefficients of order $\tau$ and the mode-zero and mode-two coefficients of order $\tau^2$. The mode-three equations are
$\widetilde J_5(3)\,\mathbf x_3+\mathbf y_3=0.$

\subsection{Solution of the expansion}\label{sub_15_solution_of_the_expansion}

For $n\ge1$, let $D_n(W)$ denote the expression \eqref{eq_d} at $k=n$, $g=1$, with ${{W_{\!I}}}=W$ and $s=W+{{\Omega_1}}d$. On the branch $s<0$, this is the dispersion function $D_n(\lambda)$ at $\lambda=s^2$. The mirror symmetry \eqref{eq_mirror_symmetry} changes the signs of $s,W,{{\Omega_1}},{{\Omega_2}}$ together and leaves $s^2$, ${{\Omega_1}}s$ and $({{\Omega_1}}-{{\Omega_2}})/W$ fixed.

\begin{lemma}\label{lem_b1}
For $n\ge1$,
\begin{equation}\label{eq_det_jn}
\det\widetilde J_5(n)=\frac{2nW\,D_n(W)}{\sinh nd\,\sinh n(h-d)} .
\end{equation}
At the fixed mode-one root $W$ of $D_1$, $\ker\widetilde J_5$ is spanned by $\widetilde v_5=(e_1,a_1,b_1,r_1,1)$ with
\[
b_1=r_1=-W,\qquad a_1=\frac{-W(G_1+K_1)-({{\Omega_1}}-{{\Omega_2}})}{\Sigma_1},\qquad e_1=-\frac{a_1}s.
\]
If the root lies in the component $\{s<0<W\}$ with $q\ne0$, then $a_1=-W\varsigma(q+t)/q$, where $t=\tanh d$, $\varsigma=\operatorname{sech}d$ and $q=q(\lambda)=s^2/(1+{{\Omega_1}}s)-t$.
\end{lemma}

\begin{proof}
\smallskip\noindent\textit{Step 1. The mode determinant.}\enspace
Subtract $W$ times the third and fourth columns from the fifth; this does not change the determinant and turns the third and fourth rows into unit rows. Expanding along them leaves the determinant of the $3\times3$ matrix formed by rows $1,2,5$ and columns $1,2,5$:
\[
\begin{aligned}
\det\widetilde J_5(n)&=\det\begin{pmatrix}s&1&0\\2(1+{{\Omega_1}}s)&2nG_ns&2n\Sigma_nsW\\0&n\Sigma_n&{{\Omega_1}}-{{\Omega_2}}+nW(G_n+K_n)\end{pmatrix}\\
&=\big({{\Omega_1}}-{{\Omega_2}}+nW(G_n+K_n)\big)\big(2nG_ns^2-2(1+{{\Omega_1}}s)\big)-2n^2\Sigma_n^2s^2W .
\end{aligned}
\]
On the other hand, by \eqref{eq_d} with $k=n$, $\widehat B/\sinh n(h-d)=K_n+({{\Omega_1}}-{{\Omega_2}})/(nW)$, so
\begin{equation*}
\begin{aligned}
\frac{2nW\,D_n}{\sinh nd\,\sinh n(h-d)}&=2nW\Big[s^2n\Big(1+G_n\frac{\widehat B}{\sinh n(h-d)}\Big)\\
&\qquad-(1+{{\Omega_1}}s)\Big(G_n+\frac{\widehat B}{\sinh n(h-d)}\Big)\Big].
\end{aligned}
\end{equation*}
Substituting $\widehat B/\sinh n(h-d)$ and using $G_n^2-\Sigma_n^2=1$, we obtain the same factored expression as the column reduction of $\det\widetilde J_5(n)$, proving \eqref{eq_det_jn}.

\smallskip\noindent\textit{Step 2. The normalized kernel and its expression in the fold coordinates.}\enspace
The first four columns of $\widetilde J_5$ are linearly independent: if 
$\widetilde J_5^{(4)}(e,a,b,r)^{\!\top}=0$,
 rows $3$ and $4$ give $b=r=0$, row $5$ gives $\Sigma_1a=0$, so $a=0$, and row $1$ gives $se=0$, so $e=0$. At the root 
$\det\widetilde J_5=0$,
 so $\widetilde J_5$ has rank four, and every nonzero kernel vector has a nonzero last component; normalized by $z_1=1$, rows $3$, $4$, $5$, $1$ give the components $b_1$, $r_1$, $a_1$, $e_1$ of $\widetilde v_5$, and row $2$ holds because the kernel is nontrivial.

For a root in \(\{s<0<W\}\) with \(q\ne0\), \autoref{lem_convex}(i) gives \(D_1=0\) if and only if \(\mathcal E=0\). That argument uses \(q\ne0\) and \(1+{{\Omega_1}}s\ne0\), the latter implicit in the definition of \(q\). Thus
\[
{{\Omega_1}}-{{\Omega_2}}=W(\rho-K_1),\qquad
\Sigma_1a_1=-W(G_1+\rho).
\]
Using \(G_1/\Sigma_1=1/\varsigma\) and \(1/\Sigma_1=t/\varsigma\) yields
$a_1=-W(G_1+\rho)/\Sigma_1=-W(1+t\rho)/\varsigma.$
Substitute \(1+t\rho=(1-t^2)(q+t)/q\) and \(\varsigma^2=1-t^2\) to obtain
\begin{equation}\label{eq_fold_kernel_amplitude}
a_1=-W\varsigma(q+t)/q.
\end{equation}
\end{proof}

\begin{lemma}\label{lem_b2}
Suppose the fixed mode-one root $W$ with $D_1(W)=0$ satisfies $\partial_WD_1\ne0$, $D_2\ne0$ and $D_3\ne0$. Regard $\widetilde\Delta(W):=\det\widetilde J_5$ as a function of $W$ along the laminar family, with $s=W+{{\Omega_1}}d$ and all other parameters fixed.
\begin{enumerate}[label=\textup{(\alph*)},leftmargin=2.3em]
\item The relations S1--I3 can be solved order by order through $\tau^3$, with $z_1=\tau$. The order-$\tau^2$ mode-zero system \eqref{eq_order2mode0} has one free parameter $\omega$, normalized by $M_2=-d\omega$. The superscript $\flat$ denotes the solution of the mode-zero system with $M_2=0$, and $\mathbf t^\flat$ is \eqref{eq_order3} evaluated with $A_2^\flat,C_2^\flat$. The order-$\tau^3$ mode-one relations are solvable exactly for
\[
\omega=-\frac{\mathfrak g}{\partial_W\widetilde\Delta},\qquad
\mathfrak g:=\det\big[\widetilde J_5^{(4)}\,\big|\,\mathbf t^\flat\big].
\]
The mode-two correction at order $\tau^2$ is uniquely determined by the first-order coefficients. The full mode-zero family before imposing third-order solvability is
\[
\begin{aligned}
(A_2,B_2,C_2,L_2,M_2,Q_2)
&=(A_2^\flat,B_2^\flat,C_2^\flat,L_2^\flat,0,Q_2^\flat)\\
&\quad+\omega\,(1,0,1,-d,-d,2s),\qquad\omega\in{\mathbb{R}}.
\end{aligned}
\]
The reference coefficients at $M_2=0$ satisfy
\begin{equation*}
\begin{aligned}
&B_2^\flat=-\mathsf g_1,\quad
A_2^\flat=(\mathsf g_3-\mathsf g_1)/d,\quad
C_2^\flat=A_2^\flat+\mathsf g_5,\quad
L_2^\flat=-\mathsf g_4,\\
&Q_2^\flat=2sA_2^\flat+\frac12\mathsf u_1^2+sa_1e_1+\frac12a_1^2.
\end{aligned}
\end{equation*}
\item The solvability value of $\omega$ gives the $\tau^2$-coefficient of the lower flux $m_2$ in \eqref{eq_lambda_s},
\[
m_{2,2}=L_2^\flat-(h-d)C_2^\flat+(h-d)\frac{\mathfrak g}{\partial_W\widetilde\Delta} .
\]
\item The cofactor formula for $\mathfrak g$ is
\[
\begin{aligned}
\mathfrak g
&=\Sigma_1\big(2(1+{{\Omega_1}}s)\mathsf t_1^\flat-s\mathsf t_2^\flat\big)
+2\big(s^2-(1+{{\Omega_1}}s)G_1\big)\mathsf t_3^\flat\\
&\quad+2\big(s^2G_1-(1+{{\Omega_1}}s)\big)
       \big(K_1\mathsf t_4^\flat+\mathsf t_5^\flat\big).
\end{aligned}
\]

\end{enumerate}
\end{lemma}
\begin{proof}
\smallskip\noindent\textit{Step 1. The kernel and second-order modes.}\enspace
By \autoref{lem_b1}, the order-$\tau$ solution normalized by $z_1=1$
is unique and the first four columns of $\widetilde J_5$ are independent.

With the first-order vector fixed, the mode-zero corrections form a one-parameter family parametrized by $M_2=-d\omega$.
The system \eqref{eq_order2mode0} is triangular: the first relation gives $B_2$, the third $A_2$ in terms of $M_2$, the fourth $L_2$, the fifth $C_2$, the second $Q_2$. The vector $(1,0,1,-d,-d,2s)$ solves the homogeneous system.

The mode-two correction is unique and independent of $\omega$. Since $D_2\ne0$, \eqref{eq_det_jn} gives invertibility of $\widetilde J_5(2)$. The forcing in \eqref{eq_order2mode2} depends only on the first-order solution.

\smallskip\noindent\textit{Step 2. Third-order solvability.}\enspace
The order-$\tau^3$ mode-one equation fixes $\omega$ by projecting the forcing onto the cokernel of $\widetilde J_5$.
By \eqref{eq_order3}, $\mathbf t$ depends on the order-$\tau^2$ mode-zero coefficients only through $A_2$ and $C_2$, linearly, and 
$\mathbf t=\mathbf t^\flat+\omega(e_1,2\mathsf u_1,1,1,0)^{\!\top}.$
Differentiating the entries of $\widetilde J_5$ in $W$ and using $\partial_Ws=1$ and $\partial_W(1+{{\Omega_1}}s)={{\Omega_1}}$, we obtain
\begin{equation*}
\begin{aligned}
\partial_W\widetilde J_5\,\widetilde v_5
&=(e_1,\ 2{{\Omega_1}}e_1+2G_1a_1-2\Sigma_1b_1,\ 1,\ 1,\ 0)^{\!\top}\\
&=(e_1,2\mathsf u_1,1,1,0)^{\!\top}.
\end{aligned}
\end{equation*}
So, substituting this derivative vector into the order-three residual, the relations read
\[
\widetilde J_5^{(4)}(e_1',a_1',b_1',r_1')^{\!\top}+\mathbf t^\flat+\omega\,\partial_W\widetilde J_5\,\widetilde v_5=0.
\]
By independence of the first four columns, solvability is equivalent to
\[
\det[\widetilde J_5^{(4)}\,|\,\mathbf t^\flat+\omega\,\partial_W\widetilde J_5\widetilde v_5]=0.
\]
In every term of the determinant derivative with one of the first four columns differentiated, use $\widetilde J_5\widetilde v_5=0$ to express the fifth column as a combination of the first four. Since $(\widetilde v_5)_5=1$, multilinearity gives
\[
\partial_W\widetilde\Delta
=\det[\widetilde J_5^{(4)}\,|\,\partial_W\widetilde J_5\,\widetilde v_5].
\]
Hence the solvability condition is \(\mathfrak g+\omega\,\partial_W\widetilde\Delta=0\). Differentiating \eqref{eq_det_jn} at \(D_1(W)=0\), we obtain
\[
\partial_W\widetilde\Delta
=\frac{2W\,\partial_WD_1(W)}{\sinh d\,\sinh(h-d)}\ne0,
\]
since \(W\ne0\) and the root is simple. This determines \(\omega\) in (a). The same determinant formula makes \(\widetilde J_5(3)\) invertible since \(D_3\ne0\), so the third-order mode-three equation is solvable.

\smallskip\noindent\textit{Step 3. Lower flux and cofactor formulas.}\enspace
The mean correction in (a) gives the lower-flux coefficient in (b). Taking the $\tau^2$-coefficient in \eqref{eq_lambda_s} and using $C_2$ for the $\tau^2$-coefficient of $C-W$, we obtain
\begin{equation*}
\begin{aligned}
m_{2,2}&=L_2-M_2-(h-d)C_2\\
&=L_2^\flat-(h-d)C_2^\flat-(h-d)\omega.
\end{aligned}
\end{equation*}
The solvability value of $\omega$ from Step~2 gives (b).

Expanding $\det[\widetilde J_5^{(4)}|\mathbf t^\flat]$ along the last column and using $G_1^2-\Sigma_1^2=1$, we obtain the cofactors in (c).
\end{proof}

\begin{lemma}\label{lem_b3}
Let $(\lambda^\ast,0)\in{\mathcal O}$ satisfy \eqref{eq_h} and \eqref{eq_hsimple} with $g=k=1$. Parametrize the graph branch of \autoref{cor_local} by $\tau$ using \eqref{eq_local_interface_amplitude}, so that $\lambda(\tau)=\lambda^\ast+\lambda_2\tau^2+O(\tau^4)$. Then
\[
\lambda_2=\frac{2\sqrt{\lambda^\ast}}{h-d}\,m_{2,2},
\]
where $m_{2,2}$ is the coefficient of \autoref{lem_b2}(b) at $W={{W_{\!I}}}(\lambda^\ast)$.
\end{lemma}

\begin{proof}
We write $[\tau^j]$ for Taylor coefficients (including $1/j!$), and use $Y=Y_0+\sum_{j=1}^3\tau^jY_j+O(\tau^4)$ in each of the two affine maps of \eqref{eq_flat}. We write $\Phi_j=[\tau^j]\Phi$ on the fixed strip. Define physical coefficients $u_j$ recursively on the \emph{flat} layer, starting with the quadratic laminar field $u_0$. Suppress composition with $(x,Y_0(x,t))$ on the right below; primes denote physical $y$-derivatives. The chain rule gives
\[
\begin{aligned}
 \Phi_1 &=u_1+Y_1u_0',\\
 \Phi_2 &=u_2+Y_1u_1'+Y_2u_0'+\tfrac12Y_1^2u_0'',\\
 \Phi_3 &=u_3+Y_1u_2'+Y_2u_1'+\tfrac12Y_1^2u_1''
             +Y_3u_0'+Y_1Y_2u_0'' .
\end{aligned}
\]
Here $\partial_y^3u_0=0$. These are triangular definitions on a closed strip since $Y_0$ is affine with positive slope. On each compact substrip of the open physical layer the inverse map is defined for all small $\tau$. Differentiating there at fixed physical coordinates identifies the recursively defined $u_j$ with the Eulerian coefficients and, we obtain $\Delta u_j=0$ for $j\ge1$, since the sources are constant. Boundary identities follow by taking the one-sided traces of the fixed-strip identities. For example $\psi_y\circ\Pi_Y=\Phi_t/Y_t$ and $\psi_x\circ\Pi_Y=\Phi_x-Y_x\Phi_t/Y_t$. Substitution at the strip endpoints gives S1 and I1--I2 by composition, S2 by squaring the two velocity traces, and I3 by subtracting the vertical traces.

For clarity about regularity, $Y_j,\Phi_j$ are coefficients in $C^{3,\alpha}$ by analyticity in \eqref{eq_spaces}. The first harmonic at order one is smooth up to the flat boundaries. At order two its products force only the mean and the second harmonic; the mode equations and nonresonance leave only these smooth modes (the fundamental is excluded by parity). Thus in the third identity all differentiated lower-order coefficients are smooth. The order-three velocity traces require at most one derivative of $u_3$, two of $u_2$, three of $u_1$, and four of the quadratic $u_0$. The recursions and the endpoint identities therefore hold, in particular, in $C^{1,\alpha}$ on each closed side for the velocity traces. This argument uses coefficients and one-sided traces only.

Taking the mean and the required cosine coefficients now, we obtain the displayed boundary systems. Normalization fixes the first interface coefficient to $\tau$ exactly. At second order the mean equations and the nonsingular mode-two system fix all coefficients up to the change of flux. At third order the first-mode compatibility condition of \autoref{lem_b2} fixes $m_{2,2}$ uniquely since the dispersion zero is simple. Finally $m_2'(\lambda^\ast)=(h-d)/(2\sqrt{\lambda^\ast})$ in \eqref{eq_lampar}; inversion gives the stated coefficient.
\end{proof}
\subsection{The sign of the fold coefficient}\label{sec_thetasign}

Consider a double root $\lambda^\ast\in\mathcal J$ with positive kernel as in \autoref{prop_pairing}. On the fold manifold \eqref{eq_foldmanifold}, formula \eqref{eq_galg} expresses the Stokes determinant $\mathfrak g$ in terms of the polynomial $\mathsf P$. Combining the sign of the denominator in \autoref{lem_galg}(iii) with the curvature sign in \autoref{lem_thetasign}, we see that it suffices to prove $\mathsf P<0$ in order to obtain ${{\Theta}}>0$.

Take $g=k=1$ and write $t=\tanh d$ and $\varsigma=\operatorname{sech}d$. Given ${{\Omega_1}}<0$, $d>0$ and $\mathsf x>0$, define $\sigma=\sqrt\lambda>0$ by $q(\lambda)=\varsigma\mathsf x$, i.e.
\[
\sigma=\frac12\Big(-{{\Omega_1}}(\varsigma\mathsf x+t)+\sqrt{{{\Omega_1}}^2(\varsigma\mathsf x+t)^2+4(\varsigma\mathsf x+t)}\Big),
\]
and put ${{W_{\!I}}}=-\sigma-{{\Omega_1}}d$ and $q_\sigma=\sigma(2-{{\Omega_1}}\sigma)/(1-{{\Omega_1}}\sigma)^2$, the derivative of $q$ with respect to $\sigma$ (\eqref{eq_qsigma}). Then $\rho=\varsigma/\mathsf x-t$ and $\partial_\sigma\rho=-q_\sigma/\mathsf x^2$, and the tangency conditions $\mathcal E=\partial_\sigma\mathcal E=0$ read
\begin{equation}\label{eq_foldmanifold}
\coth(h-d)=-t+\frac{\varsigma}{\mathsf x}+\frac{{{W_{\!I}}} q_\sigma}{\mathsf x^2},\qquad {{\Omega_2}}={{\Omega_1}}+\frac{{{W_{\!I}}}^2q_\sigma}{\mathsf x^2}.
\end{equation}
For \(f={{W_{\!I}}}\rho\), the first right-hand side in \eqref{eq_foldmanifold} is \(\partial_W f\). A fold in the window is a double zero \(\lambda^\ast\in\mathcal J\) of \(\mathcal E\). The tangency equations \(\mathcal E=\partial_\sigma\mathcal E=0\) identify its unique triple \((\mathsf x,{{\Omega_1}},d)\) with \(\mathsf x=q(\lambda^\ast)/\varsigma>0\) and
\[
{{W_{\!I}}}>0,\qquad
\partial_W f=-t+\frac{\varsigma}{\mathsf x}+\frac{{{W_{\!I}}}q_\sigma}{\mathsf x^2}>1.
\]
Conversely, each such triple determines \(h>d\) and \({{\Omega_2}}\) by \eqref{eq_foldmanifold}. The resulting double root is a nondegenerate fold in the parameter \({{\Omega_2}}\) by \autoref{rem_foldwindow}.

The double-angle identities give $G_1=1/t$, $\Sigma_1=\varsigma/t$, $G_2=(1+t^2)/(2t)$, $\Sigma_2=\varsigma^2/(2t)$ and $K_2=(1+K_1^2)/(2K_1)$.

\begin{lemma}\label{lem_thetasign}
Let $g=k=1$. At every fold in the window, $\partial_W^2\widetilde\Delta<0$, so that, by \eqref{eq_thetaformula}, ${{\Theta}}=8\lambda^\ast\mathfrak g/\partial_W^2\widetilde\Delta$ has the sign of $-\mathfrak g$.
\end{lemma}
\begin{proof}
For $g=k=1$, \eqref{eq_det_jn} and \autoref{lem_convex}(i) give
\[
\widetilde\Delta=\frac{2{{W_{\!I}}} D_1}{\sinh d\,\sinh(h-d)}=-\frac{2(1+{{\Omega_1}}s)\,q}{t}\,\mathcal E .
\]
At a fold in the window, $\mathcal E=\partial_W\mathcal E=0$, so 
$\partial^2_W\widetilde\Delta=-2(1+{{\Omega_1}}s)q\,\partial_W^2\mathcal E/t.$
By \autoref{lem_convex}(iii), $\partial_W^2\mathcal E=\partial_\sigma^2\mathcal E>0$. The fold lies in $\mathcal J$, so $q>0$. Together with $1+{{\Omega_1}}s>0$ and $t>0$, these signs give $\partial_W^2\widetilde\Delta<0$.
\end{proof}

The Stokes determinant $\mathfrak g=\det[\widetilde J_5^{(4)}|\mathbf t^\flat]$ can be expressed algebraically in the following variables. Write
\[
\sigma=\sqrt\lambda=-s,\qquad \mathsf w:={{W_{\!I}}},\qquad \mathsf p:={{\Omega_1}}s=-{{\Omega_1}}\sigma,\qquad \mathsf q:=(1+\mathsf p)\,q=\sigma^2-t(1+\mathsf p).
\]
The determinant denominator in \eqref{eq_galg} uses $\mathsf E$, and the requirement $\coth(h-d)>1$ gives the threshold $\mathsf w_{\min}$. Define
\begin{gather*}
\mathsf E:=-\sigma^5(1-t^2)\,\mathsf q+\mathsf w\,t\,(2+\mathsf p)\big[\sigma^4(t^2-3)+3\sigma^2t(1+\mathsf p)-t^2(1+\mathsf p)^2\big],\\
\mathsf w_{\min}:=\frac{(\sigma^2-1-\mathsf p)\,\big(\sigma^2-t(1+\mathsf p)\big)}{(1-t)\,\sigma\,(2+\mathsf p)} .
\end{gather*}

\begin{lemma}\label{lem_galg}
Let $g=k=1$.
\begin{enumerate}[label=\textup{(\roman*)},leftmargin=2.3em]
\item There is a polynomial $\mathsf P\in\mathbb Z[\sigma,\mathsf p,t,\mathsf w]$, the same for all folds, with $1001$ terms and degrees $19,11,8,6$ in $\sigma,\mathsf p,t,\mathsf w$, such that at every fold in the window
\begin{equation}\label{eq_galg}
\mathfrak g=-\frac{(1-t^2)\,\mathsf P(\sigma,\mathsf p,t,\mathsf w)}{4\,\sigma\,t\,\mathsf q^6\,\mathsf E}.
\end{equation}
\item On the fold manifold, $\partial_Wf-1=(1-t^2)\,q_\sigma\,(\mathsf w-\mathsf w_{\min})/q^2$, where $q=\mathsf q/(1+\mathsf p)$.
\item $\mathsf E<0$ whenever $\sigma,\mathsf p,\mathsf w>0$, $0<t<1$ and $\mathsf q>0$.
\end{enumerate}
\end{lemma}
\begin{proof}
\smallskip\noindent\textit{Step 1. The determinant in physical variables.}
Work in $\mathbb Q(\sigma,\mathsf p,t,\mathsf w)$, set $c=1-t^2$, $s=-\sigma$, $W=\mathsf w$, $\Omega_1=-\mathsf p/\sigma$, $d=\sigma(\sigma+W)/\mathsf p$. The fold relations give
$K_1=c/q-t+Wc q_\sigma/q^2,\qquad \Omega_2=\Omega_1+W^2c q_\sigma/q^2 .$
Write the first-order surface amplitudes $a_1,e_1,\mathsf u_1$ as $\sqrt c$ times rational functions, using \eqref{eq_fold_kernel_amplitude}. The first two third-order residuals likewise contain one common $\sqrt c$; their cofactors contribute the second. All second-mode multipliers are rational, since $G_2=(1+t^2)/(2t)$, $\Sigma_2=c/(2t)$ and $K_2=(K_1+K_1^{-1})/2$. Eliminate $a_2,b_2,r_2$ from rows 1, 3 and 4 of \eqref{eq_order2mode2}, solve the remaining two rows, and form the cofactor expression of \autoref{lem_b2}(c) with $M_2=0$. Perform these operations over the rational function field, collecting like powers after clearing the common denominator. Multiplying the result by $-4\sigma t\mathsf q^6\mathsf E/c$, we obtain a reduced fraction with denominator one, integer coefficients and the support and degrees asserted in (i). All intervening divisions are valid at physical folds: $q,\sigma,\mathsf p,t,d,K_1$ are positive, and the mode-two system is invertible by \autoref{lem_sl} and \eqref{eq_det_jn}.

\smallskip\noindent\textit{Step 2. The derivative threshold.}\enspace
Substitute \(\varsigma\mathsf x=q\) in \eqref{eq_foldmanifold}. With \(q=\mathsf q/(1+\mathsf p)\), \(q-1+t=(\sigma^2-1-\mathsf p)/(1+\mathsf p)\) and \(q_\sigma=\sigma(2+\mathsf p)/(1+\mathsf p)^2\) from \eqref{eq_qsigma}, this gives
\begin{equation*}
\begin{aligned}
q^2(\partial_W f-1)
&=(1-t^2)\Bigl(\mathsf w q_\sigma-\frac{q(q-1+t)}{1-t}\Bigr)
=(1-t^2)q_\sigma(\mathsf w-\mathsf w_{\min}).
\end{aligned}
\end{equation*}

\smallskip\noindent\textit{Step 3. The denominator sign.}\enspace
We set $u=t(1+\mathsf p)$. Since $0<u<\sigma^2$ and $0<t<1$, the bracket in $\mathsf E$ equals
$\sigma^4(t^2-1)-(\sigma^2-u)(2\sigma^2-u)<0$.
Its multiplier is positive and the other summand is negative, giving $\mathsf E<0$.
\end{proof}

\begin{theorem}\label{thm_thetapos}
Let $g,k>0$. Suppose that $\lambda^\ast$ is a double root of $D_k$ in the component $\{s<0<{{W_{\!I}}}\}$ and that the normalized kernel profile at $\lambda^\ast$ is positive.

Then $\lambda^\ast$ is a nondegenerate fold in $\mu={{\Omega_2}}$ with fold coefficient ${{\Theta}}>0$ in \eqref{eq_theta}.
At the fold, ${{\Omega_1}}<{{\Omega_2}}$ and ${{\Omega_1}}<0$.
\end{theorem}
\begin{proof}[Computer-assisted proof in exact integer arithmetic]
\smallskip\noindent\textit{Step 1. Reduction to a polynomial sign.}\enspace
Replacing $(x,y)$ by $(kx,ky)$ and $\psi$ by $k^{3/2}g^{-1/2}\psi$ maps the problem with constants $(g,k,{{\Omega_1}},{{\Omega_2}},h,d)$ to the one with constants $g=k=1$, ${\Omega}_i/\sqrt{gk}$, $kh$ and $kd$. The scaling maps laminar flows to laminar flows, roots and folds to roots and folds, and positive kernels to positive kernels, with $\lambda\mapsto(k/g)\lambda$ and $\zeta_1\mapsto k\zeta_1$. By \autoref{lem_theta}(ii), ${{\Theta}}$ is multiplied by the positive factor $g^{-2}$. We therefore take $g=k=1$. The velocity signs and
${{\Omega_1}}d=s-{{W_{\!I}}}$ give ${{\Omega_1}}<0$, and the second
identity in \eqref{eq_foldmanifold} gives ${{\Omega_2}}>{{\Omega_1}}$.

After scaling, define the polynomial sign domain by
\begin{equation}\label{eq_fold_polynomial_domain}
\mathcal D_{\rm fold}:=\big\{(\sigma,\mathsf p,t,\mathsf w):\ \sigma>0,\ \mathsf p>0,\ 0<t<1,\ \sigma^2>t(1+\mathsf p),\ \mathsf w>\max(0,\mathsf w_{\min})\big\}.
\end{equation}
By \autoref{rem_foldwindow} the double root is a nondegenerate fold. By \autoref{lem_thetasign}, \autoref{lem_galg}(i),(iii) and $\mathsf q>0$, we have 
\[
\operatorname{sgn}{{\Theta}}=-\operatorname{sgn}\mathfrak g=-\operatorname{sgn}\mathsf P(\sigma,\mathsf p,t,\mathsf w).
\]
By the discussion after \eqref{eq_foldmanifold} and \autoref{lem_galg}(ii), the point $(\sigma,\mathsf p,t,\mathsf w)$ of a fold in the window lies in the domain \eqref{eq_fold_polynomial_domain}.
The domain $\mathcal D_{\rm fold}$ treats $t$ as independent of $\sigma,\mathsf p,\mathsf w$ subject to \eqref{eq_fold_polynomial_domain}. At an actual fold, $t=\tanh\big(\sigma(\sigma+\mathsf w)/\mathsf p\big)$. We prove
$\mathsf P<0\quad\text{on }\mathcal D_{\rm fold}.$

Form the sparse coefficient dictionary of $\mathsf P$ as in \autoref{lem_galg}. The four exponents are the keys and the coefficients are signed integers.

\smallskip\noindent\textit{Step 2. Rational parametrizations of the fold domain.}\enspace
Let $\theta_1,\dots,\theta_4$ be variables. The following maps cover $\mathcal D_{\rm fold}$; each parameter range is mapped bijectively onto its corresponding part of $\mathcal D_{\rm fold}$.
\begin{itemize}
\item[(Ia)] $\mathcal D_{\rm fold}\cap\{\sigma\le1\}$: 
\begin{equation*}
\sigma=\frac1{1+\theta_1},\quad
\mathsf p=\theta_2,\quad
t=\frac{\theta_3}{1+\theta_3}\cdot\frac{\sigma^2}{1+\theta_2},\quad
\mathsf w=\theta_4,
\end{equation*}
 with $\theta_1\ge0$ and $\theta_2,\theta_3,\theta_4>0$. Here $\sigma^2<1+\mathsf p$, so $\mathsf w_{\min}<0$.
\item[(Ib)] $\mathcal D_{\rm fold}\cap\{\sigma\ge1,\ \sigma^2\le1+\mathsf p\}$: 
\begin{equation*}
\sigma=1+\theta_1,\quad
\mathsf p=\sigma^2(1+\theta_2)-1,\quad
t=\frac{\theta_3}{(1+\theta_3)(1+\theta_2)},\quad
\mathsf w=\theta_4,
\end{equation*}
 with $\theta_1,\theta_2\ge0$, $\theta_1+\theta_2>0$ and $\theta_3,\theta_4>0$. Here $\mathsf w_{\min}\le0$.
\item[(II)] $\mathcal D_{\rm fold}\cap\{\sigma^2>1+\mathsf p\}$: 
\begin{equation*}
\sigma=1+\theta_1,\quad
\mathsf p=\frac{(\sigma^2-1)\theta_2}{1+\theta_2},\quad
t=\frac{\theta_3}{1+\theta_3},\quad
\mathsf w=\mathsf w_{\min}+\theta_4,
\end{equation*}
 with $\theta_1,\dots,\theta_4>0$. Here $\sigma^2>1+\mathsf p>t(1+\mathsf p)$ holds automatically, and $\mathsf w_{\min}>0$.
\end{itemize}
For reproducibility, use the following denominator convention. To substitute $v=N/D$ in $P=\sum_{j=0}^m P_jv^j$, form $\sum_{j=0}^mP_jN^jD^{m-j}$; the specified $m$ may exceed the degree. Each denominator is positive. Apply the substitutions in the following order:
\begin{itemize}[nosep]
\item[(Ia)] first substitute $t=\theta_3\sigma^2/((1+\theta_3)(1+\mathsf p))$ with $m=8$, then $\sigma=1/(1+\theta_1)$ with $m=35$ and rename $\mathsf p,\mathsf w$.
\item[(Ib)] first set $\mathsf p=\sigma^2(1+\theta_2)-1$, then substitute the displayed $t$ with $m=8$ and set $\sigma=1+\theta_1$.
\item[(II)] first substitute $\mathsf w=(N+D\theta_4)/D$ with $m=6$, where $N=(\sigma^2-1-\mathsf p)(\sigma^2-t(1+\mathsf p))$ and $D=(1-t)\sigma(2+\mathsf p)$. Substitute $t=\theta_3/(1+\theta_3)$ with $m=14$, then $\mathsf p=(\sigma^2-1)\theta_2/(1+\theta_2)$ with $m=23$, and finally $\sigma=1+\theta_1$.
\end{itemize}
The three resulting polynomials have the sign of $\mathsf P$ on their respective regions and have $14\,503$, $11\,122$ and $46\,866$ nonzero coefficients.

\smallskip\noindent\textit{Step 3. The integer sign certificate and its verification.}\enspace
Expanding by integer convolution, we obtain the following coefficient tests: all coefficients of $\mathsf P_{\rm Ia}$ and $\mathsf P_{\rm Ib}$ are negative. The same expansions show that $\mathsf P_{\rm Ia}$ contains monomials free of $\theta_1$, and $\mathsf P_{\rm Ib}$ contains monomials free of $\theta_1$ and monomials free of $\theta_2$. On each allowed face $\theta_i=0$, a monomial independent of $\theta_i$ still has strictly positive value and a negative coefficient. Thus both polynomials are strictly negative throughout their parameter sets.

$\mathsf P_{\rm II}$ has $46{,}482$ negative and $384$ positive coefficients. For each positive term $c\,\theta^e$, the certificate uses pairs of monomials $\theta^{e'}$, $\theta^{e''}$ with negative coefficients and $e=\frac12(e'+e'')$. On $(0,\infty)^4$, $2\theta^e\le\theta^{e'}+\theta^{e''}$ because
$\theta^{e'}+\theta^{e''}-2\theta^e=(\theta^{e'/2}-\theta^{e''/2})^2.$
We now assign the positive coefficients to $390$ pairs and check the weights using integers. If $c'$ is a negative coefficient, the sum of the weights $\alpha_j$ assigned to pairs using this monomial does not exceed $2|c'|$. Consequently
\[
2\,\mathsf P_{\rm II}=\mathsf R-\sum_{j=1}^{390}\alpha_j\big(\theta^{e'_j/2}-\theta^{e''_j/2}\big)^2,\qquad \alpha_j\in\mathbb N,
\]
where all coefficients of $\mathsf R\in\mathbb Z[\theta]$ are $\le0$. Of these coefficients, $46{,}468$ are negative and $398$ are zero. Hence $\mathsf P_{\rm II}<0$ on $(0,\infty)^4$.

Here is a deterministic construction of the weights, using exact rational comparisons and lexicographic order on four-component exponent vectors. It also provides a way to reproduce the certificate from the polynomial alone. Initially give each negative exponent $a$ budget $B_a=-2c_a$. For each positive exponent $e$ let $P_e$ consist of pairs $(a,b)$ of negative exponents with $a+b=2e$ and $a< b$. Repeat until all positive exponents have been processed:
\begin{enumerate}[label=\arabic*.,nosep]
\item Among the remaining $e$, minimize $\sum_{(a,b)\in P_e}\min(B_a,B_b)/c_e$, breaking ties by the smallest $e$. We set $r=c_e$.
\item Order $P_e$ by decreasing $\min(B_a,B_b)$ at the start of the pass, breaking ties by increasing $a$. For each pair assign $v=\min(r,B_a,B_b)$, and subtract $v$ from each of $r,B_a,B_b$. Record the row if $v>0$; stop the pass when $r=0$.
\item Check that $r=0$, then remove $e$. At the end verify the midpoint and budget identities above. The three substitutions in Step~2 and this allocation require only arithmetic on integer coefficients. The resulting counts include both admissible boundary faces in (Ib) through the separate face tests above.
\end{enumerate}
\end{proof}

\begin{proof}[Proof of \autoref{thm_c}]

By \autoref{thm_thetapos}, each double root in the window is a nondegenerate
fold in ${{\Omega_2}}$ with ${{\Theta}}>0$. The root pair lies on the side
${{\Omega_2}}<{{\Omega_2}}^\ast$ by \autoref{prop_pairing}. The small loop
and, when both critical heights exist at the fold, the four stagnation
points on each nontrivial wave follow from \autoref{thm_b}.

\end{proof}

\subsection{The type of the bifurcations at simple roots}\label{sec_type}

We next consider a simple root of $D_1$. With the interface coefficient $\tau=\zeta_1(w)$ as amplitude, the sign of $\lambda_2$ in \eqref{eq_pitchfork} determines the direction of the local branch. We calculate this sign both near the fold and at parameters farther from it.

Let $(\lambda^\ast,0)\in{\mathcal O}$ satisfy \eqref{eq_h} and \eqref{eq_hsimple}. The first interface cosine coefficient $\tau$ is a nonzero multiple of $a$ by \eqref{eq_local_interface_amplitude}. Half-period symmetry in \autoref{cor_local} therefore makes $\lambda$ analytic and even in $\tau$, with expansion
\begin{equation}\label{eq_pitchfork}
\lambda(\tau)=\lambda^\ast+\lambda_2\tau^2+O(\tau^4).
\end{equation}
The pitchfork is nondegenerate if $\lambda_2\ne0$. At a simple positive
root $\lambda_\pm$, write $\lambda_2^{(\pm)}$. If $\lambda_2^{(\pm)}\neq0$, the nontrivial solutions near $(\lambda_\pm,0)$ lie in
\[
\begin{cases}
\{\lambda>\lambda_\pm\}&\text{if }\lambda_2^{(\pm)}>0,\\
\{\lambda<\lambda_\pm\}&\text{if }\lambda_2^{(\pm)}<0.
\end{cases}
\]
We then say that the branch \emph{bends to the right}, resp.\ \emph{to the left}. For a pair \(\lambda_-<\lambda_+\), the inequality \(\partial_\lambda{{W_{\!I}}}=-1/(2\sqrt\lambda)<0\) puts \(\lambda_+\) at the smaller zero when \(\mathcal E\) is viewed as a convex function of \({{W_{\!I}}}\). Hence
\[
\left.\partial_W\mathcal E\right|_{\lambda_+}<0<
\left.\partial_W\mathcal E\right|_{\lambda_-}.
\]
In cases A, B and C of \autoref{tab_roots}, $\lambda_-$ is the surface-mode root and $\lambda_+$ the interface-mode root. Both branches \emph{bend towards the partner} if 
$\lambda_2^{(+)}<0<\lambda_2^{(-)}.$

To compute $\lambda_2^{(\pm)}$ in \eqref{eq_pitchfork}, evaluate the dispersion factors at the simple root $\lambda_\pm$. With $g=k=1$, $t=\tanh d$, $K=\coth(h-d)$, $\sigma=\sqrt\lambda$, $\mathsf p={{\Omega_1}}s=-{{\Omega_1}}\sigma$, $\mathsf w={{W_{\!I}}}>0$, $q=q(\lambda)$, $\rho=(1-t^2)/q-t$ and $q_\sigma$ as in \eqref{eq_qsigma}, and $\partial_Wf=\rho+\mathsf w(1-t^2)q_\sigma/q^2$, the root equation $\mathcal E=0$ reads ${{\Omega_2}}={{\Omega_1}}+K{{W_{\!I}}}-{{W_{\!I}}}\rho$. The definition of $\mathcal E$ preceding \autoref{lem_convex} gives $\partial_W\mathcal E=\partial_Wf-K$.

\begin{lemma}\label{lem_lambda2}
Let $g=k=1$ and let $\lambda$ be a simple positive root of $D_1$ in the component $\{s<0<{{W_{\!I}}}\}$. Let $\mathfrak g=\det[\widetilde J_5^{(4)}|\mathbf t^\flat]$ be the Stokes determinant of \textup{\autoref{lem_b2}}, evaluated at this simple root. Then
\begin{equation}\label{eq_lambda2formula}
\lambda_2=\frac{\sigma\big(\mathsf w(K+\rho)-{{\Omega_1}}\big)}{2(h-d)}+2\sigma\,\mathsf B,\qquad \mathsf B:=\frac{t\,\mathfrak g}{2(1+\mathsf p)\,q\,(K-\partial_Wf)}+\frac{\mathsf w\,(\partial_Wf+\rho)-{{\Omega_1}}}{4d}.
\end{equation}
The first term on the right-hand side of the formula for $\lambda_2$ is positive.
\end{lemma}
\begin{proof}
\smallskip\noindent\textit{Step 1. The lower flux and the determinant derivative.}\enspace
The coefficient $m_{2,2}$ of $\tau^2$ in $m_2$ determines $\lambda_2$:
\[
\lambda_2=\frac{2\sigma}{h-d}\,m_{2,2}.
\]
The velocity signs give laminar admissibility, and \autoref{lem_sl} gives $D_n(\lambda)>0$ for $n\ge2$. Thus \eqref{eq_h} and \eqref{eq_hsimple} hold at this root. By \autoref{lem_b3} and~\autoref{lem_b2}(b), 
\[
m_{2,2}=L_2^\flat-(h-d)C_2^\flat+(h-d)\frac{\mathfrak g}{\partial_W\widetilde\Delta}.
\]

By the proof of \autoref{lem_thetasign}, $\widetilde\Delta=-2(1+\mathsf p)q\,\mathcal E/t$. At the root, $\mathcal E=0$ and $\partial_W\mathcal E=\partial_W f-K\ne0$ because the root is simple. Differentiating the factorization, we obtain
\begin{equation}\label{eq_simple_root_determinant_derivative}
\partial_W\widetilde\Delta=2(1+\mathsf p)q(K-\partial_Wf)/t\ne0.
\end{equation}

\smallskip\noindent\textit{Step 2. The explicit coefficient and its sign.}\enspace
It remains to evaluate $L_2^\flat$ and $C_2^\flat$ in the expression for the lower flux. This will give \eqref{eq_lambda2formula}.
By \autoref{lem_b2}(a), $L_2^\flat=-\mathsf g_4=-\frac12K r_1-\frac14{{\Omega_2}}$ and $C_2^\flat=A_2^\flat+\mathsf g_5$, $A_2^\flat=(\mathsf g_3-\mathsf g_1)/d$, with the kernel $b_1=r_1=-{{W_{\!I}}}=-\mathsf w$, $a_1=-\mathsf w\varsigma(q+t)/q$, $e_1=-a_1/s$ of \autoref{lem_b1}. First, $\mathsf g_5=\frac12(b_1-r_1)=0$, so $C_2^\flat=A_2^\flat$. The root satisfies $\mathcal E=0$, i.e.\ ${{\Omega_2}}={{\Omega_1}}+K\mathsf w-\mathsf w\rho$, so
\[
L_2^\flat=\frac14\big(2K\mathsf w-{{\Omega_2}}\big)=\frac14\big(\mathsf w(K+\rho)-{{\Omega_1}}\big).
\]
The mode-zero coefficient of the lower affine slope is
\begin{equation}\label{eq_mode_zero_height_correction}
C_2^\flat=A_2^\flat=\frac{{{\Omega_1}}-\mathsf w\,(\partial_Wf+\rho)}{4d}.
\end{equation}
The mode-zero equation for $A_2^\flat$ gives
\[
4d A_2^\flat=4(\mathsf g_3-\mathsf g_1)={{\Omega_1}}+2(\Sigma_1a_1-G_1b_1)-e_1\big(2(G_1a_1-\Sigma_1b_1)+{{\Omega_1}}e_1\big).
\]
Substituting \(G_1=1/t\), \(\Sigma_1=\varsigma/t\) and \(\varsigma^2=1-t^2\), we obtain
\begin{equation*}
\begin{aligned}
\Sigma_1a_1-G_1b_1
&=-\frac{\mathsf w}{t}\Bigl(\frac{(1-t^2)(q+t)}{q}-1\Bigr)
=-\mathsf w\Bigl(\frac{1-t^2}{q}-t\Bigr)=-\mathsf w\rho,\\
G_1a_1-\Sigma_1b_1
&=-\frac{\mathsf w\varsigma}{t}\Bigl(\frac{q+t}{q}-1\Bigr)
=-\frac{\mathsf w\varsigma}{q}.
\end{aligned}
\end{equation*}
Using \(q+t=\sigma^2/(1+\mathsf p)\), \(a_1=-\mathsf w\varsigma\sigma^2/((1+\mathsf p)q)\), \(e_1=-a_1/s=a_1/\sigma\) and \({{\Omega_1}}\sigma=-\mathsf p\), we obtain
\begin{equation*}
\begin{aligned}
-e_1\Bigl(-\frac{2\mathsf w\varsigma}{q}+{{\Omega_1}}e_1\Bigr)
&=\frac{2\mathsf w\varsigma a_1}{\sigma q}-\frac{{{\Omega_1}}a_1^2}{\sigma^2}
=-\frac{\mathsf w^2(1-t^2)\sigma\bigl(2(1+\mathsf p)-\mathsf p\bigr)}{(1+\mathsf p)^2q^2}
=-\frac{\mathsf w^2(1-t^2)q_\sigma}{q^2}.
\end{aligned}
\end{equation*}
By \eqref{eq_qsigma}, \(q_\sigma=\sigma(2+\mathsf p)/(1+\mathsf p)^2\), and \(\partial_Wf=\rho+\mathsf w(1-t^2)q_\sigma/q^2\). Thus
\begin{equation*}
\begin{aligned}
4dA_2^\flat
&={{\Omega_1}}-2\mathsf w\rho-\mathsf w^2(1-t^2)q_\sigma/q^2
={{\Omega_1}}-\mathsf w(\partial_Wf+\rho).
\end{aligned}
\end{equation*}
Inserting these into the formula for $m_{2,2}$, we obtain $m_{2,2}=\frac14\big(\mathsf w(K+\rho)-{{\Omega_1}}\big)+(h-d)\mathsf B$, and multiplying by $2\sigma/(h-d)$ gives \eqref{eq_lambda2formula}. 

The first term in \eqref{eq_lambda2formula} is positive since $K>1$, $\mathsf w>0$, ${{\Omega_1}}<0$ and
$\rho=(1-t^2)/q-t>-t>-1,\qquad K+\rho>0.$
\end{proof}

\begin{theorem}\label{thm_type}
Assume that ${{\Omega_1}}<\min(0,{{\Omega_2}})$. Let $\lambda_-<\lambda_+$ be the two simple positive roots of $D_1$ in $\{s<0<{{W_{\!I}}}\}$.

The existence of the two roots gives ${{\Omega_1}}<{{\Omega_2}}<{{\Omega_2}}^\ast$ by \autoref{prop_pairing}.
\begin{enumerate}[label=\textup{(\alph*)},leftmargin=2.3em]
\item The local curve at each $\lambda_\pm$ is a pitchfork with expansion \eqref{eq_pitchfork}. For $g=k=1$, the pitchfork coefficient $\lambda_2^{(\pm)}$ is given by \eqref{eq_lambda2formula}.
\item Near the fold, both branches bend towards each other. As ${{\Omega_2}}\uparrow{{\Omega_2}}^\ast$,
$\lambda_2^{(\pm)}(\lambda_\mp-\lambda_\pm)\to{{\Theta}}>0.$
By \autoref{thm_loops}, the fold is a $\mathbb Z_2$-symmetric isola.
Only trivial solutions occur locally at the fold itself. A small loop
appears as ${{\Omega_2}}$ decreases through ${{\Omega_2}}^\ast$.
\item The branch at $\lambda_+$ can bend away from the smaller root $\lambda_-$. For $({{\Omega_1}},{{\Omega_2}},h,d)=(-1,1,4,3)$ and $g=k=1$,
$\lambda_2^{(+)}>0.$
For $g=k=1$ and $(\Omega_1,h,d)=(-1,4,3)$ fixed, the coefficient $\lambda_2^{(+)}$ changes sign at least twice as $\Omega_2$ varies. Distinct sign-changing zeros lie in $(0.2100,0.2102)$ and $(1.7115,1.7117)$.
\end{enumerate}
\end{theorem}

\begin{proof}
\smallskip\noindent\textit{Step 1. Local bifurcations and the isola near the fold.}\enspace
At both roots the velocity signs give laminar admissibility, and \autoref{lem_sl} gives higher-mode nonresonance. At either simple root, the linear change of parameter \eqref{eq_local_interface_amplitude} carries the even analytic function $\lambda(a)$ of \autoref{cor_local} to the expansion \eqref{eq_pitchfork} in the physical amplitude $\tau$. For a positive root with an upper critical layer and $g=k=1$, $\lambda_2$ is given by \autoref{lem_lambda2}.

The roots approach the nondegenerate fold of \autoref{rem_foldwindow} as ${{\Omega_2}}\uparrow{{\Omega_2}}^\ast$. The bending limit follows from \autoref{lem_theta}(ii) and \autoref{thm_thetapos}, giving $\lambda_2^{(+)}<0<\lambda_2^{(-)}$ near the fold. The isola is the positive-${{\Theta}}$ case of \autoref{thm_loops}.

\smallskip\noindent\textit{Step 2. Two changes in bending at the larger root.}
For family C, apply \eqref{eq_pitchfork} and \autoref{lem_b2} to the simple root with smaller $W_I$, hence larger $\lambda$. The enclosures obtained by the rational interval procedure of \autoref{app_cap} are
\[
\begin{array}{c|rrrr}
\mu&0.2100&0.2102&1.7115&1.7117\\ \hline
10^4\lambda_2^{(+)}
&(-1.271,-1.268)&(1.825,1.828)&(2.767,2.770)&(-5.424,-5.420).
\end{array}
\]
At $\mu=1$ the same test gives $0.44609017<\lambda_2^{(+)}<0.44609019$.

For $({{\Omega_1}},h,d)=(-1,4,3)$, the pair exists for ${{\Omega_2}}\in({{\Omega_1}},{{\Omega_2}}^\ast)$ with ${{\Omega_2}}^\ast=1.83073\ldots$ (\autoref{thm_certified}), and ${{\Omega_2}}=1$ lies in this interval. On this interval the root \(\lambda_+\) is simple, so the implicit function theorem for \(\mathcal E\) gives real-analytic dependence on \({{\Omega_2}}\). The coefficient \(\lambda_2^{(+)}\) in \eqref{eq_lambda2formula} is also real-analytic. Indeed \(q>0\) on \(\mathcal J\), and \(K-\partial_Wf=-\partial_W\mathcal E>0\) at the smaller zero in \({{W_{\!I}}}\). Also \(1+\mathsf p>0\), \(d>0\), \(h-d>0\), and \(\Sigma_1\) and \(s\) are nonzero. The remaining determinant in \(\mathfrak g\) is nonzero by \autoref{lem_sl},
\[
\det\widetilde J_5(2)=\frac{4{{W_{\!I}}}D_2}{\sinh 2d\,\sinh 2(h-d)}\ne0.
\]
The sign pattern is $-,+,+,-$. Continuity gives a zero in each of the two disjoint brackets. Analyticity and the strict endpoint signs ensure at least one sign-changing zero in each interval. At these parameters the quadratic term vanishes; its vanishing alone does not determine the quartic term or the total number of zeros.
\end{proof}

\begin{proof}[Proof of \autoref{thm_d}]

At an admissible simple mode-one root satisfying \eqref{eq_h} and
\eqref{eq_hsimple}, half-period symmetry gives \eqref{eq_pitchfork} in the
interface amplitude $\tau=\zeta_1(w)$. In the positive-kernel window
\autoref{lem_lambda2} computes $\lambda_2$ for $g=k=1$. The scaling in the
proof of \autoref{thm_thetapos} sends $\lambda$ to $(k/g)\lambda$ and
$\zeta_1$ to $k\zeta_1$, so the original coefficient is $gk$ times the
coefficient in the scaled problem.

The limit and isola in \autoref{thm_type}(b) give the near-fold bending
and local loop formation. The certified sign and the change of direction
for $({{\Omega_1}},h,d)=(-1,4,3)$ follow from \autoref{thm_type}(c).

\end{proof}

\subsection{Evaluation of the fold coefficient}\label{sub_11_evaluation_of_the_fold_coefficient}

To illustrate the preceding results, we give four sets of parameters. Interval evaluation of the signed Stokes determinant encloses the fold coefficient in each case, and the laminar profile gives the two critical heights.

\begin{theorem}\label{thm_certified}
Let $g=k=1$ and $\mu={{\Omega_2}}$. Each tabulated triple $({{\Omega_1}},h,d)$ has a nondegenerate fold $(\lambda^\ast,{{\Omega_2}}^\ast)$ in the positive-kernel window. At each fold the coefficient in \eqref{eq_thetaformula} satisfies ${{\Theta}}>0$, and the laminar flow has critical heights $W(y_c^\pm)=0$ in both layers.

The two-root side is ${{\Omega_2}}<{{\Omega_2}}^\ast$. The values of the fold parameters and the coefficient are rounded to $15$ significant digits.
\begin{center}\small
\begin{tabular}{@{}c c c c@{}}
\toprule
$({{\Omega_1}},h,d)$ & $\lambda^\ast$ & ${{\Omega_2}}^\ast$ & ${{\Theta}}$\\
\midrule
$(-1,5,3)$ & $2.89361421181034$ & $1.47181181033377$ & $0.935131984450719$\\
$(-1,4,3)$ & $2.87503052443117$ & $1.83073377082644$ & $0.915354610544684$\\
$(-2,7,2.8)$ & $6.86160530682888$ & $3.53211245553920$ & $16.0798930575829$\\
$(-0.5,7,4.2)$ & $1.68353508130286$ & $1.07708491896854$ & $0.0607161614485543$\\
\bottomrule
\end{tabular}
\end{center}
\end{theorem}
\begin{proof}
Take $g=k=1$. On $q>0$, $W=W_I>0$, solve the root equation for $\mu=\Omega_1+W(K-\rho)$, where $\sigma=-W-\Omega_1d$. Its derivative is $K-\rho-W(1-t^2)q_\sigma/q^2$. For each row, \autoref{app_cap} gives a rational interval on which $q,W,\sigma>0$ and the derivative has opposite signs at the endpoints. The strictly convex reduction in \autoref{prop_pairing} and nondegeneracy in \autoref{rem_foldwindow} identify the interior fold and the two-root side. The interval evaluation gives $W-\mu(h-d)<0$, hence critical heights in both layers, and evaluates $8s^2\mathfrak g/\widetilde\Delta_{WW}$ on the entire bracket. Every denominator excludes zero. The resulting enclosures certify the displayed rounding and strict positivity in each row.
\end{proof}

\subsection{Exact verification and interval enclosures}\label{app_cap}
The polynomial and integer weights can be reconstructed from the formulas and the ordered allocation in the proof of \autoref{thm_thetapos}. In the code accompanying the article, \texttt{fold\_polynomial.py} generates the numerator, and \texttt{fold\_sign.py} checks it after substitution. The option \texttt{--write} regenerates the weights. The independent residual check projects the five boundary equations onto the first harmonic before simplification.
\paragraph{Interval enclosures for finite parameters.}
All interval endpoints are integer multiples of $10^{-55}$, and arithmetic is rounded outwards. For $x=nd$ or $n(h-d)$ the program encloses $e^x$ from below by $T_{200}=\sum_{j=0}^{200}x^j/j!$ and from above by
\[
T_{200}+\frac{x^{201}}{201!}\frac1{1-x/202},\qquad 0<x<202.
\]
In this way we obtain initial enclosures for the hyperbolic functions. Each fold bracket has radius $10^{-23}$ about the following centre in $W_I$, in the order of \autoref{thm_certified}:
\[
\begin{array}{ll}
1.2989373286647134324724071461,&1.3044085030788888294565151367,\\
2.9805333926868089736598970931,&0.8024888897189115473380617942 .
\end{array}
\]
To enclose a zero of $d\mu/dW$ we check that the endpoint signs are opposite. For a simple root the code tests opposite signs of $\widetilde\Delta$ at the endpoints and distinguishes the two roots by the sign of $d\mu/dW$. It evaluates the five-dimensional determinants on the whole interval, testing every denominator before division. For the four sign tests of \autoref{thm_type}(c), use the same radius and the following centres, ordered by increasing $\mu$:
\[
\begin{array}{ll}
0.5271859741761819262613779717393,&0.5272734224046383829717191361564,\\
1.208090172380479787016709445626,&1.208200555866705120559960475746 .
\end{array}
\]
For the larger root at $\mu=1$ use $0.8748253983120974066440877218$. All enclosures use the signed physical variables in the Stokes equations.

\section{Global graph continuation and laminar return}\label{sec_globalg}
\subsection{The vertical velocity}\label{sub_05_the_vertical_velocity}

We first derive the equations satisfied by the vertical velocity $V=-\psi_x$. The calculation applies also to overhanging waves. In the graph setting the physical equations, regularity and absence of boundary stagnation follow from \autoref{prop_operator}(iii).

\begin{lemma}\label{lem_v}
Let $\psi$ be the stream function of an even $2\pi$-periodic two-layer flow satisfying \eqref{eq_intro_model} with
\[
\psi^\pm\in C^{3,\alpha}(\overline{{\Omega}^\pm}),\qquad
\nabla\psi\ne0\quad\text{on }{{\mathcal S}}\cup{{\mathcal I}}.
\]
Let $n$ be the unit normal pointing out of the fluid on ${{\mathcal S}}$ and into ${\Omega}^+$ on ${{\mathcal I}}$. The vertical velocity $V=-\psi_x$ has the following properties.
\begin{enumerate}[label=\textup{(\roman*)},leftmargin=2.3em]
\item ${\Delta} V=0$ in ${\Omega}^+$ and in ${\Omega}^-$; $V$ is odd in $x$ about $x=0$ and about $x=\pi$, so $V=0$ on $\{x=0\}$ and $\{x=\pi\}$; $V=0$ on ${{\mathcal B}}$.
\item On ${{\mathcal I}}$, $V$ is continuous, and $[\nabla V]=0$ wherever $V=0$. Near points where $\psi_y\ne0$, ${{\mathcal I}}$ is a graph $y=b(x)$ with $b'=V/\psi_y$, and there
\begin{equation}\label{eq_vy_jump}
[V_y]=\frac{({{\Omega_1}}-{{\Omega_2}})\,b'}{1+b'^2}.
\end{equation}
\item The following equation holds on the entire free surface \({{\mathcal S}}\),
\begin{equation}\label{eq_vsurf}
(\psi_y^2-V^2)V_y-2\psi_yVV_x=(g+{{\Omega_1}}\psi_y)V .
\end{equation}
Where \(\psi_y\ne0\), the surface is locally \(y=\eta(x)\) with \(\eta'=V/\psi_y\).
\end{enumerate}
\end{lemma}
\begin{proof}
Differentiating \eqref{eq_intro_model}, we obtain the equations for $V=-\psi_x$. By the Poisson equations \eqref{eq_intro_model} we have ${\Delta} V=0$ in both layers. The symmetries imply that $V$ is odd about $x=0$ and $x=\pi$. The constant trace in \eqref{eq_intro_model} also gives $V=0$ on ${{\mathcal B}}$.

On ${{\mathcal I}}$, continuity of the gradient in \eqref{eq_intro_model} gives continuity of $V$. One tangential derivative of the gradient trace shows that $[D^2\psi]$ vanishes on tangent vectors. This symmetric matrix has trace $[{\Delta}\psi]={{\Omega_1}}-{{\Omega_2}}$ by \eqref{eq_intro_model}, which gives the jump and its negative $x$-row
\begin{equation*}
\begin{aligned}
[D^2\psi]&=({{\Omega_1}}-{{\Omega_2}})\,n\otimes n,\\
[\nabla V]&=-({{\Omega_1}}-{{\Omega_2}})\,n_x\,n.
\end{aligned}
\end{equation*}
At an interface zero of $V$, nonstagnation and the constant trace give $n_x=0$, hence $[\nabla V]=0$.

On either boundary, the implicit function theorem gives a local graph where $\psi_y\ne0$, and differentiation of the constant trace gives the corresponding slope. On a graph piece of ${{\mathcal I}}$, substituting a unit normal proportional to $(-b',1)$ into $[\nabla V]$ gives \eqref{eq_vy_jump}. For the corresponding scalar equations and interface relations, see \cite[Lemma~5.1 and (5.5a)--(5.5b)]{DMSW}.

On ${{\mathcal S}}$, the constant trace of $\psi$ makes $(\psi_y,V)$ tangent to the boundary. Differentiate the Bernoulli condition \eqref{eq_intro_model} along this vector and use $\psi_{xy}=-V_y$ and $\psi_{yy}={{\Omega_1}}+V_x$ from \eqref{eq_intro_model} to obtain
\begin{equation*}
\begin{aligned}
0&=\frac{1}{2}(\psi_y\partial_x+V\partial_y)(\psi_y^2+V^2+2gy)\\
 &=-(\psi_y^2-V^2)V_y+2\psi_yVV_x+(g+{{\Omega_1}}\psi_y)V .
\end{aligned}
\end{equation*}
This gives \eqref{eq_vsurf} at every surface point.
\end{proof}

\subsection{Orientation}\label{sub_05_orientation}

For positive orientation we require \(V=-\psi_x>0\) in the fluid, on the interface and on the free surface for \(0<x<\pi\), as well as the boundary derivative signs in \autoref{def_n}. These conditions determine the geometry of the profiles. Formula \eqref{eq_nflat} allows us to check them on the fixed strips.

We use the normalization \(k=1\). For a graph solution, the symmetry lines are \(x=j\pi\), \(j\in\mathbb Z\). We define the half-period cell between \(x=0\) and \(x=\pi\) and its boundary arcs by
\[
\begin{aligned}
D^\pm&={\Omega}^\pm\cap\{0<x<\pi\},&
\overline D&=\overline{{\Omega}_\eta}\cap\{0\le x\le\pi\},\\
{{\mathcal I}}_\pi&=\{(x,b(x)):0<x<\pi\},&
{{\mathcal S}}_\pi&=\{(x,\eta(x)):0<x<\pi\},\\
{{\mathcal B}}_\pi&={{\mathcal B}}\cap\{0<x<\pi\}.
\end{aligned}
\]

\begin{definition}\label{def_n}
A graph solution $(\lambda,v)\in\mathcal O$ is \emph{positively oriented} if its reconstructed flow satisfies the following conditions.
\begin{enumerate}[label=\textup{(N\arabic*)},leftmargin=3em]
\item $V>0$ in $D^+\cup D^-\cup{{\mathcal I}}_\pi\cup{{\mathcal S}}_\pi$;
\item Away from the bed corners, $V_x>0$ on $\{x=0\}\cap\overline D$ and $V_x<0$ on $\{x=\pi\}\cap\overline D$, including the interface and free-surface endpoints.
\item $V_y>0$ on ${{\mathcal B}}_\pi$;
\item $V_{xy}>0$ at $(0,-h)$ and $V_{xy}<0$ at $(\pi,-h)$.
\end{enumerate}
It is \emph{negatively oriented} if $-V$ has these properties; $\mathcal T_\pi$ exchanges the two.
\end{definition}

We write \(\partial_n\) for differentiation along the unit normal \(n\) of \autoref{lem_v}. The constant boundary traces give \(\psi_y=n_y\partial_n\psi\), with \(n_y>0\) on both graphs. For solutions \((\lambda,v)\) in the component of \(\mathcal O\) containing the laminar points \((\lambda,0)\) with \(s<0<{{W_{\!I}}}\), \autoref{rem_components}(a) thus gives
\begin{equation}\label{eq_signsg}
\partial_n\psi<0\ \text{on }{{\mathcal S}},\qquad \partial_n\psi>0\ \text{on }{{\mathcal I}},\qquad\text{i.e.}\qquad \psi_y<0\ \text{on }{{\mathcal S}},\quad \psi_y>0\ \text{on }{{\mathcal I}} .
\end{equation}

\begin{lemma}\label{lem_nconseq}
Assume that \({\Omega_1}<0<{\Omega_2}\). Let $(\lambda,v)\in\mathcal O$ be a positively oriented graph solution satisfying \eqref{eq_signsg}. Then its reconstructed flow has the following properties.
\begin{enumerate}[label=\textup{(\alph*)},leftmargin=2.3em]
\item $\eta'<0$ and $b'>0$ on $(0,\pi)$, so the surface falls strictly from a crest at $x=0$ to a trough at $x=\pi$, and the interface rises strictly from a trough at $x=0$ to a crest at $x=\pi$;
\item every stagnation point in $\overline{{\Omega}_\eta}\setminus{{\mathcal B}}$ lies on $\{x=j\pi\}$ for some $j\in\mathbb Z$;
\item the bed velocity $\psi_y(\cdot,-h)$ is strictly decreasing on $[0,\pi]$;
\item away from the bed corners,
\[
\begin{aligned}
\psi_{xy}&=0
&&\text{on }\{x=0\}\text{ and on }\{x=\pi\},\\
\psi_{xx}&<0
&&\text{on }\{x=0\}\cap\overline D,\\
\psi_{xx}&>0
&&\text{on }\{x=\pi\}\cap\overline D.
\end{aligned}
\]
In the upper layer, a stagnation point at $x=\pi$ is a saddle, and a stagnation point at $x=0$ is a centre if $V_x<-{{\Omega_1}}$ at that point. In the lower layer, a stagnation point at $x=0$ is a saddle, and a stagnation point at $x=\pi$ is a centre if $\psi_{xx}<{{\Omega_2}}$ at that point;
\item the open segment $\{\pi\}\times(b(\pi),\eta(\pi))$ contains exactly one stagnation point, and that point is a saddle.
\end{enumerate}
\end{lemma}
\begin{proof}
The slope identities $\eta'=V/\psi_y$ and $b'=V/\psi_y$ combine with \eqref{eq_signsg} to give profile monotonicity. On $x=0,\pi$, the Hessian $D^2\psi$ and the endpoint values of $\psi_y$ locate and classify the stagnation points.
The slopes from \autoref{lem_v}(ii),(iii), the signs in \eqref{eq_signsg}, and (N1) give
\[
\begin{aligned}
\eta'&=\frac{V}{\psi_y}\Big|_{{\mathcal S}}<0, &
b'&=\frac{V}{\psi_y}\Big|_{{\mathcal I}}>0
\qquad\text{on }(0,\pi).
\end{aligned}
\]

By reflection and periodicity, any stagnation point off the bed has an image in $\overline D\setminus{{\mathcal B}}$. At such a point $V=0$. Condition (N1) therefore places its image on $x=0$ or $x=\pi$, and the original point must lie on a line $\{x=j\pi\}$.

On the bed, (N3) gives
$\partial_x\big(\psi_y(x,-h)\big)=\psi_{xy}(x,-h)=-V_y(x,-h)<0\qquad(0<x<\pi).$

On the symmetry lines, \autoref{lem_v}(i) and (N2) give
\[
\psi_x=\psi_{xy}=0,\qquad \psi_{xx}=-V_x<0\ (x=0),\qquad
\psi_{xx}=-V_x>0\ (x=\pi),
\]
away from the bed corners. At a stagnation point $P$, \eqref{eq_intro_model} gives
\[
D^2\psi(P)=\operatorname{diag}(\psi_{xx},\psi_{yy}),\qquad
\psi_{yy}={\Omega}_i-\psi_{xx}.
\]
The signs of the Hessian entries in the four cases are
\begin{equation*}
\begin{aligned}
&{\Omega}^+,\quad x=\pi:\quad \psi_{xx}>0,\quad \psi_{yy}={{\Omega_1}}-\psi_{xx}<0,\\
&{\Omega}^+,\quad x=0,\ V_x<-{{\Omega_1}}:\quad \psi_{xx}=-V_x<0,\quad \psi_{yy}={{\Omega_1}}+V_x<0,\\
&{\Omega}^-,\quad x=0:\quad \psi_{xx}<0,\quad \psi_{yy}={{\Omega_2}}-\psi_{xx}>0,\\
&{\Omega}^-,\quad x=\pi,\ \psi_{xx}<{{\Omega_2}}:\quad \psi_{xx}>0,\quad \psi_{yy}={{\Omega_2}}-\psi_{xx}>0.
\end{aligned}
\end{equation*}
The indefinite Hessians give nondegenerate saddles. The negative- and positive-definite Hessians give a strict local maximum and minimum, respectively, both centres.

On the segment $\{\pi\}\times[b(\pi),\eta(\pi)]$, stagnation is equivalent to $\psi_y=0$ because $\psi_x=0$. By Step 2 and \eqref{eq_intro_model}, at every point of the open segment,
$\partial_y\psi_y={{\Omega_1}}-\psi_{xx}<0.$
By \eqref{eq_signsg} the endpoint signs are $\psi_y(\pi,b(\pi))>0$ and $\psi_y(\pi,\eta(\pi))<0$. Strict monotonicity gives a unique zero in the open segment. By Step 2 this point is a saddle.
\end{proof}

Under the assumptions of \autoref{lem_nconseq}, every streamline in $D^\pm$ is strictly monotone in $y$, since its tangent is parallel to $(\psi_y,V)$ with $V>0$. If the segment $\{0\}\times(b(0),\eta(0))$ contains exactly one stagnation point and it is a centre, then this centre and the saddle on $x=\pi$ from \autoref{lem_nconseq}(e) give a cat's eye in the upper layer; this holds near the folds (\autoref{lem_catseye}), but we do not prove it in general.

For $(\lambda,v)\in\mathcal O$, use $w=w_0(\lambda)+v$ to define the \emph{flattened vertical velocities}
\begin{equation}\label{eq_vflatg}
\widetilde V^\pm(x,t):=V\big(x,Y^\pm(x,t)\big)=-\Big(\Phi^\pm_x-\frac{\Phi^\pm_t\,Y^\pm_x}{T^\pm}\Big)\qquad\text{on }[0,\pi]\times[0,1],
\end{equation}
by \eqref{eq_strip_vel}, where $\Phi^\pm$ are given by \eqref{eq_strip_phi} and reconstructed by \eqref{eq_reconstruct_phi}. The fields $\Phi^\pm$, $Y^\pm$ and $T^\pm$ depend real-analytically on the complete variables with values in $C^{3,\alpha}$. Since $T^\pm>0$ on $\mathcal O$, \eqref{eq_vflatg} defines a real-analytic map in the complete variables and, by analytic reconstruction, in the perturbation variables, with values in
$C^{2,\alpha}([0,\pi]\times[0,1])^2\subset C^2([0,\pi]\times[0,1])^2.$
At laminar flows, \(\widetilde V^\pm=0\). Evenness and $2\pi$-periodicity of $\eta,b$ give $Y^\pm_x=T^\pm_x=0$ at $x=0,\pi$, and hence
\begin{equation*}
\begin{aligned}
\partial_x\widetilde V^\pm&=V_x &&\text{at }x=0,\pi,\\
\partial_t\widetilde V^-&=T^-V_y &&\text{at the bed},\\
\partial_{xt}\widetilde V^-&=T^-V_{xy} &&\text{at the bed corners}.
\end{aligned}
\end{equation*}
The map $\Pi_{Y^+}$ takes $(0,\pi)\times[0,1]$ onto $D^+\cup{{\mathcal I}}_\pi\cup{{\mathcal S}}_\pi$, $\Pi_{Y^-}$ takes $(0,\pi)\times(0,1]$ onto $D^-\cup{{\mathcal I}}_\pi$, and both take $\{0\}\times[0,1]$ and $\{\pi\}\times[0,1]$ onto the corresponding segments of the symmetry lines. Hence (N) is equivalent to
\begin{equation}\label{eq_nflat}
\begin{aligned}
&\widetilde V^+>0\ \text{on }(0,\pi)\times[0,1],\qquad \widetilde V^->0\ \text{on }(0,\pi)\times(0,1],\\
&\partial_x\widetilde V^+>0\ \text{on }\{0\}\times[0,1],\qquad \partial_x\widetilde V^+<0\ \text{on }\{\pi\}\times[0,1],\\
&\partial_x\widetilde V^->0\ \text{on }\{0\}\times(0,1],\qquad \partial_x\widetilde V^-<0\ \text{on }\{\pi\}\times(0,1],\\
&\partial_t\widetilde V^->0\ \text{on }(0,\pi)\times\{0\},\qquad \partial_{xt}\widetilde V^-(0,0)>0,\quad \partial_{xt}\widetilde V^-(\pi,0)<0 .
\end{aligned}
\end{equation}

\begin{lemma}\label{lem_open}
Let $\mathcal V\subset C^2([0,\pi]\times[0,1])^2$ consist of the pairs with
\[
\widetilde V^\pm=0\quad\text{on }\{0,\pi\}\times[0,1],
\qquad
\widetilde V^-=0\quad\text{on }[0,\pi]\times\{0\}.
\]
The strict signs in \eqref{eq_nflat} define an open subset of $\mathcal V$ in the $C^2$ topology. If a sequence of pairs in this subset converges uniformly to $(\widetilde V^+,\widetilde V^-)$, then
\[
\widetilde V^+\ge0,\qquad \widetilde V^-\ge0
\quad\text{on }[0,\pi]\times[0,1].
\]
\end{lemma}
\begin{proof}
We first factor out the boundary zeros of $\widetilde V^\pm$. This gives stability of the strict signs in \eqref{eq_nflat} up to the bed corners. We then pass to the uniform limit in the nonnegativity condition.

Away from the bed corners, divide by $x$, $\pi-x$ or $t$ as in the openness argument of \cite[Lemma~5.3]{DMSW}. The strict signs in \eqref{eq_nflat} give positive continuous factors up to the corresponding zero boundary, stable under $C^1$ perturbations preserving that boundary zero. At $(0,0)$,
\[
\widetilde V^-=xt\,\Xi_2,\qquad \Xi_2(x,t)=\int_0^1\!\!\int_0^1\partial_{xt}\widetilde V^-(\theta x,\rho t)\,d\theta\,d\rho ,
\]
where the factor depends continuously on $\widetilde V^-$ from $C^2$ to $C^0$ and satisfies
\begin{equation*}
\begin{aligned}
\Xi_2(0,0)&=\partial_{xt}\widetilde V^-(0,0)>0,\\
\partial_x\widetilde V^-(0,t)&=t\,\Xi_2(0,t),\qquad
\partial_t\widetilde V^-(x,0)=x\,\Xi_2(x,0).
\end{aligned}
\end{equation*}
The factor $\Xi_2$ stays positive in a neighbourhood of $(0,0)$ under small $C^2$ perturbations, so the factorization and the boundary identities preserve the strict signs in \eqref{eq_nflat} there. Factoring out $(\pi-x)t$ treats $(\pi,0)$ in the same way. Positive minima on the remaining compact sets and a finite covering complete the openness argument.

The strict signs in \eqref{eq_nflat} and the prescribed boundary zeros give nonnegative velocities on the closed strips. Uniform convergence preserves these inequalities.
\end{proof}

\subsection{The nodal theorem}\label{sub_05_the_nodal_theorem}

Suppose that a nontrivial graph wave satisfies \(V\ge0\) on \(\overline D\). The boundary nodal lemma below gives (N1)--(N4) of \autoref{def_n}. In the fixed strips these inequalities take the form \eqref{eq_nflat}.

\begin{lemma}\label{lem_nodalbdry}
Let ${\Omega}$, ${{\mathcal S}}$, ${{\mathcal I}}$, ${\Omega}^\pm$ be as in \textup{\autoref{lem_transmission}}, with $L=2\pi$, except that ${{\mathcal S}}$ and ${{\mathcal I}}$ are only required to be of class $C^2$, and let ${\Omega}$ and ${{\mathcal I}}$ be symmetric about $\{x=0\}$. Put
\[
\begin{aligned}
&D^\pm={\Omega}^\pm\cap\{0<x<\pi\},\quad {{\mathcal S}}_\pi={{\mathcal S}}\cap\{0<x<\pi\},\\
&{{\mathcal I}}_\pi={{\mathcal I}}\cap\{0<x<\pi\},\quad {{\mathcal B}}_\pi={{\mathcal B}}\cap\{0<x<\pi\}\ \text{and}\ \overline D=\overline{\Omega}\cap\{0\le x\le\pi\}.
\end{aligned}
\] Assume:
\begin{enumerate}[label=\textup{(G\arabic*)},leftmargin=3em]
\item $D^+$ and $D^-$ are connected; ${{\mathcal S}}\cap\{0\le x\le\pi\}$ and ${{\mathcal I}}\cap\{0\le x\le\pi\}$ are arcs that meet $\{x=0\}$ and $\{x=\pi\}$ only at their endpoints. Denote their endpoints on $x=0$ by $Q_{{\mathcal S}}$ and $Q_{{\mathcal I}}$, respectively, and those on $x=\pi$ by $Q'_{{\mathcal S}}$ and $Q'_{{\mathcal I}}$;
\item $\overline{\Omega}\cap\{x=0\}$ is the segment from $(0,-h)$ to $Q_{{\mathcal S}}$, which contains $Q_{{\mathcal I}}$. The open segment from $(0,-h)$ to $Q_{{\mathcal I}}$ lies in ${\Omega}^-$, and the open segment from $Q_{{\mathcal I}}$ to $Q_{{\mathcal S}}$ lies in ${\Omega}^+$. The same holds at $x=\pi$ with $(\pi,-h)$, $Q'_{{\mathcal I}}$, $Q'_{{\mathcal S}}$;
\item the tangents of ${{\mathcal I}}$ at $Q_{{\mathcal I}},Q'_{{\mathcal I}}$ and of ${{\mathcal S}}$ at $Q_{{\mathcal S}},Q'_{{\mathcal S}}$ are horizontal.
\end{enumerate}
Let $n$ be the unit normal of ${{\mathcal S}}$ pointing out of ${\Omega}$ and that of ${{\mathcal I}}$ pointing into ${\Omega}^+$. Let $V:\overline D\to{\mathbb{R}}$ be continuous, with $V|_{\overline D\cap\overline{{\Omega}^\pm}}\in C^2(\overline D\cap\overline{{\Omega}^\pm})$, harmonic in $D^+$ and in $D^-$, and assume:
\begin{enumerate}[label=\textup{(V\arabic*)},leftmargin=3em]
\item $V\ge0$ on $\overline D$, $V\not\equiv0$ on $D^+$ and $V\not\equiv0$ on $D^-$;
\item $V=0$ on $\overline D\cap\big(\{x=0\}\cup\{x=\pi\}\cup{{\mathcal B}}\big)$;
\item at every point $P\in{{\mathcal I}}\cap\overline D$ with $V(P)=0$, the one-sided gradients of $V$ agree: $[\nabla V](P)=0$;
\item at every point $P\in{{\mathcal S}}_\pi$ with $V(P)=0$, $\partial_nV(P)=0$;
\item if $V_x(Q)=0$ at $Q\in\{Q_{{\mathcal I}},Q'_{{\mathcal I}}\}$, then $[V_{xy}](Q)=0$; if $V_x(Q)=0$ at $Q\in\{Q_{{\mathcal S}},Q'_{{\mathcal S}}\}$, then $V_{xy}(Q)=0$.
\end{enumerate}
Then
\[
\begin{aligned}
&V>0\ \text{in}\ D^+\cup D^-\cup{{\mathcal I}}_\pi\cup{{\mathcal S}}_\pi;\\
&V_x>0\ \text{at every point of}\ \overline D\cap\{x=0\}\ \text{other than}\ (0,-h),\\
&\qquad\text{and}\ V_x<0\ \text{at every point of}\ \overline D\cap\{x=\pi\}\ \text{other than}\ (\pi,-h);\\
&V_y>0\ \text{on}\ {{\mathcal B}}_\pi;\ \text{and}\\
&V_{xy}(0,-h)>0>V_{xy}(\pi,-h).
\end{aligned}
\]
At $Q_{{\mathcal I}}$ and $Q'_{{\mathcal I}}$, the one-sided values of $V_x$ agree by \textup{(V3)}.
\end{lemma}
\begin{proof}
We treat the interior of each half-layer by the maximum principle and the open boundary arcs by Hopf's lemma. At the corners we use Serrin's lemma together with (V5).
The interface argument is that of \cite[Lemmas~5.2--5.3]{DMSW}; for the surface and bed see \cite[Lemma~5.2]{CS}, and for the free surface see \cite[Lemma~15]{CSV}; compare also \cite[Lemma~4.8]{GWYZ}. The $C^2$ boundaries and (G1)--(G3) provide the interior-ball and orthogonal-corner geometry required by Hopf's boundary lemma \cite{GT} and~Serrin's corner lemma \cite[Lemma~2]{Serrin}.

By the strong maximum principle we have $V>0$ in $D^\pm$. If the velocity vanishes at a boundary point, Hopf's lemma yields
\begin{equation*}
\begin{aligned}
V(P)=0,\quad P\in{{\mathcal I}}_\pi
&\ \Longrightarrow\ \partial_nV^+(P)>0>\partial_nV^-(P)\\
&\ \Longrightarrow\ [\partial_nV](P)>0,\\
V(P)=0,\quad P\in{{\mathcal S}}_\pi
&\ \Longrightarrow\ \partial_nV(P)<0,
\end{aligned}
\end{equation*}
contradicting (V3) and (V4), respectively. Thus $V>0$ also on ${{\mathcal I}}_\pi\cup{{\mathcal S}}_\pi$. Hopf's lemma further gives
\begin{equation*}
\begin{aligned}
V_x&>0 &&\text{on }\{x=0\}\cap({\Omega}^+\cup{\Omega}^-),\\
V_x&<0 &&\text{on }\{x=\pi\}\cap({\Omega}^+\cup{\Omega}^-),\\
V_y&>0 &&\text{on }{{\mathcal B}}_\pi.
\end{aligned}
\end{equation*}

We show \(V_x>0\) at the left endpoints of the interface and surface. At $Q\in\{Q_{{\mathcal I}},Q_{{\mathcal S}}\}$, (V2), harmonicity and nonnegativity give
\begin{equation*}
V(Q)=V_y(Q)=V_{yy}(Q)=V_{xx}(Q)=0,\qquad V_x(Q)\ge0,
\end{equation*}
with one-sided values at $Q_{{\mathcal I}}$; there $V_x$ agrees across ${{\mathcal I}}$ by (V3). If $V_x(Q_{{\mathcal I}})=0$, Serrin's corner lemma \cite[Lemma~2]{Serrin}, applied in the directions $e_x+e_y$ above ${{\mathcal I}}$ and $e_x-e_y$ below it, yields
\begin{equation*}
\begin{aligned}
&0<\partial^2_{e_x+e_y}V^+(Q_{{\mathcal I}})=2V^+_{xy}(Q_{{\mathcal I}}),\\
&0<\partial^2_{e_x-e_y}V^-(Q_{{\mathcal I}})=-2V^-_{xy}(Q_{{\mathcal I}}),\\
&[V_{xy}](Q_{{\mathcal I}})&=V^+_{xy}(Q_{{\mathcal I}})-V^-_{xy}(Q_{{\mathcal I}})>0,
\end{aligned}
\end{equation*}
which contradicts (V5). If $V_x(Q_{{\mathcal S}})=0$, Serrin's corner lemma \cite[Lemma~2]{Serrin} inside ${\Omega}^+$ gives
\begin{equation*}
0<\partial^2_{e_x-e_y}V(Q_{{\mathcal S}})=-2V_{xy}(Q_{{\mathcal S}}),
\end{equation*}
again contradicting (V5). This proves the asserted signs at the left endpoints.
Reflecting by $x\mapsto\pi-x$ gives the right-endpoint signs.

At each bed corner, (V2) makes all first and pure second derivatives vanish. Serrin's corner lemma \cite[Lemma~2]{Serrin} therefore gives
\begin{equation*}
\begin{aligned}
0&<\partial^2_{(e_x+e_y)/\sqrt2}V(0,-h)=V_{xy}(0,-h),\\
0&<\partial^2_{(-e_x+e_y)/\sqrt2}V(\pi,-h)=-V_{xy}(\pi,-h),
\end{aligned}
\end{equation*}
which proves the remaining corner inequalities.
\end{proof}

\begin{theorem}\label{thm_nodalg}
A nontrivial graph solution $(\lambda,v)\in\mathcal O$ whose reconstructed flow has $V\ge0$ on $\overline D$ is positively oriented. If a subset of nontrivial graph solutions in $\mathcal O$ is connected in ${\mathbb{R}}\times X$ and contains a positively oriented solution, then all its solutions are positively oriented.
\end{theorem}
\begin{proof}
For the continuation argument, compare \cite[Corollary~5.4]{DMSW}.
If \(V\equiv0\) in \({\Omega}^+\), the slope identities in \autoref{lem_v} and boundary nonstagnation give \(\eta'=b'=0\). The mean conditions give \(\eta=0\), \(b=-d\), hence \(v=0\) by \autoref{prop_operator}(iii).

If \(V\equiv0\) in \({\Omega}^-\), \autoref{lem_v} and interface nonstagnation give \(b'=0\), hence \(b=-d\). By \autoref{lem_v}(i),(ii), \(V^+\) is harmonic and shares the zero value and zero normal derivative of \(V^-\) along the flat interface:
\[
{\Delta} V^+=0\quad\text{in }{\Omega}^+,\qquad
V^+=\partial_y V^+=0\quad\text{on }y=-d.
\]
Odd reflection across the flat segment yields a harmonic extension; its zero Cauchy data force all Taylor coefficients there to vanish by the Laplace equation. Unique continuation gives \(V^+\equiv0\), reducing to the preceding case. Thus a nontrivial wave has \(V\not\equiv0\) in each layer, and oddness and periodicity give \(V\not\equiv0\) in each \(D^\pm\).

With $V$ nonzero in each layer, we verify (G1)--(G3) and (V1)--(V5) for the graph half-cell in order to apply \autoref{lem_nodalbdry}.
The smooth separated graphs give (G1)--(G2); evenness and periodicity give $\eta'(a)=b'(a)=0$ for $a\in\{0,\pi\}$, hence (G3). \autoref{lem_reg} and~\autoref{lem_v}, together with the preceding paragraph, give the required regularity, harmonicity and (V1)--(V3). At $P\in{{\mathcal S}}_\pi$ with $V(P)=0$, nonstagnation and the kinematic condition give $\psi_y(P)\ne0$ and a horizontal tangent. Thus \eqref{eq_vsurf} gives (V4):
\[
n(P)=e_y,\qquad \psi_y(P)^2V_y(P)=0,\qquad
\partial_nV(P)=V_y(P)=0.
\]
At the endpoints, with $a\in\{0,\pi\}$ put $Q_{{{\mathcal I}},a}=(a,b(a))$ and $Q_{{{\mathcal S}},a}=(a,\eta(a))$. By \autoref{lem_v}(i), $V(a,y)=0$ on each layer's segment of the symmetry line, hence $V_y(a,y)=0$ there. Continuity in each closed layer gives $V=V_y=0$ at the endpoints, and nonstagnation then gives $\psi_y\ne0$. Differentiating $b'=V/\psi_y$ and \eqref{eq_vy_jump} along ${{\mathcal I}}$, we obtain
\begin{equation*}
\begin{aligned}
b''(a)&=\frac{V_x(Q_{{{\mathcal I}},a})}{\psi_y(Q_{{{\mathcal I}},a})},\\
[V_{xy}](Q_{{{\mathcal I}},a})
&=({{\Omega_1}}-{{\Omega_2}})b''(a)
=({{\Omega_1}}-{{\Omega_2}})\frac{V_x(Q_{{{\mathcal I}},a})}{\psi_y(Q_{{{\mathcal I}},a})}.
\end{aligned}
\end{equation*}
Differentiating \eqref{eq_vsurf} along ${{\mathcal S}}$ at $Q_{{{\mathcal S}},a}$, we obtain
\begin{equation*}
\begin{aligned}
\psi_y(Q_{{{\mathcal S}},a})^2V_{xy}(Q_{{{\mathcal S}},a})
&=2\psi_y(Q_{{{\mathcal S}},a})V_x(Q_{{{\mathcal S}},a})^2\\
&\quad+\bigl(g+{{\Omega_1}}\psi_y(Q_{{{\mathcal S}},a})\bigr)V_x(Q_{{{\mathcal S}},a}).
\end{aligned}
\end{equation*}
Thus $V_x=0$ gives $[V_{xy}]=0$ at either interface endpoint and $V_{xy}=0$ at either surface endpoint, proving (V5). \autoref{lem_nodalbdry} yields all of (N1)--(N4).

Let $\mathcal Z_\ast$ be the nontrivial solution set in $\mathcal O$ and $\mathcal N^+$ its positively oriented subset. By \eqref{eq_vflatg} and \autoref{lem_open}, $\mathcal N^+$ is relatively open in $\mathcal Z_\ast$. A sequence in $\mathcal N^+$ converging in ${\mathbb{R}}\times X$ to a point of $\mathcal Z_\ast$ has $C^2$-convergent flattened velocities, so its limit has $V\ge0$ and is positively oriented by the first part. Hence $\mathcal N^+$ is relatively closed.

Every connected subset of $\mathcal Z_\ast$ meeting $\mathcal N^+$ therefore lies in $\mathcal N^+$.
\end{proof}

\subsection{Nodal limits and the orientation of the local branch}\label{sub_05_nodal_limits_and_the_orientation_of_the_local_branch}

We now consider a sequence of waves approaching a laminar flow. By normalizing $V=-\psi_x$, we identify the possible return points and determine which half of the local branch has positive orientation.

\begin{lemma}\label{lem_linnodal}
Let $(\lambda,0)\in\mathcal O$ and let $N\subset\{n\ge1:D_n(\lambda)=0\}$ be finite. On the flat layers $-h<y<-d$ and $-d<y<0$, take
\[
{\widehat\Psi}(x,y)=\sum_{n\in N}c_n\varphi_n(y)\cos nx,
\qquad
{\widehat V}=-\partial_x{\widehat\Psi},
\]
where $\varphi_n$ is the profile \eqref{eq_f} for $k=n$ at $\lambda$ and the $c_n$ are real and not all zero. If ${\widehat V}\ge0$ on $[0,\pi]\times[-h,0]$, then
$D_1(\lambda)=0,\qquad \varphi_1>0\quad\text{on }(-h,0].$
\end{lemma}
\begin{proof}
We use the first sine coefficient of ${\widehat V}$ at the interface and the surface, as in \cite[proof of Lemma~4.9]{GWYZ}. With $c_1=0$ if $1\notin N$ and $\varphi_1$ given by \eqref{eq_f} at $k=1$, orthogonality gives
\[
\frac{2}{\pi}\int_0^\pi{\widehat V}(x,y)\sin x\,dx=c_1\varphi_1(y)\ge0
\qquad(-h\le y\le0).
\]
At $y=-d$ at least one of the coefficients $nc_n\varphi_n(-d)$ is nonzero since $\varphi_n(-d)=\sinh(n(h-d))>0$. Hence ${\widehat V}(\cdot,-d)$ is nonnegative and not identically zero. The integral is strictly positive there, from which we deduce $c_1>0$, $1\in N$ and $D_1(\lambda)=0$.

At $y=0$, the projection gives $\varphi_1(0)\ge0$. If $\varphi_1(0)=0$, the dispersion relation \eqref{eq_d} and $s\ne0$ give $\varphi_1'(0)=0$. Uniqueness for $\varphi_1''=\varphi_1$ on $(-d,0)$ would give $\varphi_1(-d)=0$, contradicting $\varphi_1(-d)=\sinh(h-d)>0$. Thus $\varphi_1(0)>0$, and \autoref{lem_pos}(i) gives $\varphi_1>0$ on $(-h,0]$.
\end{proof}

\begin{lemma}\label{lem_nodallimitg}
Let $(\lambda_j,v_j)$ be nontrivial positively oriented graph solutions in the set $\mathcal O_\delta$ defined by \eqref{eq_odelta}, converging in ${\mathbb{R}}\times X$ to a trivial solution $(\lambda,0)$. Then $D_1(\lambda)=0$ and $\lambda$ is a positive root, i.e.\ the mode-one profile $\varphi$ of \eqref{eq_f} at $\lambda$ is positive on $(-h,0]$.
\end{lemma}
\begin{proof}
To obtain convergence of the normalized perturbations, we estimate their component in a complement of the kernel. For related arguments at a laminar limit, see \cite[proof of Lemma~4.9]{GWYZ} and \cite[proof of Theorem~13]{CSV}.
The limiting margins are at least $\delta$, so $(\lambda,0)$ and, for large $j$, the segments from $(\lambda_j,0)$ to $(\lambda_j,v_j)$ lie in $\mathcal O$. Set
\[
\varepsilon_j=\|v_j\|_X>0,\qquad
{\widehat w}_j=\frac{v_j}{\varepsilon_j}.
\]
Define
\[
L_j:=\int_0^1\partial_v F(\lambda_j,\theta v_j)\,d\theta,\qquad L:=\partial_v F(\lambda,0).
\]
The mean-value formula applied to $F(\lambda_j,v_j)=F(\lambda_j,0)=0$, the normalization, and analyticity give
\[
L_j{\widehat w}_j=0,\qquad \|{\widehat w}_j\|_X=1,\qquad
L_j\longrightarrow L\quad\text{in }\mathcal L(X,Y).
\]
We choose a closed complement $Z$ of $\ker L$ and write ${\widehat w}_j=k_j+z_j$. Fredholmness (\autoref{thm_fred}) gives
\begin{equation*}
\begin{aligned}
\|z_j\|_X
&\le C\|Lz_j\|_Y
=C\|(L-L_j){\widehat w}_j\|_Y
\le C\|L-L_j\|_{\mathcal L(X,Y)}\longrightarrow0.
\end{aligned}
\end{equation*}
Bounded projections and finite-dimensionality of $\ker L$ give a convergent subsequence of $k_j$. On that subsequence,
${\widehat w}_j\longrightarrow{\widehat w}\in\ker L,\qquad \|{\widehat w}\|_X=1.$

It remains to obtain convergence of the velocities in \(C^2\) as stated in \eqref{eq_graph_normalized_velocity}. The same calculation gives the variation used in \autoref{prop_local_ng}. Once nonnegativity of the limit is established, we may apply \autoref{lem_linnodal}.
By \autoref{prop_modes}, the nonzero kernel limit has only nonzero Fourier modes, with frequencies in the finite set $N=\{n\ge1:D_n(\lambda)=0\}$. Let $\varphi_n$ be the profile \eqref{eq_f} at $k=n$. Under the conjugation of \autoref{lem_conj}, the physical perturbation has coefficients $c_n$, not all zero, and satisfies
\begin{equation*}
\begin{aligned}
{\widehat m}_1&={\widehat Q}=0,\qquad
{\widehat\Psi} =\sum_{n\in N}c_n\varphi_n(y)\cos nx,\qquad
{\widehat V} =-{\widehat\Psi}_x=\sum_{n\in N}nc_n\varphi_n(y)\sin nx.
\end{aligned}
\end{equation*}

We write $\widetilde V_j=\widetilde V(\lambda_j,v_j)$ in perturbation coordinates and let $\Pi_{Y_0}$ denote the pair of maps $\Pi_{Y_0^\pm}$ from the fixed strips to the flat physical layers. The map $(\lambda,v)\mapsto\widetilde V$ is analytic by \eqref{eq_vflatg} and vanishes on the laminar family, so
\begin{equation}\label{eq_graph_normalized_velocity}
\begin{aligned}
\frac{\widetilde V_j}{\varepsilon_j}
&=\int_0^1\partial_v\widetilde V
\big(\lambda_j,\theta v_j\big){\widehat w}_j\,d\theta
\longrightarrow
\partial_v\widetilde V(\lambda,0)\,{\widehat w}
={\widehat V}\circ\Pi_{Y_0}
\quad\text{in }C^2([0,\pi]\times[0,1])^2.
\end{aligned}
\end{equation}
Indeed, at the laminar flow $\Phi_x^\pm=Y_x^\pm=0$ and $\Phi_t^\pm/T^\pm=W\circ Y_0^\pm$, so \eqref{eq_vflatg} and \eqref{eq_psihat} give
\begin{equation*}
\begin{aligned}
\partial_v\widetilde V^\pm(\lambda,0)\,{\widehat w}
&=-{\widehat\Phi}_x^\pm+(W\circ Y_0^\pm){\widehat Y}_x^\pm
=-{\widehat\Psi}_x^\pm\circ\Pi_{Y_0^\pm}
={\widehat V}^\pm\circ\Pi_{Y_0^\pm}.
\end{aligned}
\end{equation*}
Since $\widetilde V_j/\varepsilon_j$ satisfies \eqref{eq_nflat}, \autoref{lem_open} gives nonnegativity of each limit ${\widehat V}^\pm\circ\Pi_{Y_0^\pm}$ on its closed strip. Thus ${\widehat V}\ge0$ on the flat $\overline D$, and \autoref{lem_linnodal} yields $D_1(\lambda)=0$ and $\varphi_1>0$ on $(-h,0]$.
\end{proof}

\begin{proposition}\label{prop_local_ng}
Let $\lambda^\ast$ satisfy \eqref{eq_h} and \eqref{eq_hsimple} and let $(\lambda(a),w(a))$ be the local branch of \autoref{cor_local}. If the mode-one profile $\varphi_1$ at $\lambda^\ast$ is positive on $(-h,0]$, there is $a_1\in(0,a_0]$ such that
\[
(\lambda(a),w(a))\text{ is }
\begin{cases}
\text{positively oriented} & \text{for }0<a<a_1,\\
\text{negatively oriented} & \text{for }-a_1<a<0.
\end{cases}
\]
Under this positivity hypothesis, if ${{W_{\!I}}}(\lambda^\ast)>0$, the positively and negatively oriented halves have $\tau<0$ and $\tau>0$, respectively, by \eqref{eq_local_interface_amplitude}.
If the mode-one profile $\varphi_1$ changes sign on $(-d,0)$, the branch has neither orientation for small $a\ne0$.
\end{proposition}
\begin{proof}
For the kernel velocity ${\widehat V}_1=\varphi_1(y)\sin x$, the local branch satisfies
$\widetilde V(a)/a\longrightarrow{\widehat V}_1\circ\Pi_{Y_0}\quad\text{in }C^2.$
By \eqref{eq_vflatg}, the map $(\lambda,v)\mapsto\widetilde V(\lambda,v)$ is analytic into $C^2$ and vanishes on the laminar line $v=0$. In particular $\partial_\lambda\widetilde V(\lambda^\ast,0)=0$. Along the local branch $v(a)=a(w_1+z(a))$, we have $v'(0)=w_1$. The chain rule gives
\begin{equation*}
\left.\frac{\mathrm d}{\mathrm da}\widetilde V(\lambda(a),v(a))\right|_{a=0}
=\partial_v\widetilde V(\lambda^\ast,0)w_1
={\widehat V}_1\circ\Pi_{Y_0}.
\end{equation*}
The last equality uses the derivative in \eqref{eq_graph_normalized_velocity} for the physical kernel perturbation ${\widehat\Psi}=\varphi_1(y)\cos x$ (\autoref{prop_modes}).

For $\varphi_1>0$ on $(-h,0]$, the limiting velocity satisfies \eqref{eq_nflat}. The interior and side signs are
\begin{equation*}
\begin{aligned}
&{\widehat V}_1(x,y)>0 \qquad(0<x<\pi,\ -h<y\le0),\\
&\partial_x{\widehat V}_1(0,y)=\varphi_1(y)>0,
\quad \partial_x{\widehat V}_1(\pi,y)=-\varphi_1(y)<0\quad(-h<y\le0).
\end{aligned}
\end{equation*}
These include $y=0$ and $y=-d$. At the bed and its corners, the normalization $\varphi_1'(-h)=1$ gives
\begin{equation*}
\begin{aligned}
&\partial_y{\widehat V}_1(x,-h)=\varphi_1'(-h)\sin x=\sin x>0 &&(0<x<\pi),\\
&\partial_{xy}{\widehat V}_1(0,-h)=1>0,\qquad
\partial_{xy}{\widehat V}_1(\pi,-h)&=-1.
\end{aligned}
\end{equation*}

The flat-layer maps only rescale the vertical coordinate, so ${\widehat V}_1\circ\Pi_{Y_0}$ satisfies \eqref{eq_nflat}. By \autoref{lem_open}, $\widetilde V(a)/a$ satisfies \eqref{eq_nflat} for $0<|a|<a_1$, which is (N) for $w(a)$ when $a>0$ and (N) for $-V$ when $a<0$.

Suppose that $\varphi_1(y_0)<0$ for some $y_0\in(-d,0)$. The normalized velocity limit gives, at two fixed strip points,
\begin{equation*}
\begin{aligned}
&\widetilde V^+\bigl(a;\pi/2,(y_0+d)/d\bigr)=a\varphi_1(y_0)+o(a),\quad \varphi_1(y_0)<0,\\
&\widetilde V^-\bigl(a;\pi/2,1/2\bigr)=a\varphi_1(y_1)+o(a),\quad \varphi_1(y_1)=\sinh\frac{h-d}2>0,
\end{aligned}
\end{equation*}
where $y_1=-h+\frac12(h-d)$. For either sign of sufficiently small nonzero $a$, these two velocities have opposite signs. Thus neither orientation holds.
\end{proof}

\subsection{Analytic continuation}\label{sub_06_the_analytic_global_bifurcation_theorem}

We apply analytic global bifurcation to the local curve of \(F(\lambda,v)=0\). Uniqueness allows us to compare the continuations in nested admissible sets. We also use it to continue an incoming arc through a point of \(v=0\) at which local bifurcation holds.

We use analytic global bifurcation in the corrected form \cite[Theorem~6 and Remarks~7--8]{CSV}, based on the analytic theory of \cite{Dancer} and \cite[Theorem~9.1.1]{BT}; see also \cite[Theorem~4.1]{GWYZ}. Throughout this subsection and the next, $s$ denotes the parameter along the global curve; the surface velocity of the laminar flows is written $s(\lambda)$.

\begin{theorem}\label{thm_bt}
Let $X,Y$ be real Banach spaces, $U\subset{\mathbb{R}}\times X$ open and $F:U\to Y$ an operator. Assume
\begin{enumerate}[label=\textup{(H\arabic*)},leftmargin=3em]
\item $F$ is real-analytic, and $(\lambda,0)\in U$, $F(\lambda,0)=0$ for all $\lambda$ in an open interval $J_{\rm tr}$;
\item at some $\lambda^\ast\in J_{\rm tr}$, $\partial_v F(\lambda^\ast,0)$ is Fredholm of index zero with one-dimensional kernel $\operatorname{span}\{v_1\}$, and
\[
\partial_\lambda\partial_v F(\lambda^\ast,0)v_1\notin\operatorname{Ran}\partial_v F(\lambda^\ast,0);
\]

\item $\partial_v F(\lambda,v)$ is Fredholm of index zero at every $(\lambda,v)\in U$ with $F(\lambda,v)=0$;
\item there are bounded closed sets $Q_j\subset U$ with $U=\bigcup_j\operatorname{int}_U Q_j$ such that $\{F=0\}\cap Q_j$ is compact for every $j$.
\end{enumerate}
Then there is a continuous curve in the admissible solution set,
\[
\mathcal R=\{\mathcal R(s)=(\Lambda(s),\mathcal W(s)):s\in{\mathbb{R}}\}\subset\{F=0\}\cap U,
\]
with the following properties.
\begin{enumerate}[label=\textup{(C\arabic*)},leftmargin=3em]
\item $\mathcal R(0)=(\lambda^\ast,0)$, and there are a neighbourhood $N$ of $(\lambda^\ast,0)$ and $\varepsilon>0$ such that the nontrivial solutions in $N$ are precisely the points $\mathcal R(s)$, $0<|s|<\varepsilon$; these lie on the local curve of \textup{\autoref{thm_cr}}, with $s$ of the same sign as the coefficient $a$ of the kernel vector $v_1$ in (H2);
\item around each of its points $\mathcal R$ admits a real-analytic reparametrization, and $s\mapsto\mathcal R(s)$ is locally injective;
\item one of the following holds: \textup{(I)} for every $j$ there is $s_j>0$ such that $\mathcal R(s)\notin Q_j$ for all $|s|>s_j$;

\textup{(II)} $\mathcal R$ is periodic: there is $T>0$ with $\mathcal R(s+T)=\mathcal R(s)$ for all $s$.
\end{enumerate}
If $\mathcal R'$ is another continuous curve of solutions in $U$ with the properties \textup{(C1)}--\textup{(C3)}, then $\mathcal R'=\mathcal R\circ\phi$ for an increasing homeomorphism $\phi$ of ${\mathbb{R}}$ with $\phi(0)=0$.
\end{theorem}
\begin{proof}
Existence and (C3) follow from \cite[Theorem~6 and Remark~8]{CSV}, correcting the alternatives of \cite[Theorem~9.1.1]{BT}; the interval form of (H1) is justified by \cite[Theorem~2.4, Remark~2.5 and Appendix~A]{DMSW}.
The continuation supplied by the analytic bifurcation theorem has locally injective analytic parametrizations \cite[Corollary~7.5.3 and proof of Theorem~9.1.1]{BT}. We record why its incoming germ determines the outgoing germ; this is the uniqueness needed below. Suppose that $A(t),B(u)$ are two such parametrizations at a junction, both with positive parameter on the common incoming half. A continuous linear functional nonconstant on the incoming arc exists by Hahn--Banach. Subtract its value at the junction and choose its sign so that its leading coefficient on that half is positive. Absorbing the positive analytic units into local analytic parameters makes this observable $t^m$ and $u^n$. We set $a=n/\gcd(m,n)$ and $b=m/\gcd(m,n)$. The observable is strictly monotone on each incoming half, so $A(z^a)=B(z^b)$ for small $z>0$. The identity theorem applied after every continuous linear functional extends the equality to $z<0$. If exactly one exponent were even, one of the curves would take the same value at opposite nonzero parameters, contradicting local injectivity. Coprimality thus makes both odd. The outgoing germs agree, in their order along the curves.

Starting from (C1), continue the increasing parameter correspondence for as long as the curves agree. At finite corresponding endpoints continuity and the germ argument extend it. If one parameter were unbounded while the other tended to a finite endpoint, its tail would converge to an admissible point. In the periodic case local injectivity rules this out; in the escape case it is ruled out on any neighbourhood of that point contained in an exhaustion set. For this last argument take the exhaustion locally absorbing, as are the fixed-margin sets of \autoref{prop_compact} used below. This proves uniqueness up to the stated increasing reparametrization.
\end{proof}

Hypothesis (H1) requires admissible laminar zeros on an open interval $J_{\rm tr}$ around $\lambda^\ast$. For the graph admissible set \eqref{eq_odelta}, the interface margin excludes a neighbourhood of $\lambda={{\Omega_1}}^2d^2$, where ${{W_{\!I}}}(\lambda)=0$. The admissible trivial line is likewise restricted in \cite[(3.14)]{CSV} to $\lambda\ne0$. In \cite[Section~4.1]{GWYZ}, $\lambda=\lambda_0+e^\tau$ parametrizes the entire admissible trivial line $\{\lambda>\lambda_0\}$.

\begin{lemma}\label{lem_btcons}
In the setting of \autoref{thm_bt}, the continuation $\mathcal R$ has the following properties.
\begin{enumerate}[label=\textup{(\alph*)},leftmargin=2.3em]
\item \emph{Restriction.} Let $U'\subset U$ be open and contain $(\lambda^\ast,0)$. Assume that $U'$ has closed bounded sets $Q'_j$ as in \textup{(H4)}, each contained in some $Q_{j'}$. Let $\mathcal R'$ be the continuation in $U'$ and $I$ the connected component of $\{s:\mathcal R(s)\in U'\}$ containing $0$. Then
$\mathcal R'=\mathcal R\circ\phi$
for an increasing homeomorphism $\phi:{\mathbb{R}}\to I$ with $\phi(0)=0$.

\item \emph{Symmetry.} Let $T$ be a linear isometry of $X$ and $\widehat T$ an isomorphism of $Y$, with
\[
\begin{aligned}
(\lambda,Tv)&\in U,\qquad
F(\lambda,Tv)=\widehat T F(\lambda,v)
&&\text{for }(\lambda,v)\in U,\\
TQ_j&=Q_j&&\text{for all }j.
\end{aligned}
\]
Here $T$ acts on the state component of pairs. If $T$ maps the local curve of \autoref{thm_cr} at $\lambda^\ast$ to itself by $a\mapsto-a$, then
$T\{\mathcal R(s):s>0\}=\{\mathcal R(s):s<0\}.$

\item \emph{Crossing.} Suppose $\mathcal R(s_0)=(\lambda^\circ,0)$ and \textup{(H2)} holds at $\lambda^\circ$. Write the local curve of \autoref{thm_cr} there as
\[
\mathcal C^\circ=\{(\lambda^\circ(a),v^\circ(a)):|a|<a_0\},
\qquad a=\ell(v),
\]
where $\ell$ is a bounded linear functional. If $\mathcal R(s)\in\mathcal C^\circ$ on either one-sided neighbourhood $(s_0-\eta,s_0)$ or $(s_0,s_0+\eta)$ for some $\eta>0$, then $\mathcal R(s)\in\mathcal C^\circ$ for all $|s-s_0|$ sufficiently small. The local kernel coefficient $s\mapsto\ell(\mathcal W(s))$ is strictly monotone there, so $\mathcal R$ enters along one half of $\mathcal C^\circ$ and leaves along the other.
\end{enumerate}
\end{lemma}
\begin{proof}
Conditions (C1) and (C2) hold after shrinking $N$ into $U'$. For (C3), write $I=(s_-,s_+)$. If $\mathcal R$ is periodic and $I={\mathbb{R}}$, the restriction is periodic. If $\mathcal R$ is periodic and $I$ is proper, a point outside $U'$ recurs in both directions, so both endpoints of $I$ are finite. At any finite endpoint $s_e$, continuity gives $\mathcal R(s_e)\in U\setminus U'$. Hence $\mathcal R(s)$ eventually leaves each closed $Q'_j\subset U'$ as $s$ tends to $s_e$ inside $I$. An infinite endpoint in the exhaustion alternative inherits exit from each $Q'_j$ because $Q'_j\subset Q_{j'}$ for some $j'$. After reparametrizing $I$ by ${\mathbb{R}}$, the restricted curve satisfies (C3), and uniqueness identifies the two continuations.

Reversing the parameter and applying $T$ gives the curve $s\mapsto(\Lambda(-s),T\mathcal W(-s))$ in $U$. The action on the local curve gives (C1), reversal and the isometry preserve (C2), and $TQ_j=Q_j$ preserves (C3). Uniqueness gives
$(\Lambda(-s),T\mathcal W(-s))=\mathcal R(\phi(s))$
for an increasing homeomorphism $\phi$ of ${\mathbb{R}}$ fixing $0$. This identifies the two oriented halves, as in \cite[Lemma~10(ii)]{CSV}.

If $T^2=\operatorname{id}$, the decreasing homeomorphism $\sigma(s)=\phi(-s)$ fixes $0$ and satisfies $\mathcal R(\sigma^2(s))=\mathcal R(s)$. The fixed-point set of $\sigma^2$ is closed and is open by local injectivity of $\mathcal R$ at each fixed point. Hence $\sigma^2=\operatorname{id}$ on ${\mathbb{R}}$. Keeping the original parameter for $s\ge0$ and using $\sigma(-s)$ for $s<0$ gives a parametrization in which $T$ acts by $s\mapsto-s$.

By (C2) there are an analytic curve $\gamma:(-1,1)\to{\mathbb{R}}\times X$ and a homeomorphism $\varrho$ of a neighbourhood of $s_0$ onto $(-1,1)$, with $\varrho(s_0)=0$, such that $\mathcal R=\gamma\circ\varrho$; we may take $\varrho$ increasing. We write $\gamma(t)=(\lambda_\gamma(t),v_\gamma(t))$ and define the analytic function $a(t):=\ell(v_\gamma(t))$. On the one-sided interval in the hypothesis, the local parameter $a(t)$ gives
$\gamma(t)=(\lambda^\circ(a(t)),v^\circ(a(t))).$
By analyticity in $t$ near $0$ the equality extends to small $|t|$. Local injectivity of $\mathcal R$ then implies that $a$ is injective near $0$. In particular $a$ is strictly monotone there, so the continuation passes from one half of the local arc to the other.
\end{proof}

\subsection{The global curve through a positive root}\label{sub_06_the_global_curve_through_a_positive_root}

We start at a positive root of \(D_1(\lambda)=0\), so that the profile \(\varphi\) defined in \eqref{eq_f} is positive on \((-h,0]\). We use the sign of \(V=-\psi_x\) on \(0<x<\pi\) to identify a possible return to \(v=0\).

Let $\lambda^\ast$ satisfy \eqref{eq_h} and \eqref{eq_hsimple}. For \autoref{thm_bt} use the graph residual $F$ on $U={{\mathcal O}}_\delta$. The interval $J_{\rm tr}$ is the component containing $\lambda^\ast$ of the set of $\lambda$ for which the laminar flow has all margins in \eqref{eq_odelta} larger than $\delta$.

Hypotheses (H1), (H2) and (H3) follow from \autoref{prop_operator}(i), \autoref{thm_trans} and \autoref{cor_fredall}, respectively. For (H4), use the closed bounded exhaustion sets $Q_j(\delta)$ of \autoref{prop_compact}.

For each $\delta>0$ with $(\lambda^\ast,0)\in U$, let $\mathcal R_{{\rm g},\delta}$ be the resulting curve, with $w_1$ normalized as in \autoref{cor_local}. For $\delta'<\delta$, the definitions give
\begin{equation*}
\begin{aligned}
{{\mathcal O}}_\delta&\subset{{\mathcal O}}_{\delta'},\\
Q_j(\delta)&\subset Q_{j'}(\delta')\quad\text{if }1/j'\le\delta-\delta'+1/j,\ j'\ge j.
\end{aligned}
\end{equation*}
By \autoref{lem_btcons}(a), $\mathcal R_{{\rm g},\delta}$ is, after reparametrization, the restriction of $\mathcal R_{{\rm g},\delta'}$ to an open interval $I_\delta\ni0$. With compatible parametrizations, take the union
\[
\mathcal R_{\rm g}=\bigcup_{\delta>0}\mathcal R_{{\rm g},\delta}
\quad\text{on }\bigcup_\delta I_\delta
\]
and reparametrize over ${\mathbb{R}}$.

The half-period translation $\mathcal T_\pi$ satisfies \autoref{lem_btcons}(b) by \autoref{rem_components} and \autoref{cor_local}. The norms and margins are translation invariant, and $\mathcal T_\pi w_0(\lambda)=w_0(\lambda)$. Since $\mathcal T_\pi^2=\operatorname{id}$, the reparametrization in the proof of \autoref{lem_btcons}(b) gives
\[
\mathcal R_{\rm g}(-s)=\mathcal T_\pi\mathcal R_{\rm g}(s)\qquad(s\in{\mathbb{R}}).
\]
Norms and margins along the continuation are evaluated on the complete states $\mathcal A_{\rm g}(\mathcal R_{\rm g}(s))$, where $\mathcal A_{\rm g}(\lambda,v)=(\lambda,w_0(\lambda)+v)$. The surface and interface signs in \eqref{eq_signsg} hold on these original graph states by \autoref{rem_components}(a).

In the graph perturbation coordinates \eqref{un_know}, write
$\mathcal R_{\rm g}(s)=(\lambda(s),v(s))\in\mathcal O.$

\begin{theorem}\label{thm_globalg}
Assume that ${{\Omega_1}}<0<{{\Omega_2}}$. Suppose that $\lambda^\ast$ satisfies \eqref{eq_h} and \eqref{eq_hsimple}, that the mode-one profile $\varphi$ in \eqref{eq_f} at $\lambda^\ast$ is positive on $(-h,0]$, and that $s(\lambda^\ast)<0<{{W_{\!I}}}(\lambda^\ast)$.

Let $s_1$ be the first $s>0$ with $v(s)=0$, and set $s_1=\infty$ if there is no such $s$. Then the flow represented by $\mathcal R_{\rm g}(s)$ is positively oriented for $0<s<s_1$. The component $\{s(\lambda)<0<{{W_{\!I}}}(\lambda)\}$ contains exactly one other positive root $\lambda^{\ast\ast}$, and $\lambda^{\ast\ast}$ satisfies \eqref{eq_h} and \eqref{eq_hsimple}. Exactly one of the following holds.
\begin{enumerate}[label=\textup{(\roman*)},leftmargin=2.5em]
\item \textup{(Closed loop, $s_1<\infty$.)}
The first return occurs at
$\mathcal R_{\rm g}(s_1)=(\lambda^{\ast\ast},0).$
The closed loop has exactly two trivial points. Removing them gives
\[
\mathcal R_{\rm g}({\mathbb{R}})\setminus\{(\lambda^\ast,0),(\lambda^{\ast\ast},0)\}
=\mathfrak B\sqcup\mathcal T_\pi\mathfrak B,
\]
where the positively oriented bridge $\mathfrak B=\mathcal R_{\rm g}((0,s_1))$
joins the two bifurcation points and its half-period translate is negatively
oriented.

\item \textup{(No return, $s_1=\infty$.)}
In each direction along $\mathcal R_{\rm g}$, with $(\lambda,w)=\mathcal A_{\rm g}(\mathcal R_{\rm g}(s))$,
either $|\lambda|+\|w\|_X$ is unbounded or
at least one margin tends to zero along a sequence $s_n\to\pm\infty$. The margins in
\eqref{eq_odelta} are
\[
\lambda,\qquad \inf T^+,\qquad \inf T^-,\qquad
\inf_{{\mathcal S}}|\nabla\psi|,\qquad \inf_{{\mathcal I}}|\nabla\psi^+|.
\]
The signed estimate \eqref{eq_signed_flux_bounds} bounds $\lambda$ throughout this continuation.
\end{enumerate}
\end{theorem}
\begin{proof}
Nodal positivity restricts the first laminar return to two roots. Reversing the continuation rules out the starting root. If there is no return, we use the exhaustion $Q_j(\delta)$.

By \autoref{lem_convex}(i), $D_1=-c\,\mathcal E$ on $\mathcal J$ with $c>0$, so the roots of $D_1$ and $\mathcal E$ in $\mathcal J$ have the same multiplicities. Since $\partial_\lambda D_1(\lambda^\ast)\ne0$, $\lambda^\ast$ is not a double root, so ${{\Omega_2}}<{{\Omega_2}}^\ast$ and \autoref{prop_pairing} gives exactly two positive roots in the component, both simple: $\lambda^\ast$ and its partner $\lambda^{\ast\ast}$. At $\lambda^{\ast\ast}$, $D_n\neq0$ for $n\ge2$ by \autoref{lem_sl}, and ${{W_{\!I}}}(\lambda^{\ast\ast})\ne0$; so \eqref{eq_h} and \eqref{eq_hsimple} hold there.
By (C1), the continuation parameter $s$ and the graph kernel coefficient $a$ have the same sign near $s=0$, so \autoref{prop_local_ng} gives positive orientation for small $s>0$. The set
$\mathcal R_{\rm g}((0,s_1))$ is connected and consists of
nontrivial solutions in $\mathcal O$, so \autoref{thm_nodalg} preserves their
positive orientation. If $s_1<\infty$, write $\mathcal R_{\rm g}(s_1)=(\lambda^\circ,0)$ and choose $\delta>0$ with $\mathcal R_{\rm g}(s_1)\in\mathcal O_\delta$. By continuity and openness,
\[
\mathcal R_{\rm g}(s)\longrightarrow(\lambda^\circ,0)
\qquad (s\uparrow s_1),\qquad
\mathcal R_{\rm g}(s)\in\mathcal O_\delta\quad\text{for }s<s_1\text{ close to }s_1.
\]
By \autoref{lem_nodallimitg}, $\lambda^\circ$ is a positive root, and \eqref{eq_signsg}
puts it in the component $\{s(\lambda)<0<{{W_{\!I}}}(\lambda)\}$.
Hence
$\lambda^\circ\in\{\lambda^\ast,\lambda^{\ast\ast}\}.$

At a finite return $\mathcal R_{\rm g}(s_1)=(\lambda^\circ,0)$, we show that
$\lambda^\circ=\lambda^{\ast\ast}.$
Let $\mathcal C^\circ$ be the local curve of \autoref{cor_local} in $(\lambda,v)$.
Near $(\lambda^\circ,0)$ the solutions are the trivial
ones and those on $\mathcal C^\circ$. By \autoref{prop_local_ng}, the continuation approaches the return on the positively oriented half $a>0$. By \autoref{lem_btcons}(c), $\mathcal R_{\rm g}$ leaves $(\lambda^\circ,0)$ along the half $a<0$.

We fix $\delta$ so small that 
$\mathcal R_{\rm g}([0,s_1])\subset\mathcal O_\delta,$
 then $\lambda^\circ$ lies in the interval $J_{\rm tr}$ of (H1), since the laminar flows at $\lambda^\ast$ and $\lambda^\circ$ have margins larger than $\delta$ and both lie in $\{s(\lambda)<0<{{W_{\!I}}}(\lambda)\}$, where the margin conditions define an interval in $\lambda$. 
We choose the parametrization of $\mathcal R_{{\rm g},\delta}$ to agree pointwise
with $\mathcal R_{\rm g}$ on $[0,s_1]$, using the restriction in
\autoref{lem_btcons}(a). The reversed continuation at $s_1$ is
$\mathcal R^\flat(s):=\mathcal R_{{\rm g},\delta}(s_1-s).$
Near $s=0$, the reversed curve $\mathcal R^\flat$ follows the positive half of the local arc for $s>0$ and the negative half for $s<0$, so it satisfies (C1) at $\lambda^\circ$. Reversal preserves (C2) and (C3). By uniqueness, $\mathcal R^\flat$ is the continuation of \autoref{thm_bt} at $\lambda^\circ$ in $\mathcal O_\delta$.

Suppose that $\lambda^\circ=\lambda^\ast$. Uniqueness gives
$\mathcal R^\flat=\mathcal R_{{\rm g},\delta}\circ\phi$ for an increasing
homeomorphism $\phi$ with $\phi(0)=0$. The first positive trivial return of
both parametrizations is at $s_1$, so $\phi(s_1)=s_1$. On $[0,s_1]$ the chosen
parameters give the pointwise identity
$\mathcal R_{\rm g}(s_1-u)=\mathcal R_{\rm g}(\phi(u)).$
The decreasing homeomorphism $\rho(s):=\phi(s_1-s)$ of $[0,s_1]$ has a fixed point $s_m\in(0,s_1)$, and 
$\mathcal R_{\rm g}(s)=\mathcal R_{\rm g}(\rho(s))$
 with $\rho(s)\neq s$ for $s\neq s_m$ near $s_m$. This contradicts the local injectivity of $\mathcal R_{\rm g}$ at $s_m$, so $\lambda^\circ=\lambda^{\ast\ast}$. In \cite[proof of Theorem~13]{CSV}, condition (C4) excludes the same configuration.

Let $\mathfrak B=\{\mathcal R_{\rm g}(s):0<s<s_1\}$ be the \emph{bridge}, and define its $2s_1$-periodic completion by
\[
\mathcal P(s)=\begin{cases}
\mathcal R_{\rm g}(s),&0\le s\le s_1,\\
\mathcal T_\pi\mathcal R_{\rm g}(2s_1-s),&s_1\le s\le2s_1,
\end{cases}
\qquad\mathcal P(s+2s_1)=\mathcal P(s).
\]
The periodic curve $\mathcal P$ consists of solutions in $\mathcal O_\delta$. Near $s=s_1$ it arrives along the half $a>0$ of the local curve at $\lambda^{\ast\ast}$ and leaves along its image under $\mathcal T_\pi$, the half $a<0$ (\autoref{cor_local}). Thus it is locally a monotone reparametrization of that analytic curve, and the same holds near $s=0\equiv2s_1$ at $\lambda^\ast$. Away from these points $\mathcal P$ reparametrizes pieces of $\mathcal R_{\rm g}$ or $\mathcal T_\pi\mathcal R_{\rm g}$. Hence $\mathcal P$ satisfies (C1), (C2), and (C3) in the form (II). By uniqueness there is an increasing homeomorphism $\chi$ fixing $0$ such that
$\mathcal P=\mathcal R_{{\rm g},\delta}\circ\chi.$

The image of $\mathcal P$ is compact in $\mathcal O_\delta$.
For $\delta'<\delta$ the interval in \autoref{lem_btcons}(a) is all of ${\mathbb{R}}$.
At a finite endpoint $\mathcal R_{{\rm g},\delta'}$ would have a limit
outside $\mathcal O_\delta$, whereas the compact image of the
restriction is contained in that set. Thus $\mathcal R_{{\rm g},\delta'}$
and $\mathcal R_{{\rm g},\delta}$ differ by an increasing reparametrization.
Taking the union over the margins, we obtain
$\mathcal R_{\rm g}({\mathbb{R}})=\mathcal P([0,2s_1]).$
The two trivial points of this loop are $(\lambda^\ast,0)$ and $(\lambda^{\ast\ast},0)$, and its nontrivial part consists of $\mathfrak B$ and $\mathcal T_\pi\mathfrak B$. The half-period relation $\mathcal R_{\rm g}(-s)=\mathcal T_\pi\mathcal R_{\rm g}(s)$ from \autoref{lem_btcons}(b) identifies $\mathcal T_\pi\mathfrak B$ with $\mathcal R_{\rm g}((-s_1,0))$. This bridge is negatively oriented and therefore disjoint from $\mathfrak B$.

\smallskip\noindent\textit{Step 4. The exhaustion in the absence of return.}\enspace
If $s_1=\infty$, then by \autoref{lem_btcons}(b), $\mathcal R_{\rm g}(s)$ is nontrivial for all $s\neq0$. A period $T$ of $\mathcal R_{{\rm g},\delta}$ would give a return $\mathcal R_{{\rm g},\delta}(T)=(\lambda^\ast,0)$, so each $\mathcal R_{{\rm g},\delta}$ satisfies (I). It suffices to treat $s>0$ because $\mathcal T_\pi$ exchanges the two halves. Let $m(\lambda,v)>0$ be the minimum of the five margins in \eqref{eq_odelta}, evaluated using $w=w_0(\lambda)+v$. It is continuous on $\mathcal O$ by \autoref{prop_compact}. When $|\lambda|+\|w\|_X$ is bounded on the positive half-curve, we show that there is a sequence with
$m(\mathcal R_{\rm g}(s_n))\to0,\qquad s_n\to\infty.$
Suppose that $|\lambda|+\|w\|_X\le j_0$ on $\mathcal A_{\rm g}(\mathcal R_{\rm g}((0,\infty)))$.
For each sufficiently small $\delta>0$ and each $j\ge j_0$, (I) gives a point on the positive half of $\mathcal R_{{\rm g},\delta}$ outside $Q_j(\delta)$. The norm bound holds there, so $m<\delta+1/j$. Letting $\delta\downarrow0$ and $j\to\infty$ gives
$\inf_{s>0}m(\mathcal R_{\rm g}(s))=0.$
Continuity and positivity of $m(\mathcal R_{\rm g}(s))$ on compact parameter intervals give the required sequence $s_n\to\infty$. Passing to a subsequence, one fixed margin tends to zero.
\end{proof}

\subsection{Signed fluxes and uniform bounds}
The following estimates are obtained directly in the periodic physical domain and apply in either representation.

Let $n$ denote the normals specified in \autoref{lem_transmission} and let $\partial_n$ denote differentiation along them. In all integrals over ${\Omega}^\pm,{{\mathcal S}},{{\mathcal I}}$ we take the quotient by $(L,0)$, where $L=2\pi/k$. We indicate this quotient by the subscript $L$.

Throughout the estimates, ${{\Omega_1}}<0<{{\Omega_2}}$. Constants depend on $g,k,h,d,{{\Omega_1}},{{\Omega_2}},\alpha$ and on the indicated quantities. The explicit bounds and the dependence of the elliptic estimates on ellipticity and coefficient bounds make the constants uniform for $(g,k,h,d,{{\Omega_1}},{{\Omega_2}})$ in compact subsets of the parameter range.

\begin{lemma}\label{lem_fluxid}
Let ${\Omega},{{\mathcal S}},{{\mathcal I}},{\Omega}^\pm$ be as in \autoref{lem_transmission}, with boundaries ${{\mathcal S}}$ and ${{\mathcal I}}$ of class $C^{3,\alpha}$. Suppose that $\psi^\pm\in C^{3,\alpha}(\overline{{\Omega}^\pm})$ solve \eqref{eq_intro_model} with $m_2\in{\mathbb{R}}$ and
\[
\partial_n\psi<0\quad\text{on }{{\mathcal S}},\qquad
\partial_n\psi>0\quad\text{on }{{\mathcal I}},\qquad
|{\Omega}^+_L|=Ld.
\]
Let $\omega_1$ and $\omega_2$ be periodic harmonic functions on ${\Omega}^+$ and ${\Omega}^-$, respectively, with traces
\[
\omega_1\big|_{{\mathcal S}}=0,\qquad
\omega_1\big|_{{\mathcal I}}=\omega_2\big|_{{\mathcal I}}=1,\qquad
\omega_2\big|_{{\mathcal B}}=0.
\]
With the capacities of the two harmonic functions defined by
\[
\operatorname{cap}_1=\int_{{\Omega}^+}|\nabla\omega_1|^2,\qquad
\operatorname{cap}_2=\int_{{\Omega}^-}|\nabla\omega_2|^2,
\]
the fluxes satisfy
\begin{gather}
m_1\operatorname{cap}_1=m_2\operatorname{cap}_2-{{\Omega_1}}\int_{{\Omega}^+}\omega_1-{{\Omega_2}}\int_{{\Omega}^-}\omega_2,\label{eq_m1id}\\
\int_{{\mathcal S}}|\nabla\psi|\,ds=m_2\operatorname{cap}_2-{{\Omega_2}}\int_{{\Omega}^-}\omega_2-{{\Omega_1}}Ld ,\label{eq_fluxid}
\end{gather}
The lower-layer capacity satisfies
\[
\operatorname{cap}_2\le\frac{L}{\operatorname{dist}({{\mathcal I}},{{\mathcal B}})}.
\]
With the two prescribed areas \eqref{eq_intro_mean_levels}, the capacities obey
$\operatorname{cap}_1\ge L/d,\qquad \operatorname{cap}_2\ge L/(h-d).$

\end{lemma}
\begin{proof}
By the Lax--Milgram theorem on the bounded quotients there exist unique periodic harmonic functions $\omega_1$ in ${\Omega}^+_L$ and $\omega_2$ in ${\Omega}^-_L$ with the prescribed traces. It also gives the periodic $H^1$ solutions $\chi^\pm$ with zero traces and Poisson sources ${{\Omega_1}}$ and ${{\Omega_2}}$, respectively. For $\omega_j$, first subtract a smooth extension of the boundary values. Global Dirichlet Schauder estimates \cite[Theorems~6.14 and 6.19]{GT} give $C^{2,\alpha}$ and then $C^{3,\alpha}$ regularity up to the $C^{3,\alpha}$ boundaries. The maximum principle gives $0<\omega_1<1$ in ${\Omega}^+$ and $0<\omega_2<1$ in ${\Omega}^-$.

The expressions below solve the same Poisson equations and have the same boundary traces as $\psi^\pm$. By uniqueness of the Dirichlet problem,
\begin{equation}\label{eq_stream_dirichlet_decomposition}
\psi^+=m_1\omega_1+\chi^+,\qquad
\psi^-=m_1+m_2(1-\omega_2)+\chi^-.
\end{equation}

The lateral boundary terms cancel by periodicity. For ${\Omega}^+$ the outward normal is $n$ on ${{\mathcal S}}$ and $-n$ on ${{\mathcal I}}$. For ${\Omega}^-$ it is $n$ on ${{\mathcal I}}$. We apply Green's identity to $(\omega_1,\chi^+)$ in ${\Omega}^+$ and to $(\omega_2,\chi^-)$ in ${\Omega}^-$. On integrating $|\nabla\omega_j|^2$ by parts and using the boundary traces, we also obtain the following capacity identities.
\begin{equation*}
\begin{aligned}
{{\Omega_1}}\int_{{\Omega}^+}\omega_1&=-\int_{{\mathcal I}}\partial_n\chi^+,&
\operatorname{cap}_1&=\int_{\partial{\Omega}^+}\omega_1\partial_\nu\omega_1=-\int_{{\mathcal I}}\partial_n\omega_1,\\
{{\Omega_2}}\int_{{\Omega}^-}\omega_2&=\int_{{\mathcal I}}\partial_n\chi^-,&
\operatorname{cap}_2&=\int_{{\mathcal I}}\partial_n\omega_2.
\end{aligned}
\end{equation*}
Substituting the Poisson flux and capacity identities into \eqref{eq_stream_dirichlet_decomposition}, we obtain the interface fluxes
\begin{equation*}
\begin{aligned}
\int_{{\mathcal I}}\partial_n\psi^+\,ds&=-m_1\operatorname{cap}_1-{{\Omega_1}}\int_{{\Omega}^+}\omega_1,\\
\int_{{\mathcal I}}\partial_n\psi^-\,ds&=-m_2\operatorname{cap}_2+{{\Omega_2}}\int_{{\Omega}^-}\omega_2,
\end{aligned}
\end{equation*}
where continuity of the gradient, $[\nabla\psi]=0$, equates the two interface fluxes and gives \eqref{eq_m1id}.

The divergence theorem in ${\Omega}^+$ gives
\[
{{\Omega_1}}|{\Omega}^+_L|=\int_{{\Omega}^+}{\Delta}\psi=\int_{{\mathcal S}}\partial_n\psi\,ds-\int_{{\mathcal I}}\partial_n\psi^+\,ds.
\]
Here $|{\Omega}^+_L|=Ld$. On ${{\mathcal S}}$, $\psi$ is constant and $\partial_n\psi<0$, so $\nabla\psi=\partial_n\psi\,n$ and $\partial_n\psi=-|\nabla\psi|$. Substituting the upper-layer interface flux and applying \eqref{eq_m1id}, we obtain \eqref{eq_fluxid}.

For the capacity upper bound take $v=\min\{1,(y+h)/\varrho\}$, where $\varrho=\operatorname{dist}({{\mathcal I}},{{\mathcal B}})$. This Lipschitz periodic function is $0$ on ${{\mathcal B}}$ and $1$ on ${{\mathcal I}}$, since every interface point has $y+h\ge\varrho$. On the strip $\{-h<y<-h+\varrho\}$ of area $L\varrho$ per period, $|\nabla v|\le1/\varrho$; elsewhere $\nabla v=0$. Thus
\begin{equation*}
\begin{aligned}
\operatorname{cap}_2
&\le\int_{{\Omega}^-_L}|\nabla v|^2 
\le\frac{L\varrho}{\varrho^2}
 =L/\varrho.
\end{aligned}
\end{equation*}
For each layer apply the divergence theorem to $\omega_j e_y$. Only the interface contributes. The integral of the vertical component of its oriented normal is $\pm L$: for a regular lift the integrand times arclength is the signed horizontal differential, and its displacement in one period is $L$. Thus $|\int_{\Omega^j_L}\partial_y\omega_j|=L$, with $\Omega^j=\Omega^+$ or $\Omega^-$ as appropriate. Cauchy--Schwarz gives $\operatorname{cap}_j\ge L^2/|\Omega^j_L|$. The area constraints give the asserted bounds, also for overhanging profiles.
\end{proof}

\begin{corollary}\label{cor_signedflux}
Under the signs of \autoref{lem_fluxid}, the two prescribed areas and $\Omega_1<0<\Omega_2$, set $a=-\Omega_1$, $H_-=h-d$, and $A_0=(3aLd/(4\sqrt{2g}))^{2/3}$. Let $A=\max_{\mathcal S}y-\min_{\mathcal S}y$. Then
\begin{gather}
 |m_1|<ad^2,\quad -adH_-<m_2<\Omega_2H_-^2,\qquad
 0<\sqrt\lambda<ad+\tfrac12\Omega_2H_-,\label{eq_signed_flux_bounds}\\
 A<A_0,\qquad -2gh<Q<2gA_0+(ad)^2.\label{eq_height_bound_new}
\end{gather}
The bound for $\lambda$ uses the prescribed lower flux \eqref{eq_lampar}.
\end{corollary}
\begin{proof}
We set $I=\int_{\mathcal I}\partial_n\psi\,ds$ and $S=\int_{\mathcal S}|\nabla\psi|\,ds$. Both are strictly positive. Integration of the upper Poisson equation gives $I+S=aLd$. The two interface identities in the proof of \autoref{lem_fluxid} give
\[
 m_1\operatorname{cap}_1=a\int_{\Omega^+}\omega_1-I,\qquad
 m_2\operatorname{cap}_2=\Omega_2\int_{\Omega^-}\omega_2-I.
\]
Use $0<\omega_j<1$, $I<aLd$ and the two capacity lower bounds, then substitute $m_2=H_-(\sqrt\lambda-ad)+\Omega_2H_-^2/2$.

Let $y_-,y_+$ be the surface extrema. The strip between the bed and $y_-$ lies in the fluid, so its area is at most $Lh$ and $y_-\le0$. Bernoulli gives $Q\ge2gy_+$. On each of the two per-period surface arcs between extrema use $ds\ge|dy|$; the fundamental theorem applied to a primitive of $\sqrt{2g(y_+-y)}$ gives
\[
 aLd>S\ge 2\int_{y_-}^{y_+}\sqrt{2g(y_+-y)}\,dy
 =\tfrac43\sqrt{2g}\,A^{3/2}.
\]
For a flat surface the same inequality holds with right-hand side zero. Surface length is at least $L$, whence the minimum surface speed, attained at $y_+$, is at most $S/L<ad$. Since $y_+\le A$ and $y_+>-h$, Bernoulli now gives \eqref{eq_height_bound_new}.
\end{proof}

\subsection{A physical alternative in graph coordinates}
\begin{proposition}\label{prop_collars}
We fix $M\ge1$ and $\delta>0$. Along graph solutions in the sign component of \eqref{eq_signsg}, assume
\[
 \|\eta'\|_\infty+\|b'\|_\infty\le M,\quad
 \inf T^\pm\ge\delta,\quad \min_{{\mathcal S}\cup{\mathcal I}}|\nabla\psi|\ge\delta,
 \quad\lambda\ge\delta .
\]
Their complete states are bounded in $X$ and belong to a common fixed-margin set.
\end{proposition}
\begin{proof}
All constants below depend only on $M,\delta$ and the fixed physical data. The graph slopes and the vertical separation give Euclidean separation at least $\min(\delta/2,\delta/(2M))$: for abscissae within $\delta/(2M)$ use the vertical gap and the slope bound, and otherwise use the horizontal distance. At either moving boundary $|\psi_y|\ge c_0=\delta/\sqrt{1+M^2}$ by \eqref{eq_trace_reduction} and its surface counterpart. These estimates already give the fixed margins in \eqref{eq_odelta}.

\emph{Uniform collars.} The bounds on $Q,m_1,m_2$ in \autoref{cor_signedflux} put the fluid in a horizontal strip of bounded height. On this strip a quadratic function of $y$ with negative constant Laplacian and nonnegative trace bounds both signs of the zero-Dirichlet Poisson parts of $\psi$. The harmonic parts are bounded by their constant boundary data. Thus $\|\psi\|_\infty\le C$. Across the interface the distributional Laplacian is bounded since the normal derivatives agree. Interior $W^{2,4}$ estimates and Morrey's inequality give a uniform $C^{1,1/2}$ bound on balls of fixed radius separated from the surface and bed \cite[Theorems~9.11 and 7.26]{GT}. Near the bed subtract its constant trace and $\Omega_2(y+h)^2/2$ and reflect the harmonic remainder oddly. We choose a horizontal-graph collar of vertical depth $r$ smaller than one quarter of the separation. Its lower edge is at distance at least $r/\sqrt{1+M^2}$ from the surface and a fixed distance from the interface. In this periodic collar both velocity components are harmonic; Bernoulli bounds their surface traces and the interior estimate bounds their traces on the lower edge. The maximum principle yields $|\nabla\psi|\le C$ everywhere.

At the interface $\psi_y\ge c_0$ persists in a uniform collar by the interior estimate. At a surface point $P$, $-\psi_y\ge c_0$ on the surface. A uniform exterior cone of half-angle $\theta_0=\arctan(1/M)$ points vertically upward. In a ball of fixed radius $R$ meeting only the upper layer, with angle $\theta$ measured from this direction, the harmonic barrier
\[
 B=(r/R)^\gamma\sin(\gamma(\theta-\theta_0/2))/\sin(\gamma\theta_0/2),
 \qquad \gamma=\pi/(2\pi-\theta_0)
\]
is nonnegative in the fluid and at least one on its spherical edge. Comparison with $c_0-(c_0+C)B$ shows $-\psi_y\ge c_0/2$ in a smaller ball of uniform radius. Thus $q=x$, $p=-\psi$ at the surface, or $p=\psi-m_1$ at the interface, give an inverse $y=H(q,p)$ on a uniform half-box $p\le0$ or full box about $p=0$. Indeed one integrates $p_y\ge c_0/2$ on vertical segments of fixed length about the boundary graph and uses its Lipschitz bound. Conversely the image of each reduced box contains a uniformly sized physical neighbourhood by $|\nabla p|\le C$. In these boxes $H_p$ is bounded above and below by positive constants and $|H_q|\le C$.

\emph{First H\"older bound.} Define $E(\xi)=(1+\xi_1^2)/(2\xi_2)$ for $\xi_2>0$ and $A=\nabla E$. The inverse-coordinate equation and surface condition are
\[
 \operatorname{div} A(\nabla H)=-f(p),\qquad
 A_2(\nabla H)=gH-Q/2\quad(p=0\text{ at the surface}),
\]
where $f=\Delta p$ is piecewise constant. At the interface $H,\nabla H$ are continuous, so the divergence equation holds across it. The Hessian $DA$ is positive definite with determinant $\xi_2^{-4}$, uniformly on the compact rectangle of gradients just obtained. A tangential difference quotient $v=D_q^tH$ solves a homogeneous divergence equation with the segment average of $DA$ as coefficient and is uniformly bounded. At the surface its conormal data are $gv$. Reflect $v$ evenly and the coefficient by conjugation with $\mathrm{diag}(1,-1)$. Indeed, on a symmetric box $B$ compactly supported tests satisfy
\[
 \int_B \widetilde A_t\nabla\widetilde v\cdot\nabla\phi
 =2\int_{p=0}gv(q,0)\phi(q,0)\,dq
 =\int_B (0,-2gv(q,0)\mathbf1_{p>0})\cdot\nabla\phi .
\]
The vector on the right is uniformly bounded by the Lipschitz bound for $H$; no derivative of its trace is required. Interior De Giorgi--Nash estimates then bound $v$ in $C^\beta$, uniformly in $t$, for some $0<\beta\le\alpha$; at the interface the interior estimate applies directly \cite[Theorem~8.24]{GT}.

To recover the normal derivative put $F_0(p)=\int_0^p f$ and integrate the closed form defining $Z_q=A_2+F_0$, $Z_p=-A_1$. We write $J(u,v)=(-v,u)$. Then $\nabla H=A^{-1}(J\nabla Z-(0,F_0))$, where
$A^{-1}(a)=(a_1,1)(-2a_2-a_1^2)^{-1/2}.$
The radicand is bounded away from zero on segments joining attained values of $A$ by concavity. The derivative of $A^{-1}$ is therefore uniformly positive definite on those segments. The compatibility of the two components of $\nabla H$ gives a divergence equation for $D_q^t Z$ with coefficient $J^T\overline{DA^{-1}}J$, the bar denoting the segment average. Its surface Dirichlet trace is the interval average of $gH-Q/2$ and has uniformly bounded Lipschitz norm. The interior and boundary H\"older estimates \cite[Theorems~8.24, 8.27, 8.29]{GT} give $Z_q\in C^\beta$ uniformly on smaller boxes. On a strictly smaller box the bounds pass to $H_q$ and $Z_q$ as $t\to0$ by compactness and identification of the distributional derivatives. Here $F_0$ is Lipschitz across $p=0$. Since $-2A_2=(1+H_q^2)/H_p^2$ is uniformly separated from zero and $H_p>0$, the positive square root in
From $H_p=((1+H_q^2)/[-2(Z_q-F_0)])^{1/2}$ we obtain the same $C^\beta$ bound for the normal derivative.

\emph{Higher derivatives.} The weak equation for $u=H_q$ is $\operatorname{div}(DA(\nabla H)\nabla u)=0$. Its coefficients are uniformly $C^\beta$ across the interface, and its surface conormal datum is $gu\in C^\beta$. Local divergence-form Schauder estimates give $u\in C^{1,\beta}$ up to the flat surface and through the interface (for the surface datum subtract the vertical vector $(0,gu(q,0))$ in the weak equation). The remaining derivative follows on either closed side from
$(1+H_q^2)H_{pp}=2H_qH_p H_{qp}-H_p^2 H_{qq}-f(p)H_p^3.$
Thus $H$ is piecewise $C^{2,\beta}$. Differentiating the equation tangentially preserves the transmission of the trace and of the conormal flux. At the surface the differentiated conormal datum is $g$ times the differentiated trace. More explicitly, choose nested boxes $B_4\Subset B_3\Subset B_2\Subset B_1\Subset B_0$ by shrinking the vertical and lateral radii by fixed fractions, retaining the flat boundary in a half-box. The preceding H\"older and first Schauder estimates give $H\in C^{2,\beta}$ on the closed sides of $B_2$. For $j=2,3$ set $u_j=\partial_q^j H$. Its differentiated flux is $DA(\nabla H)\nabla u_j+R_j$, where $R_j$ contains only derivatives of $H$ through order $j$; it is bounded in piecewise $C^\beta$ when $H$ is bounded in piecewise $C^{j,\beta}$. This flux has continuous normal trace across the interface and normal component $g u_j$ at the surface. The divergence-form Schauder estimate bounds $u_j$ in piecewise $C^{1,\beta}$ on $B_{j+1}$, with the bound from $B_j$. Differentiating the displayed identity for $H_{pp}$ on each side recovers, in increasing order of the number of normal derivatives, all remaining derivatives of order $j+1$. Thus $H$ is bounded in piecewise $C^{4,\beta}$ on $B_4$ \cite{ADN1,ADN2,GT}. The common principal trace satisfies the complementing condition verified in Section~2. Uniform physical collars and a finite covering of one period transfer these estimates to the graphs and the inverse stream functions. Away from the collars use interior and reflected-bed estimates. Finally $C^{4,\beta}$ embeds continuously in $C^{3,\alpha}$ for every $\alpha<1$. The affine strip reconstruction and \eqref{eq_signed_flux_bounds} now bound the complete state in $X$.
\end{proof}
\begin{proof}[Proof of \autoref{thm_a}]
By \autoref{thm_globalg} and the signed estimates of \autoref{cor_signedflux}, the flux parameter stays bounded. Suppose that none of the five physical alternatives occurs on a tail of the curve. Then \autoref{prop_collars} applies on the whole tail with common constants $M,\delta$. Its norm is bounded and its margins stay positive, contradicting the escape alternative. The same argument applies in the opposite direction by half-period symmetry.
\end{proof}

\subsection{Numerical continuation away from the fold}
\begin{figure}[htbp]
\centering\includegraphics[width=\textwidth]{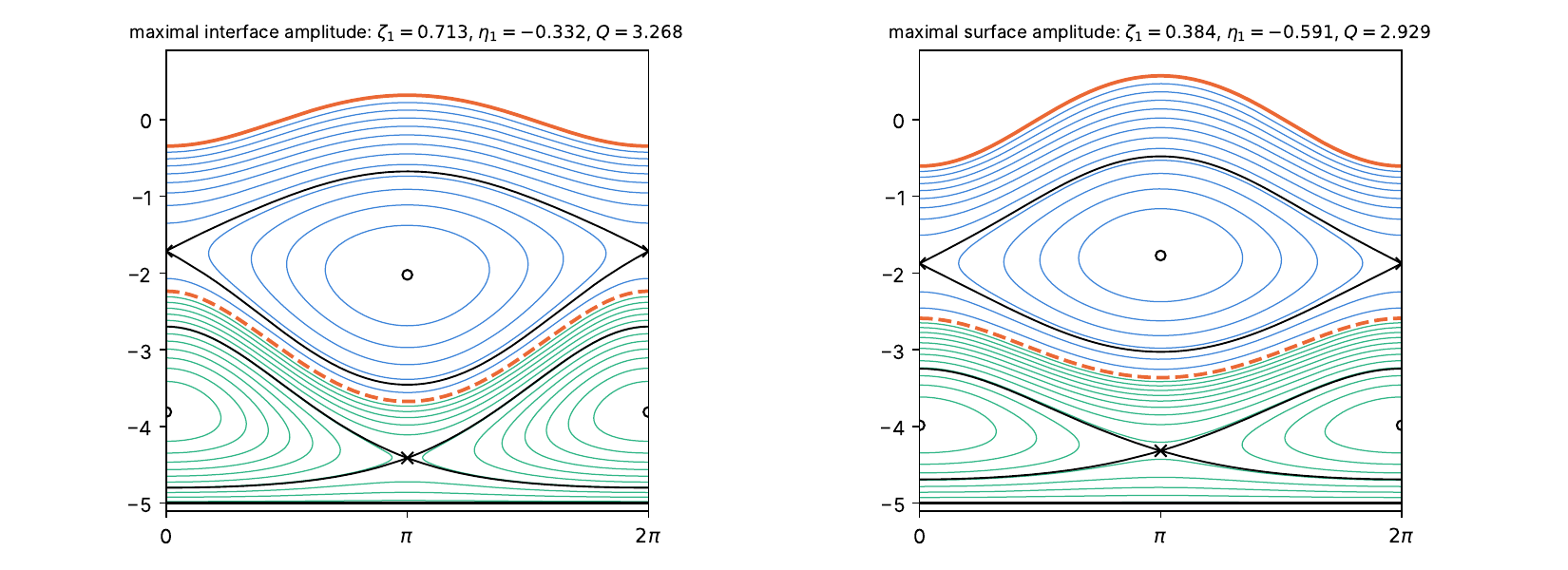}
\caption{Case A, streamlines over one period at the fold in $\zeta_1$ (left) and at the maximal surface amplitude (right), for the half-period translates of the computed solutions on the bridge. Upper layer blue, lower layer green, surface and interface orange, bed black; black curves are the separatrices; $\circ$ centres, $\times$ saddles.}\label{fig_eyes_a}
\end{figure}

We set $g=k=1$ and $(\Omega_1,\Omega_2,h,d)=(-1,1,5,3)$. The approximations in \autoref{fig_eyes_a} use cosine collocation in $x$ ($N=32$) and Chebyshev collocation in the two affine strips ($M=20$). The two starts have $\lambda_-\simeq2.6651$ (surface mode) and $\lambda_+\simeq4.0067$ (interface mode). The continuation starts with the normalized kernel at the chosen root and follows the first interface coefficient up to a turn, where pseudo-arclength is used. The collocation nodes are $x_j=\pi j/N$ and $M+1$ Chebyshev--Gauss--Lobatto points on each vertical strip. The Newton stopping rule is $\|F_{N,M}\|_\infty<10^{-10}$, with centred Jacobian differences of size $10^{-7}\max(1,|u_j|)$. The arclength norm is Euclidean on $(\eta_n,\zeta_n,m_1,m_2,Q,\zeta_1)$, counting $\zeta_1$ twice. After a failed step the step size is halved. It is multiplied by $1.25$ after at most three Newton iterations, up to $0.02$. A sign change in $\zeta_1$ triggers linear interpolation of $m_2$ to zero amplitude for comparison with the partner flux; the speed cutoffs are $\min_{\mathcal I}|\nabla\psi|<0.004$ and $\min_{\mathcal S}|\nabla\psi|^2<0.02$.

At the matched state $\zeta_1=-0.6$ on the interface-mode branch, the recorded refinement values are
\[
\begin{array}{c|ccc}
(N,M)&\lambda&Q&\eta_1\\\hline
(12,10)&3.825174000&3.584847890&0.1738249070\\
(16,12)&3.825174708&3.584848135&0.1738253657\\
(24,16)&3.825174649&3.584848103&0.1738253335\\
(32,20)&3.825174649&3.584848103&0.1738253337\\
(40,24)&3.825174649&3.584848103&0.1738253337
\end{array}
\]
The mesh comparison concerns the state at $\zeta_1=-0.6$. Along the computed branch, $|\zeta_1|$ reaches $0.713$ and the branch approaches the second laminar root. The largest observed crest-to-trough height $H$ satisfies $H/L\simeq0.19$. At the sampled states the stagnation points retain the pattern obtained for small amplitudes. These computations illustrate the continuation away from the fold. For an exact return the identity of the second laminar state follows from \autoref{thm_globalg}.
Near the fold, both halves of the small closed branch reach the second laminar state. At finite amplitude the theorem identifies that state provided the continuation returns. Whether a given branch instead approaches stagnation or reaches the boundary of the graph coordinates remains open.

\end{document}